\documentclass[10pt]{amsart}
\usepackage{egastyle}

\usepackage{amssymb,amsmath,amscd,amsthm,amsfonts,wasysym,mathrsfs,amsxtra,stmaryrd,array}

\usepackage{graphicx,array,url}
\usepackage[vcentermath]{genyoungtabtikz}
\usepackage{tikz}
\usepackage[section]{placeins}
\usepackage{afterpage}
\usepackage{verbatim}
\usepackage{tikz-cd}
\input xy
\xyoption{all}
\usepackage[normalem]{ulem}

\usepackage{algorithm}
\usepackage{algorithmic}
\usepackage{caption}

\newcommand{\OO}{\mathcal{O}}

\newcommand{\Spec}{\textnormal{Spec }}

\newcommand{\image}{\textnormal{im}\,}

\newcommand{\cokernel}{\textnormal{coker}\,}

\newcommand{\HHom}{\mathscr{H}om} 
\newcommand{\Hom}{\textnormal{Hom}}
\newcommand{\dimension}{\textnormal{dim}\,}

\newcommand{\Ext}{\textnormal{Ext}}
\newcommand{\EExt}{\mathcal{E}xt}

\newcommand{\Bc}{\mathcal{B}}

\newcommand{\Fc}{\mathcal{F}}
\newcommand{\Tc}{\mathcal{T}}

\newcommand{\Coh}{\mathrm{Coh}}

\newcommand{\arinj}{\ar@{^{(}->}}
\newcommand{\arsurj}{\ar@{->>}}
\newcommand{\areq}{\ar@{=}}

\newcommand{\ch}{\mathrm{ch}}
\newcommand{\td}{\mathrm{td}}

\newcommand{\wt}{\widetilde}

\newcommand{\bsm}{\begin{smallmatrix}}
\newcommand{\esm}{\end{smallmatrix}}

\newcommand{\projs}{{\mathbb{P}^2}}

\newcommand{\B}{\mathcal{B}}
\newcommand{\T}{\mathcal{T}}
\newcommand{\F}{\mathcal{F}}

\newcommand{\HN}{\text{HN}}

\newcommand{\cP}{\mathcal{P}}

\newcommand{\cO}{\mathcal{O}}

\newcommand{\Eff}{\mathrm{Eff}}

\makeatletter
\newtheorem*{rep@theorem}{\rep@title}
\newcommand{\newreptheorem}[2]{%
\newenvironment{rep#1}[1]{%
 \def\rep@title{#2 \ref{##1}}%
 \begin{rep@theorem}}%
 {\end{rep@theorem}}}
\makeatother

  {\begin{genthm}{Question}{true}{#1}{paragraph}}%
  {\end{genthm}}

  {\begin{genthm}{Question}{true}{#1}{equation}}%
  {\end{genthm}}

\newreptheorem{theorem}{Theorem}

\usepackage{scalerel,stackengine}
\stackMath
\newcommand\reallywidehat[1]{%
\savestack{\tmpbox}{\stretchto{%
  \scaleto{%
    \scalerel*[\widthof{\ensuremath{#1}}]{\kern-.6pt\bigwedge\kern-.6pt}%
    {\rule[-\textheight/2]{1ex}{\textheight}}
  }{\textheight}%
}{0.5ex}}%
\stackon[1pt]{#1}{\tmpbox}%
}
\begin{document}

\title[]{Stability of the tangent bundle on del Pezzo surfaces}

\author[Tristan C. Collins]{Tristan C. Collins}
\address{Department of Mathematics \\
University of Toronto\\
40 St. George St.\\
Toronto, ON \\
Canada}
\email{tristanc@math.toronto.edu}

\author[Jason Lo]{Jason Lo}
\address{Department of Mathematics \\
California State University, Northridge\\
18111 Nordhoff Street\\
Northridge CA 91330 \\
USA}
\email{jason.lo@csun.edu}

\author[Yun Shi]{Yun Shi}
\address{Department of Mathematics \\
Missouri State University\\
901 S. National Ave.\\
Springfield, MO 65897
USA}
\email{ys694e@missouristate.edu}

\author[Shing-Tung Yau]{Shing-Tung Yau}
\address{Yau Mathematical Sciences Center\\
Tsinghua University\\
Haidian District, Beijing\\
China}
\email{styau@tsinghua.edu.cn}


\begin{abstract}
In this article, we study the Bridgeland stability of the tangent bundle on del Pezzo surfaces.  We obtain a complete description of the region where the tangent bundle is (un)stable with respect to the standard Bridgeland stability conditions associated to the anticanonical polarization. 
\end{abstract}

\maketitle
\tableofcontents

\section{Introduction}

This paper is the second in a series examining the Bridgeland stability of the tangent bundle of compact K\"ahler surfaces with definite canonical bundle.  Let us briefly recall the motivation for this question; see \cite{CLSY261} for more discussion.

The Leung-Yau-Zaslow (LYZ) equation is a partial differential equation describing BPS branes in type IIB string theory. Precisely, if $E\rightarrow (X,\omega)$ is a holomorphic vector bundle over a compact K\"ahler manifold $X$ with $\dim_{\mathbb{C}}X=n$, the LYZ equation seeks a hermitian metric $h$ on $E$ such that the curvature $F_h$ of the Chern connection on $E$ satisfies the partial differential equation
\begin{equation}\label{eq: dHYMIntro}
\begin{aligned}
    {\rm Im}\left(e^{-i\hat{\theta}}(\omega \otimes I_{E} -F_h)^n\right)&=0\\
    {\rm Re}\left(e^{-i\hat{\theta}}(\omega \otimes I_{E} -F_h)^n\right)&>0
\end{aligned}
\end{equation}
where the real and imaginary parts are defined by $h$.  The system~\eqref{eq: dHYMIntro} was derived independently (for vector bundles with abelian gauge groups) by Mari\~no-Minasian-Moore-Strominger \cite{MMMS} from the action functional point-of-view, and by Leung-Yau-Zaslow \cite{LYZ} by applying a real Fourier-Mukai transform to a special Lagrangian in the mirror type IIA theory. 

Following work of Douglas \cite{DFR} and Bridgeland \cite{StabTC}, BPS branes in type IIB can also be modeled by stable objects in a Bridgeland stability condition on $D^{b}{\rm Coh}(X)$. The common starting point for these two descriptions is the Donaldson-Uhlenbeck-Yau Theorem \cite{DonaldsonASDYMC, UYHYMC} which gives an equivalence between the algebro-geometric notion of slope polystability and the  existence of a solution to the Hermitian-Yang-Mills equation
\begin{equation}\label{eq: HYM}
i\Lambda F_{h} = \mu(E)I_{E}. 
\end{equation}
Indeed, in the ``large radius limit", that is setting $\omega_t=t\omega$ and taking $t\rightarrow +\infty$ we find that~\eqref{eq: dHYMIntro} reduces to~\eqref{eq: HYM} modulo terms formally decaying in $t$. 

These two descriptions of BPS branes lead naturally to the conjecture that a holomorphic vector bundle $E\rightarrow (X,\omega)$ admits a solution of the LYZ equation~\eqref{eq: dHYMIntro} if and only if it is Bridgeland stable as an object in $D^{b}{\rm Coh}(X)$.  While this conjecture is false as stated \cite{CLSY, CoShdHYM}, the existence of solutions to the LYZ equation and Bridgeland stability are at least closely related and it is expected that some refined or modified version of stability is equivalent to the existence of solutions to~\eqref{eq: dHYMIntro}. This can be viewed as a holomorphic version of the Thomas-Yau-Joyce conjecture \cite{Th, ThY, Joyce}.

The central charge relevant to the LYZ equation is given by
\[
Z(E)= -\int_{X}e^{-\frac{i}{2\pi}\omega}ch(E).
\]
The theory of the LYZ equation, and its connections with notions of stability is relatively well-developed for line bundles \cite{Ballal, ChanJacob, GaoChen, ChuLee, ChuLeeTak, CollinsYau, CollinsYaugeo, CCL, CXY, CJY, CoShdHYM, CLSY, DatarPing, DatMetSong, Jacob, JacobSheu, JacobYau, Mete, Murakami, Murakami2, Song, Sto25, Tak1, Tak2}. In several settings it has been shown that the existence of solutions to the LYZ equation for a line bundle implies its stability with respect to standard Bridgeland stability conditions. For example, this is known for for certain surfaces of Picard rank two \cite{ARCARA20161655}, and del Pezzo surfaces of Picard rank three \cite{MY25}.  The authors studied a class of line bundles on Weierstra\ss{} elliptic K3 surfaces in \cite{CLSY} using the symmetry of the derived category. In that paper we showed that the existence of solutions to the LYZ equation for these line bundles implies Bridgeland stability, but that the converse fails in general.  Similar results for some del Pezzo surfaces were obtained in \cite{CoShdHYM}.

In contrast to the case of line bundles, very little is known for higher-rank bundles on smooth projective surfaces $X$ from either the analytic or algebraic points of view; we refer the reader to \cite{DMS, DNST, LaraSaEarp, Correa, COSD, KellerScarpa, CollinsZhang, Pingali, Takahashi, Ghosh}; see the companion paper \cite{CLSY261} for more discussion.


The starting point for our investigation is the observation that if $X$ is a surface admitting a K\"ahler-Einstein metric, then, thanks to the fact that the LYZ equation asymptotes to the Hermitian-Yang-Mills equation in the large volume limit, the LYZ equations holds asymptotically in the volume; see \cite{DMS} for precise results in this direction.  On the other hand, as one decreases the volume the quadratic curvature terms appearing in the LYZ equation \eqref{eq: dHYMIntro} may give rise to non-trivial obstructions to analytic solvability.  While this question seems difficult to attack analytically, we can ask the same question in algebro-geometric terms.  Namely

\noindent{\bf Question 1}: If $X$ is a surface with $ \pm K_{X}>0$ of $K_{X}=0$, for what values of $t>0$ is $T^{1,0}_{X}$ Bridgeland stable with respect to the K\"ahler class $tK_{X}$?

The companion paper \cite{CLSY261} investigates this question when $K_{X}<0$ or $K_{X}=0$ after replacing $T^{1,0}X$ with the cotangent bundle $\Omega_{X}$.  The present paper focuses on the Fano case.  Unlike the case of $K_{X}<0$ or $K_{X}=0$, not all Fano surfaces admit K\"ahler-Einstein metrics.  Indeed, the blow-up of $\mathbb{P}^2$ at $1$ or $2$ points is not $K$-stable and hence does not admit a K\"ahler-Einstein metric \cite{Matsushima57, Mabuchi87, Futaki83}.  Perhaps more seriously, the tangent bundle to the del Pezzo surface of degree $8$ is only slope semi-stable, and hence does not even admit a Hermitian-Yang-Mills metric.

In this paper we obtain the sharp semi-stability thresholds for the tangent bundle of an anti-canonically polarized del Pezzo surface.  Our main result can be summarized in the Theorem below. 
\begin{thm}[Theorem \ref{thm:mainthm}]
Let $X=X_n$ be the del Pezzo surface obtained by blowing up $\mathbb{P}^2$ at $n$ points in general position, where $0 \leq n \leq 8$, and let $t_0$ be given in Table \ref{tab:outermostwallintro}. Then $T_{X}$ is $\sigma_{-tK}$-semistable for $t >t_0$ and unstable for $t<t_0$. Furthermore, if $n\neq 1$, then $T_{X}$ is $\sigma_{-tK}$-stable for $t>t_0$, whereas if $n=1$, it is strictly $\sigma_{-tK}$-semistable for all $t\geq t_0$.

\begin{table}[htbp]
\small
\centering
\caption{The precise location of the outermost wall for $T_{X_n}$ along the family $\sigma_{-tK_{X_n}}$.}
\label{tab:outermostwallintro}
\begin{tabular}{cccl}
\hline
$n$ & $K^2$ & $t_0$ & Remark \\
\hline
0 & 9& 0& Theorem \ref{thm:AG54-53-1}\\ 
 1& 8& $\sqrt{\frac12}$ & Remark \ref{rem:n=2}\\ 
 2& 7& $\sqrt{\frac27}$ & Corollary \ref{cor:X2stable}\\ 
 3& 6& 0& Corollary \ref{cor:X3stable}\\ 
 4& 5& 0& Corollary \ref{cor:X4stable}\\ 
 5& 4& 0& Corollary \ref{cor:X5stable}\\ 
 6& 3& 0& Corollary \ref{cor:X6stable}\\ 
 7& 2& 0& Theorem \ref{thm:TXn-stab-78}\\ 
 8& 1& 0& Theorem  \ref{thm:TXn-stab-78}\\ 
\hline
\end{tabular}
\end{table}

\end{thm}

We also describe actual destabilizers in Section 4 for the del Pezzo surfaces $X_1$ and $X_2$. 
Strikingly, for the del Pezzo surfaces $X_n$ we note that the stability threshold $t_0$ is positive if and only if $X_n$ is $K$-unstable.

\paragraph[Organization of this paper] In Section 2, we discuss some preliminaries on Bridgeland stability conditions and the outermost walls. In Section 3, we describe some general techniques for estimating  the outermost wall for the tangent sheaf, including Matlab scripts that we wrote to assist with numerical computations. The main estimates from Section 3 are summarized in Table \ref{tab:wall-est-maximizers}. In Section 4, we analyze the walls located below the general bounds in Section 3, and obtain the precise locations of the outermost walls. The main result is Theorem \ref{thm:mainthm}. We also present the destabilizers for the outermost wall when $n=1$ and $n=2$.

\paragraph[Statement on AI use] The mathematical results were obtained prior to July 2026. AI tools were used for general proofreading. 

\paragraph[Acknowledgements] We thank Jie Liu and Zhenbo Qin for helpful discussions. T.C.C.~is supported in part by NSERC Discovery grant RGPIN-2024-03853.  J.L. was partially supported by NSF grant DMS-2100906.

\section{Preliminaries}

\paragraph Let $X$ be a smooth projective surface.  For a coherent sheaf $E$ on $X$, we will write $E^\ast$ to denote the sheaf dual $\mathcal{H}om(E,\OO_X)$.  We will also write $D^b(X)$ to denote $D^b(\Coh (X))$, the bounded derived category of coherent sheaves on $X$.  Also recall that the Hirzebruch-Riemann-Roch theorem says \cite[Example 4.1.2, Appendix 4]{Harts}
    \begin{equation}\label{eq:HRR}
    1+p_a = \tfrac{1}{12}(K_X^2+c_2(T_X)).
    \end{equation}
where $p_a$ denotes the arithmetic genus of $X$, $K_X$ the canonical divisor on $X$, and $T_X$ the tangent bundle on $X$.

\paragraph[Preliminaries on Bridgeland stability conditions]
\label{section:Bridstat}
Let $X$ be a smooth projective surface over $\mathbb C$.
We briefly recall the definition of a Bridgeland stability condition and the standard family of geometric stability conditions on $D^b(X)$.

\begin{defn}\label{defn: slicing}
A {\em slicing} $\cP$ of $D^{b}(X)$ is a collection of subcategories $\cP(\phi) \subset D^{b}(X)$ for all $\phi \in \mathbb{R}$ such that
\begin{enumerate}
\item $\cP(\phi)[1] = \cP(\phi+1)$,
\item if $\phi_1 > \phi_2$ and $A\in \cP(\phi_1)$, $B \in \cP(\phi_2)$, then ${\rm Hom}(A,B) =0$,
\item every nonzero $E\in D^{b}(X)$ admits a Harder-Narasimhan (HN) filtration by exact triangles
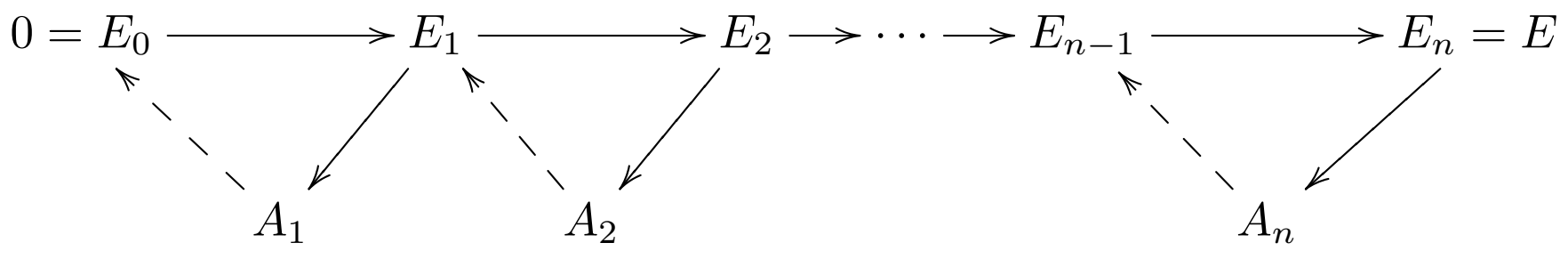
\[
\xymatrix@C=1.5em@R=2.5em{
0=E_0 \ar[rr]
&& E_1 \ar[rr] \ar[dl]
&& E_2 \ar[r] \ar[dl]
& \cdots \ar[r]
& E_{n-1} \ar[rr]
&& E_n=E \ar[dl]
\\
& A_1 \ar@{-->}[ul]
&& A_2 \ar@{-->}[ul]
&&&& A_n \ar@{-->}[ul]
&
}
\]
where $A_i \in \mathcal{P}(\phi_i)$ for $1 \leq i \leq n$ for some real numbers
\[
\phi_1 > \phi_2 > \cdots > \phi_n.
\]
\end{enumerate}
\end{defn}
\begin{defn}\label{defn: BrStab}
A Bridgeland stability condition $\sigma$ on $D^{b}(X)$  is a slicing $\cP$ together with a central charge $Z$ (group homomorphism) $K(D^b(X)) \to \mathbb{C}$ satisfying the following properties:
\begin{enumerate}
\item For any non-zero $E\in \cP(\phi)$ we have
\[
Z(E) \in \mathbb{R}_{>0} e^{\sqrt{-1}\phi},
\]
\item
\[
C := \inf \left\{ \frac{|Z(E)|}{\|\ch(E)\|} : 0 \ne E \in \cP(\phi), \phi \in \mathbb{R} \right\} >0
\]
where $\| \cdot \|$ is any norm on the finite dimensional vector space $H^{\text{even}}(X, \mathbb{R})$.
\end{enumerate}
\end{defn}

Let $\mathbb{H}$ denote the upper-half complex plane together with the negative real axis
\[
  \{ re^{i\pi \phi} : r \in \mathbb{R}_{>0}, \phi \in (0,1]\},
\]
and let $\mathbb{H}_0 = \mathbb{H}\cup \{0\}$.  Given a Bridgeland stability condition $\sigma$, objects $E$ in $\mathcal{P}(\phi)$ are said to be $\sigma$-semistable (or simply semistable) of phase $\phi$.  We also write $\phi (E)$ to denote its phase $\phi$.  We say an object $E$ is stable if it is a semistable object satisfying $\Hom (E, E)=\mathbb{C}$.  Bridgeland stability conditions also have an equivalent formulation in terms of

\begin{prop}[Bridgeland, \cite{StabTC}]\label{prop: BrStab}
	\label{pro1}
	A Bridgeland stability condition on $D^b(X)$ is an ordered pair $\sigma = (\mathcal{A}, Z)$ where $\mathcal{A}$ is the heart of a bounded t-structure on $D^b(X)$, and $Z: K(\mathcal{A})\rightarrow \mathbb{C}$ is a central charge, such that for every nonzero object $E\in \mathcal{A}$ we have 
 
 (i) $Z(E)\in \mathbb{H}$,
 
 (ii) $E$ has a finite filtration
	\begin{equation*}
	0=E_0\subset E_1\subset...\subset E_{n-1}\subset E_n=E
	\end{equation*} 
	such that the $\HN_i(E)=E_i/E_{i-1}$ are semistable objects in $\mathcal{A}$ with decreasing phase $\phi$. Furthermore, the central charge satisfies Definition \ref{defn: BrStab} (2).
\end{prop}


Using the formulation in  Proposition \ref{prop: BrStab}, semistable and stable objects can equivalently be characterized as follows: An object $0 \neq E \in \mathcal{A}$ is said to be $\sigma$-semistable (or simply semistable) if, for all nonzero proper subobject $E'$ of $E$ in $\mathcal{A}$, we have $\phi (E') \leq \phi (E)$; we say $E$ is $\sigma$-stable if we always have $\phi (E') < \phi (E)$.  An object $E\in \mathcal{A}$ that is $\sigma$-semistable but not $\sigma$-stable is said to be strictly $\sigma$-semistable.

Let $\omega\in \mathrm{NS}_\mathbb{R}(X)$ be an ample class, and $B\in \mathrm{NS}_\mathbb{R}(X)$. For any object $E\in \Coh (X)$, we define 
\begin{equation*}
\mu_{\omega, B}(E)=\frac{\omega\cdot \ch_1^B(E)}{\ch_0(E)}.
\end{equation*}
If $E$ is a torsion sheaf, we set $\mu_{\omega, B}(E)=\infty$. Then for any torsion-free sheaf $E\in \Coh(X)$, $E$ admits a Harder-Narasimhan filtration:
\begin{equation*}
	0=E_0\subset E_1\subset E_2 ...\subset E_n=E
\end{equation*}
with $\mu_{\omega, B}(E_{i+1}/E_i)$ monotonically decreasing. The factors $E_{i+1}/E_i$ are called the Harder-Narasimhan factors of $E$, and are denoted by $\HN_i(E)$.
Consider the following full subcategories of $\Coh(X)$

\begin{center}
$\mathcal{T}_{\omega, B}=\{E\in \Coh(X)| \mu_{\omega, B}(\HN_i(E))>0 \text{ for any } i\}$, 
\end{center}
\begin{center}
$\mathcal{F}_{\omega, B}=\{E\in \Coh(X)| \mu_{\omega, B}(\HN_i(E))\leq 0 \text{ for any } i\}$, 
\end{center}
together define an abelian subcategory of $D^b(X)$ by tilting $\Coh(X)$ at the torsion pair $(\T_{{\omega, B}}, \F_{\omega, B})$:
\begin{equation*}
\Coh^{\omega, B}=\langle \F_{\omega, B}[1], \T_{\omega, B}\rangle
\end{equation*}
Equivalently, we can describe $\Coh^{\omega, B}$ as
\begin{equation*}
\Coh^{\omega, B}= \{E\in D^b(X)| H^0(E)\in \mathcal{T}_{\omega, B}, H^{-1}(E)\in \mathcal{F}_{\omega, B}, H^i(E)=0 \text{ for } i\neq -1, 0\}.
\end{equation*}

 Consider the central charge taking the formula: 
\[
Z_{\omega, B}(E)=-\int e^{-i\omega}\ch^B(E).
\]
\begin{thm}[\cite{SCK3}, \cite{arcara2013minimal}]
The pair $\sigma_{\omega, B}:=(Z_{\omega, B}, \Coh^{\omega, B})$ defines a Bridgeland stability condition on $D^b(X)$.
\end{thm}

We will consider stability condition of the form $\sigma_{\omega, B}$ in this paper. 
We denote the Bridgeland slope an object $E\in D^b(X)$ by 
\[
\rho_{\omega, B}(E):=-\frac{\Re Z_{\omega, B}(E)}{\Im Z_{\omega, B}(E)}.
\]
Using the definition of phase of $E$ as in Definition \ref{defn: BrStab}, we have $\rho_{\omega, B}(E)=-\text{cot}(\pi\phi(E))$.  Also note that for nonzero objects $E, F$ in $\Coh^{\omega, B}$, we have
\[
\rho_{\omega, B}(E) \leq \rho_{\omega, B}(F) \text{\quad iff \quad} \phi (E) \leq \phi (F).
\]


\noindent
\textbf{Notation.} In all the notation introduced above, we often drop the subscript (or superscript) $B$ whenever $B=0$.

\paragraph To understand the Bridgeland stability of the tangent bundle, we begin by understanding its destabilizing sequences.  To this end, it helps to recall the definition of twisted Gieseker stability which ensures a coherent sheaf is Bridgeland semistable at large volume.  Twisted Gieseker stability is a notion originally defined by Matsuki and Wentworth \cite{matsuki1994mumford}.

\begin{defn}[twisted Gieseker stability]\cite[Definition 4.1]{LM2}\label{def:twistedGstab}
On a smooth projective surface $X$ with ample class $\omega$ and $\mathbb{Q}$-divisor $B$, we say a torsion-free sheaf $E$ is \emph{$B$-twisted Gieseker semistable} if and only if 
\begin{enumerate}
    \item[(i)] $E$ is $\mu_\omega$-semistable, and
    \item[(ii)] for every subsheaf $A \subseteq E$ such that $\mu_\omega (A) = \mu_\omega (E)$, we have
    \[
      \frac{\chi(A \otimes L)}{\ch_0(A)} \leq \frac{\chi (E \otimes L)}{\ch_0(E)}
    \]
    where $L = -B + \tfrac{K_X}{2}$ and $\chi (- \otimes L)$ is computed using the Hirzebruch-Riemann-Roch formula.
\end{enumerate}
When $B=0$, we simply say $E$ is twisted Gieseker semistable.
\end{defn}

\paragraph[Outermost wall] \label{para:outermostwall-existence} Let $X$ be a smooth projective surface, $\omega$ an ample $\mathbb{Q}$-divisor on $X$, and $E$ a torsion-free, twisted Gieseker semistable sheaf with respect to $\omega$ and with $\omega \ch_1(E)>0$.   Then there exists $t_1>0$ such that $E$ is $\sigma_{t\omega}$-semistable for all $t \geq t_1$ by \cite[Lemma 4.8]{LM2}.  By Lemma \ref{lem:goingup-v3} below, unless $E$ is $\sigma_{t\omega}$-semistable for all $t>0$,  there must be an `outermost mini-wall' $t_0 \in (0,t_1]$ such that
\begin{itemize}
    \item[(a)] $E$ is $\sigma_{t\omega}$-semistable for all $t \geq t_0$;
    \item[(b)] there exists a nonzero proper subobject $A$ of $E$ in $\Bc_\omega$ such that $\rho_{t_0\omega}(A)=\rho_{t_0\omega}(E)$ and $A$ is  $\sigma_{t_0\omega}$-semistable, while  $\rho_{t\omega} (A) > \rho_{t\omega} (E)$ for all $0<t<t_0$.  
\end{itemize}
Note that if we write $Q$ to denote  the cokernel  of the injection $A \to E$ in $\Bc_\omega$ at the outermost mini-wall, then  $\rho_{t_0\omega}(A)=\rho_{t_0\omega}(E)=\rho_{t_0\omega}(Q)$ and $Q$ is necessarily $\sigma_{t_0\omega}$-semistable.

The following lemma seems well-known and is a slight variation of \cite[Lemma 7.11]{CLSY}, but we include the statement and a proof here  for completeness.

\begin{lem}\label{lem:goingup-v3}
Let $X$ be a smooth projective surface, $\omega$ an ample divisor on $X$, and $E \in \mathcal{T}_\omega$ a coherent sheaf that is $\sigma_{t\omega}$-semistable for $t \gg 0$.  Then either $E$ is  $\sigma_{t\omega}$-semistable for all $t>0$, or there exists some $t_0>0$ such that $E$ is $\sigma_{t\omega}$-semistable for $t>t_0$, strictly $\sigma_{t_0\omega}$-semistable and $\sigma_{t\omega}$-unstable for all $0 < t < t_0$.
\end{lem}

\begin{proof}
Take any short exact sequence $0 \to A \to E \to Q \to 0$ in $\mathcal{B}_\omega$ where $A, Q \neq 0$.  We have (e.g.\ see \eqref{eq:rhoAEdiff})
\[
\rho_{t\omega}(A) - \rho_{t\omega}(E) = \frac{1}{t}\left( \frac{t^2\omega^2}{2} \left(  - \frac{1}{\mu_\omega (A)} + \frac{1}{\mu_\omega(E)}\right) + \left( \frac{\ch_2(A)}{\omega\ch_1(A)} - \frac{\ch_2(E)}{\omega\ch_1(E)}\right) \right) 
\]
and so the sign of $\rho_{t\omega}(A)-\rho_{t\omega}(E)$ is the same as that of 
\[
f_{A,E}(t) := t(\rho_{t\omega}(A)-\rho_{t\omega}(E)) = \left( \frac{t^2\omega^2}{2} \left(  - \frac{1}{\mu_\omega (A)} + \frac{1}{\mu_\omega(E)}\right) + \left( \frac{\ch_2(A)}{\omega\ch_1(A)} - \frac{\ch_2(E)}{\omega\ch_1(E)}\right) \right). 
\]
Since
\[
\frac{d}{dt} f_{A,E}(t) = \frac{\omega^2t}{\mu_\omega (A)\mu_\omega (E)} ( \mu_\omega (A) - \mu_\omega (E))
\]
which is a quadratic polynomial in $t$.

If $\mu_\omega (A) = \mu_\omega (E)$, then $f_{A,E}(t)$ is a constant polynomial in $t$ in which case the sign of $\rho_{t\omega}(A)-\rho_{t\omega}(E)$ is constant for all $t>0$.   If $\mu_\omega (A) \neq \omega (E)$, then $f_{A,E}(t)$ either has no positive root in $t$ or exactly one positive root in $t$ (and hence the same holds for $\rho_{t\omega}(A)-\rho_{t\omega}(E)$).  By the $\sigma_{t\omega}$-semistability of $E$ for $t \gg 0$, it follows that, overall, we must be in one of these three  cases:
\begin{itemize}
    \item[(i)] $\rho_{t\omega}(A)=\rho_{t\omega}(E)$ for all $t>0$;
    \item[(ii)] $\rho_{t\omega}(A) < \rho_{t\omega}(E)$ for all $t>0$;
    \item[(iii)] there exists $t'>0$ such that $\rho_{t\omega}(A)-\rho_{t\omega}(E)$ is positive for $0<t<t'$, zero for $t=t'$, and negative for $t>t'$.
\end{itemize}

Suppose now that $E$ is not $\sigma_{t\omega}$-semistable for all $t>0$.  Then  by the boundedness and local finiteness of mini-walls on the ray $\{\sigma_{t\omega}\}_{t>0}$ \cite[Theorem 1.1]{LQ}, there exists an outer-most mini-wall $t_0>0$.  The  destabilizing subobject $A$ of $E$ that gives rise to the mini-wall  $t_0$ must then be a nonzero proper subobject that is $\sigma_{t_0\omega}$-semistable  \cite[Definition 3.2]{LQ}; this definition also states that there exists some $\epsilon >0$ such that $\rho_{(t_0-\epsilon)\omega}(A) > \rho_{(t_0-\epsilon)\omega}(E)$.  Hence $A$ must be in case (iii), thus finishing the proof of the lemma.
\end{proof}

\section{General bound for the outermost wall}



A del Pezzo surface (or simply a `del Pezzo') is a smooth projective surface $X$ where the anti-canonical class $-K_X$ is ample.  (We will refrain from using the term `Fano surface' in this manuscript as it is sometimes used to refer to certain types of surfaces of general type.)

\paragraph[The tangent bundle on a del Pezzo surface] \label{para:Fanosurf-general}  For an overview on del Pezzo surfaces, the reader may refer to  \cite[$\S 0$]{mori1983fano}.  We collect here some basic facts on these  surfaces.   Recall that for any del Pezzo $X$, we have $1 \leq K_X^2 \leq 9$, and $X$ must fall into one of the following cases:
\begin{itemize}
    \item $K^2_X=9$: $X=\projs$.
    \item $K^2_X=8$: $X$ is $\mathbb{P}^1 \times \mathbb{P}^1$ or $\projs$ blown up at one point.
    \item $1 \leq K_X^2 \leq 7$: $X$ is $\projs$ blown up at $(9-K_X^2)$ points in general position.
\end{itemize}
Since $p_a (\projs)=0$ and a monoidal transformation does not alter the arithmetic genus of a smooth projective surface \cite[Corollary 3.5]{Harts}, we have $p_a(X)=0$ for every del Pezzo $X$.

In the case of $\projs$, the canonical class is $K_\projs=-3l$ where $l$ denotes the  class of a line.  From the Hirzebruch-Riemann-Roch \eqref{eq:HRR}, we have
\[
1 + p_a(\projs) = \tfrac{1}{12}(K^2_{\projs}+c_2(T_{\projs}))
\]
and so $c_2(T_{\projs})=3$, giving $\ch_2(T_{\projs})=\tfrac{1}{2}c_1(T_{\projs})^2 - c_2(T_{\projs})=\tfrac{3}{2}$.  Overall, we have
\[
\ch(T_{\projs}) = (2, 3l, \tfrac{3}{2}).
\]

For $0 \leq n \leq 8$, let us write $X_n$ to denote $\projs$ blown up at $n$ points in general position.  From above, we know $K_{X_n}^2=9-n$, and that $p_a(X_n)=p_a(\projs)=0$ for all $0 \leq n \leq 8$.  As a result, Hirzebruch-Riemann-Roch gives 
\[
 c_2(X_n) = 12-K_{X_n}^2=3+n, 
 \]
 meaning 
\[
\ch_2(T_{X_n})=\tfrac{9-n}{2}-(3+n)=\tfrac{3}{2}(1-n).
\]
We also know
\[
K_{X_n} = -3l+e_1+\cdots+e_n
\]
where (by abuse of notation) $l$ denotes the pullback of the line class from $\projs$, and $e_1, \cdots, e_n$ denote the exceptional divisors from the $n$ successive blowups.  The intersection numbers are \cite[Theorem 3.2]{pons2010acm}
\begin{gather*}
l^2=1, \text{\quad} e_i^2=-1, \text{\quad } le_i=0 \text{ for all $1 \leq i \leq n$}, \text{ and}  \\
e_i e_j =0 \text{ for all $1 \leq i \neq j\leq n$}.
\end{gather*}
Overall, we have
\[
\ch (T_{X_n}) = (2, -K_{X_n}, \tfrac{3}{2}(1-n))
\]
and 
\[
\mathrm{td}(T_{X_n}) = 1 + \tfrac{c_1(T_{X_n})}{2}+ \tfrac{1}{12}( c_1(T_{X_n})^2 + c_2(T_{X_n})) = (1, -\tfrac{1}{2}K_{X_n}, 1). 
\]

Given a del Pezzo $X$, we also know that  the tangent bundle $T_X$ is always slope stable with respect to the anti-canonical class $-K_X$ except when $X = \mathbb{P}^1 \times \mathbb{P}^1$ or $X_1$, in which case $T_X$ is only strictly slope semistable \cite[p.318]{peternell2001subsheaves}.

\noindent
\textbf{Notation.} Given a surface $X$, we will often just write $K$ to denote the canonical class $K_X$ when the context is clear.  On a del Pezzo surface of the form $X_n$, we will often write $\phi_t$ to denote the phase of an object with respect to the Bridgeland stability condition $\sigma_{-tK}$, and write $\rho_t$ to denote the Bridgeland slope $\rho_{-tK}$.

Next, we give an estimate of the rank of the destabilizing subobject at the outermost wall for the tangent bundle. First recall the following general result from \cite{CLSY261}.

\begin{prop}[\cite{CLSY261} Proposition 5.6]
\label{prop:uppboundrkdestabsubobj}
Let $X$ be a smooth projective surface.  Suppose $\omega$ is an ample class on $X$, and $E$ is a torsion-free sheaf in $\mathcal{T}_\omega$.  Then for any short exact sequence in $\Bc_\omega$ of the form
\[
0 \to A \to E \to Q \to 0
\]
where $A$ is $\sigma_\omega$-semistable and $\phi (A) = \phi (E) = \phi (Q)$, we have
\begin{equation}\label{eq:uppboundrkdestabsubobj}
\ch_0(A) < \frac{\omega c_1(E)}{\omega^2 (D+\sqrt{D^2+1})}
\end{equation}
where 
\[
D = -\frac{\Re Z_\omega (E)}{\Im Z_\omega (E)} = \frac{\ch_2 (E) - \tfrac{\omega^2}{2}\ch_0(E)}{\omega \ch_1 (E)}.
\]
\end{prop}


\begin{prop}\label{prop:AG54-46-1}
Let $X = X_n$ where $0 \leq n \leq 8$.  Suppose  $0 \to A \to T_X \to Q \to 0$ is a short exact sequence in $\mathcal{B}_{-K}$ corresponding to the outermost wall  as in \ref{para:outermostwall-existence}.  Then 
\begin{equation}\label{eq:uppboundrkdestabsubobj-v2}
\ch_0(A) < \tfrac{1}{t}( \sqrt{D^2+1}-D)
\end{equation}
where 
\[
D = \tfrac{3}{2t}\tfrac{(1-n)}{(9-n)}-t.
\]
In particular, when $n=0$ (i.e.\ $X=\mathbb{P}^2$) we have $\ch_0(A) \leq 2$ always; when $n=1$ we have $\ch_0(A) \leq 2$ whenever $t^2 \geq \tfrac{1}{3}$.
\end{prop}


\begin{proof}
Applying Proposition \ref{prop:uppboundrkdestabsubobj} to $X=X_n$ with $E=T_{X}$,  $\omega = -tK$ where $t>0$, and writing 
\[
D = \frac{\ch_2(T_{X_n})-\tfrac{\omega^2}{2}\ch_0(T_{X_n})}{\omega c_1(T_{X_n})} = \frac{\tfrac{3}{2}(1-n)-t^2(9-n)}{t(9-n)}= \frac{3}{2t}\frac{(1-n)}{(9-n)}-t,
\]

we have from  \eqref{eq:uppboundrkdestabsubobj}
\begin{align*}
\ch_0(A) &< \frac{\omega c_1(T_{X_n})}{\omega^2 (D+\sqrt{D^2+1})} \\
&=\frac{s}{D+\sqrt{D^2+1}} \text{\quad where $s := \tfrac{1}{t}$}\\
&= s (\sqrt{D^2+1}-D)  \text{ as in the proof of Corollary 5.7 in \cite{CLSY261}}.
\end{align*}
Hence $\ch_0(A) \leq 2$ would follow if we impose $s (\sqrt{D^2+1}-D) \leq 3$, i.e.\
\[
\sqrt{s^2D^2 + s^2} - sD \leq 3.
\]
Note that 
\[
sD = \tfrac{3}{2}\tfrac{(1-n)}{(9-n)}s^2-1 = Ms^2 -1
\]
if we set $M := \tfrac{3}{2}\tfrac{(1-n)}{(9-n)}$.  Then the above inequality can be rewritten as  
\begin{align*}
    \sqrt{(Ms^2-1)^2+s^2} - (Ms^2-1) \leq 3 &\Leftrightarrow  \sqrt{(Ms^2-1)^2+s^2} \leq Ms^2 + 2. 
\end{align*}
When $n=0$ or $1$, we have $M\geq 0$ and hence $Ms^2+2 > 0$, and the above inequality is equivalent to 
\begin{align*}
    &(Ms^2-1)^2 + s^2 \leq (Ms^2+2)^2  \\
    &\Leftrightarrow (1-6M)s^2 \leq 3.
\end{align*}
Since 
\[
M =\begin{cases} 
\tfrac{1}{6} &\text{ when $n=0$} \\
0 &\text{ when $n=1$} 
\end{cases}
\]
the proposition follows.
\end{proof}

Next, we focus on the case of $\projs$.  In the case of $\projs$, the tangent bundle $T_{\projs}$ is slope stable with respect to $-K_{\projs}$ and so is Bridgeland semistable with respect to $\omega = -tK_{\projs}$ for $t \gg 0$ \cite[Lemma 4.8]{LM2}.



Recall the definition of the torsion class $\T_{-K}$ in \ref{section:Bridstat}:
\[\mathcal{T}_{-K}=\{E\in \Coh(X)| \mu_{-K}(\HN_i(E))>0 \text{ for any } i\}\]

\begin{lem}\label{lem:AG54-34-1}
Suppose $0 \to A \to T_{\projs} \to Q \to 0$ is a short exact sequence in $\mathcal{B}_{-K_{\projs}}$ corresponding to the outermost wall  as in \ref{para:outermostwall-existence}.  Let 
\[
0 \to Q^{-1} \to A \overset{\beta}{\to} T_{\projs} \to Q^0 \to 0
\]
be the corresponding long exact sequence of cohomology where $Q^i := H^i(Q)$.  Then $\image \beta$ is a $\mu_{-K_{\projs}}$-semistable sheaf.
\end{lem}

\begin{proof}
Since $\image \beta$ is a subsheaf of $T_{\projs}$ and hence torsion-free, it suffices to consider the case when it is of rank 2.  Suppose $\image \beta$ is of rank 2 and is not $\mu_{-K}$-semistable.  Then it admits a 2-term HN filtration
\[
0 \neq F_0 \subsetneq F_1 = \image \beta
\]
where both HN factors $F_0, F_1/F_0$ are rank-one torsion-free sheaves with $\mu_{-K}(F_0) > \mu_{-K}(F_1/F_0) > 0$ or equivalently $\mu_l(F_0)>\mu_l(F_1/F_0)>0$,  where the second inequality holds because $\image \beta$ lies in the torsion class $\mathcal{T}_{-K}$.

Since $T_{\projs}$ is $\mu_{-K}$-stable and $-K=-3l$, it is also $\mu_l$-stable.  Then, since  $F_0$ is a rank-one subsheaf of $T_{\projs}$,  the $\mu_l$-stability of $T_{\projs}$ gives
\[
lc_1(F_0)=\mu_l(F_0) < \mu_l(T_{\projs})=\tfrac{3}{2}.
\]
Since $lc_1(F_0)$ is a positive integer, it must be 1, but then we have $1 > \mu_l(F_1/F_0)=lc_1(F_1/F_0)>0$ which is impossible.
\end{proof}

\begin{thm}\label{thm:AG54-53-1}
The tangent bundle $T_{\projs}$ on the projective plane $\projs$ is $\sigma_{-tK}$-semistable for all $t >0$ where $K = K_{\projs}$.
\end{thm}

\begin{proof}
Since $K_\projs=-3l$, we will show the equivalent statement that $T_{\projs}$ is $\sigma_{tl}$-semistable for all $t>0$.  Suppose 
\[
0 \to A \to T_{\projs} \to Q \to 0
\]
is a short exact sequence in $\mathcal{B}_l=\mathcal{B}_{-K}$ corresponding to the outermost wall  as in \ref{para:outermostwall-existence}.  Consider the corresponding long exact sequence of sheaves
\[
0 \to Q^{-1} \to A \overset{\beta}{\to} T_{\projs} \to Q^0 \to 0.
\]
By Proposition \ref{prop:AG54-46-1}, we know $\ch_0(A) \leq 2$.  We divide into two cases according to the rank of $\image \beta$:

\noindent
\textbf{Case I: $\ch_0(\image \beta)=1$.}  In this case, the $\mu_l$-stability of $T_{\projs}$ gives 
\[
0 < lc_1(\image \beta) = \mu_l (\image \beta) < \mu_l (T_{\projs}) = \tfrac{3}{2}
\]
forcing $lc_1(\image \beta)=1$.    On the other hand, we can assume $\rho_{tl}(Q) < \infty$, and so $lc_1(Q)>0$, meaning $0 < lc_1(A) < lc_1(T_{\projs})=3$, giving us $lc_1(A)=1$ or $2$.  Since we also have $lc_1(A)=lc_1(Q^{-1})+lc_1(\image \beta)$ and $lc_1(Q^{-1})\leq 0$, it follows that $lc_1(A)=1$.

Note that in the case $\ch_0(A)=2$,  since $lc_1(A)$ already attains the minimum positive value $1$ for an object in the torsion class  $\mathcal{T}_{\projs}$, any rank-one torsion-free quotient sheaf of $A$ (which is necessarily in $\mathcal{T}_{\projs}$) must have slope at least 1, meaning   the sheaf $A$ must be $\mu_l$-semistable.  That is, regardless of the rank of $A$, it is a $\mu_l$-semistable sheaf.  The Bogomolov-Gieseker (BG) inequality and the Hodge index theorem then give
\[
  \ch_2(A) \leq \frac{(c_1(A))^2}{2\ch_0(A)}  \leq \frac{(lc_1(A))^2}{2\ch_0(A)l^2} = \frac{1}{2\ch_0(A)}.
\]

With respect to the Bridgeland stability $\sigma_{tl}$ and writing $\rho_t$ to denote the corresponding slope $- \frac{\Re Z_{tl}}{\Im Z_{tl}}$, we have
\begin{align*}
    \rho_t (A) - \rho_t (T_\projs) &= \frac{\ch_2(A)-t^2\tfrac{\ch_0(A)}{2}}{t} - \frac{\tfrac{3}{2}-t^2}{3t} \\
    &=  \frac{1}{3t} \left( 3\ch_2(A)-\tfrac{3}{2} + \left( 1-\tfrac{3}{2}\ch_0(A)\right)t^2\right).
\end{align*}
Therefore, in order to achieve $\rho_t(A) \leq \rho_t(T_\projs)$, we just need
\[
3\ch_2(A)-\tfrac{3}{2} + \left( 1-\tfrac{3}{2}\ch_0(A)\right)t^2 \leq 0
\]
or equivalently
\[
t^2 \geq \frac{3\ch_2(A)-\tfrac{3}{2}}{\tfrac{3}{2}\ch_0(A)-1}
\]
which is true for any $t>0$ since the right-hand side is non-positive.

\noindent
\textbf{Case II: $\ch_0(\image \beta)=2$.} In this case, we necessarily have $\ch_0(A)=2$ and $Q^{-1}=0$, and so $A$ is a subsheaf of $T_\projs$.  

We claim that $A$ must be $\mu_l$-semistable in this case: Suppose $A$ is not $\mu_l$-semistable; then it admits a 2-term HN filtration
\[
0 \neq F_0 \subsetneq F_1 = A
\]
where $F_0, F_1/F_0$ are rank-one torsion-free sheaves with  $\mu_l(F_0)>\mu_l(F_1/F_0)>0$.  Since $F_0$ is a rank-one subsheaf of $A$ and hence $T_\projs$, the $\mu_l$-stability of $T_\projs$ gives
\[
0 < \mu_l (F_0) < \mu_l (T_\projs) = \tfrac{3}{2}
\]
and so $lc_1(F_0)=\mu_l(F_0)$ must be 1.  This is impossible, however, since
\[
1 = \mu_l(F_0) > lc_1(F_1/F_0)=\mu_l(F_1/F_0)>0
\]
from above.  Hence $A$ is $\mu_l$-semistable and the Bogomolov-Gieseker inequality and the Hodge index theorem together give
\begin{equation}\label{eq:temp1}
\ch_2(A) \leq \frac{(c_1(A))^2}{4} \leq \frac{(lc_1(A))^2}{4}.
\end{equation}
Note that $lc_1(A)=1$ or $2$ by the same argument as in Case I.

Now 
\begin{align*}
    \rho_t(A) - \rho_t(T_\projs) &= \frac{\ch_2(A) - t^2}{tlc_1(A)} - \frac{\tfrac{3}{2}-t^2}{3t} \\
&\leq \frac{\tfrac{(lc_1(A))^2}{4} - t^2}{tlc_1(A)} - \frac{\tfrac{3}{2}-t^2}{3t}
\end{align*}
and so, in order to achieve $\rho_t(A)\leq \rho_t(T_\projs)$, we just need 
\[
\frac{\tfrac{(lc_1(A))^2}{4}-t^2}{lc_1(A)} \leq \frac{\tfrac{3}{2}-t^2}{3}
\]
or equivalently
\[
t^2 \geq \frac{\tfrac{3}{4}(lc_1(A))^2-\tfrac{3}{2}lc_1(A)}{3-lc_1(A)} = \begin{cases} -\tfrac{3}{8} &\text{ if $lc_1(A)=1$} \\
0 &\text{ if $lc_1(A)=2$} \end{cases}.
\]
Hence $\rho_t(A) \leq \rho_t(T_\projs)$ for all $t >0$ in this case, too.
\end{proof}

\begin{lem}\label{lem:X1Txsslarget}
On $X=X_1$, the tangent bundle $T_{X}$ is twisted Gieseker semistable with respect to the anti-canonical class $-K$, and so is $\sigma_{-tK}$-semistable for $t \gg 0$.
\end{lem}

\begin{proof}
Since the tangent bundle $T_X$ on $X=X_1$ is $\mu_{-K}$-semistable, we just need to check condition (ii) in Definition \ref{def:twistedGstab}.  That is, we need to take a subsheaf $A$ of $T_X$ such that $\mu_{-K}(A)=\mu_{-K}(E)$ and try to show 
\[
\frac{\chi (A \otimes L)}{\ch_0(A)} \leq \frac{\chi (E \otimes L)}{\ch_0(E)}
\]
where $L = K_X/2$; note that since a twisted Gieseker semistable sheaf is required to be slope semistable, we can assume $A$ is a saturated subsheaf of $T_X$.  

Observe also that $\td (T_X) = (-1, -\tfrac{K}{2}, 1)$ and so for any coherent sheaf $F$ on $X$ we have, by Hirzebruch-Riemann-Roch,
\[
\chi (F) = \int_X \ch (F) \td (T_X) = \ch_0(F) - \tfrac{1}{2} Kc_1(F) + \ch_2(F).
\]
From this we see that we can assume $A$ has rank one, in which case $-Kc_1(A)=\mu_{-K}(A)=\mu_{-K}(T_X)=4$.  Hence the quotient sheaf $Q:= T_X/A$ is rank-one, torsion-free, and has $-Kc_1=4$.  This means $Q^\ast$ is a saturated rank-one subsheaf of $\Omega_X$ with $-Kc_1=-4$.  From \cite[Section 3]{fahlaoui1989stabilite}, we know the only saturated rank-one subsheaf of $\Omega_X$ on $X=X_1$ with degree $-Kc_1=-4$ is  $\OO_X (-2l+2e_1)$, and so $Q^\ast = \OO_X(-2l+2e_1)$ and 
\[
c_1(A) = c_1(T_X)-c_1(Q)=(-K)+(-2l+2e_1)=l+e_1.
\]

Now write $A = N \otimes I_Z$ where $N$ is a line bundle and $I_Z$ is the ideal sheaf of a 0-dimensional subsheaf $Z \subset X$.  Then 
\[
\ch (A) = (1, l+e_1, \tfrac{c_1(A)^2}{2}-l(Z))
\]
where $l(Z)$ denotes the length of $Z$ and hence
\begin{align*}
 \chi (A \otimes L) &= 1 - \tfrac{1}{2}K ( c_1(A) + c_1(L)) + \left( \tfrac{c_1(A)^2}{2} - l(Z) + c_1(A)c_1(L) + \tfrac{c_1(L)^2}{2} \right) \\
 &= 1 +0 + \left(0 -l(Z) + \tfrac{Kc_1(A)}{2} + \tfrac{K^2}{8} \right) \\
 &= -l(Z) \\
 &\leq 0.
\end{align*}
Also,
\begin{align*}
    \chi (T_X \otimes L) &= 2 - \tfrac{1}{2}K (-K + 2c_1(L)) + \left( \ch_2(T_X) + (-K)c_1(L) + c_1(L)^2 \right) \\
    &= 2+0+(0 -\tfrac{K^2}{2} + \tfrac{K^2}{4} ) \\
    &= 0.
\end{align*}
Hence 
\[
\frac{\chi (A \otimes L)}{\ch_0(A)} \leq 0 = \frac{\chi (E \otimes L)}{\ch_0(E)}
\]
and we are done.
\end{proof}

Lemma \ref{lem:X1Txsslarget}, together with the $\mu_{-K}$-stability of $T_X$ when $X=X_n$ for $n \neq 1$, gives the following corollary.

\begin{cor}\label{cor:TXsslargetalln}
Let $X=X_n$ where  $0 \leq n \leq 8$.  Then the tangent bundle $T_X$ is $\sigma_{-tK}$-semistable for $t \gg 0$.
\end{cor}


\paragraph[The effective monoid of $X_n$] Recall that an \emph{exceptional curve} on a smooth projective surface is a smooth rational curve with self-intersection number $-1$.

For $2 \leq n \leq 7$, the effective monoid 
$\Eff (X_n)$ is generated by classes of negative rational curves; for $n=8$, $\Eff (X_n)$ is generated by the negative rational curves on $X_n$ as well as $-K_{X_n}$ \cite[Proposition 2.11]{ottem2009cox}.  By adjunction, a negative smooth rational curve $C$ on $X_n$ satisfies
$C^2=-2-KC \geq -1$, and hence must be exceptional.  Since exceptional curves have been classified on del Pezzos, we obtain explicit descriptions of the effective monoid $\Eff (X_n)$.  The following formulation is taken from \cite[Theorem 2.13]{ottem2009cox}:

\begin{thm}\cite{manin1986cubic, ottem2009cox}\label{thm:fanosurf-effconegenerators}
Suppose $2 \leq n \leq 7$.  The exceptional curves on $X_n$ are strict transforms of:
\begin{itemize}
    \item[(a)] points, i.e.\ the exceptional divisors $e_i$;
    \item[(b)] lines through pairs of points, i.e.\ $l-e_i-e_j$;
    \item[(c)] conics through 5 points, i.e.\ $2l-\sum_5 e_i$;
    \item[(d)] cubics through 7 points, vanishing doubly at $e_j$, i.e.\ $3l-e_j-\sum_7 e_i$.
\end{itemize}
\end{thm}
Note that adjunction also gives that for every exceptional curve $C$ on $X_n$ we have $-KC=1$.


The next lemma will also be useful later in finding upper bounds of $\ch_2$ of sheaves:



\begin{lem}\label{lem:nKDDsmod2}
Suppose $D$ is a divisor on $X_n$ where $1 \leq n \leq 7$.  Then 
\[
-K_{X_n}D \equiv D^2 \quad\text{ mod } 2
\]
\end{lem}

\begin{proof}
In the case $2 \leq n \leq 7$, the effective monoid $\Eff (X_n)$ is generated by the exceptional curves on $X_n$ (see the proof of \cite[Proposition 2.11]{ottem2009cox}), hence so is the divisor class group $\mathrm{Div}(X_n)$ itself.   In the case $n=1$, we know the divisor class group $\mathrm{Div}(X_1)$ is generated by $l$ and $e_1$ \cite[Proposition V.3.2]{Harts}; we also have  $-Kl = 3$ and $-Ke_1=1$.  Overall, for any fixed $1 \leq n \leq 7$, the divisor class group of $X_n$ is generated by divisor classes, say $G_1, \cdots, G_m$, such that $G_i^2$ and $-KG_i$ are odd integers for all $1 \leq i\leq m$.

Take any divisor class on $X_n$ and write it as  $D = \sum_{i=1}^m a_iG_i$ for some integers $a_1, \cdots, a_m$. Then
\begin{align*}
    -KD &= \sum_{i=1}^m a_i (-KG_i) \\
    &\equiv \sum_{i=1}^m a_i \text{\quad mod }2\text{\quad since $-KG_i$ is odd for each $i$}\\
    &\equiv \sum_{i=1}^m a_i^2 \text{\quad mod }2 \\
    &\equiv \left( \sum_{i=1}^m a_iG_i \right)^2 \text{\quad mod }2 \text{\qquad since $G_i^2$ is odd for each $i$} \\
    &\equiv D^2 \text{\quad mod }2.
\end{align*}
\end{proof}


\begin{lem}\label{lem:negcurverig}
Let $C$ be a nonsingular negative curve on a smooth projective surface $X$.  Then  the only effective divisor on $X$ linearly equivalent to $C$ is $C$ itself.
\end{lem}

\begin{proof}
Since the normal bundle $N_{C/X}$ is isomorphic to $\OO_X(C)\otimes \OO_C$ \cite[Example V.1.4.1]{Harts}, tensoring the structural exact sequence
\[
0 \to \OO_X(-C) \to \OO_X \to \OO_C \to 0
\]
with $\OO_X(C)$ gives 
\[
0 \to \OO_X \to \OO_X(C) \to N_{C/X} \to 0.
\]
Consider the long exact sequence of cohomology
\[
0 \to H^0(X,\OO_X) \to H^0(X,\OO_X(C)) \to H^0(X,N_{C/X}) \to \cdots.
\]
Since the normal sheaf $N_{C/X}$ is isomorphic to $\OO_X(C) \otimes \OO_C$ which has degree $C^2 <0$ on $C$ \cite[Example V.1.4.1]{Harts}, it follows that $H^0(X,N_{C/X})\cong \Hom (\OO_X, \OO_X(C) \otimes \OO_C)=0$ and so  $H^0(X,\OO_X(C))\cong H^0(X,\OO_X)\cong \mathbb{C}$.  Hence  the complete linear system $|\OO_X(C)|$ consists of only one point, i.e.\ $C$ is the only effective divisor on $X$ that is linearly equivalent to itself.
\end{proof}

\paragraph \label{para:Xntangbdlresexccurve} Suppose $C$ is an exceptional curve on $X := X_n$ for some $1 \leq n \leq 8$.  We now compute the restriction $T_X|_C$.  From the short exact sequence of sheaves on $C$ \cite[Chap. 2, Sec. 8]{Harts}
\[
0 \to T_C \to T_X|_C \to N_{C/X} \to 0
\]
and noting $T_C \cong \omega_C^\ast \cong \OO_C(2)$ ($C$ is isomorphic to $\mathbb{P}^1)$ and $N_{C/X} \cong \OO_C (-1)$ (since $\mathrm{deg}_C(N_{C/X})=C^2=-1)$, we have $T_X|_C \cong \OO_C(2) \oplus \OO_C(-1)$ (every vector bundle on $C$ splits as a sum of line bundles).
\begin{rem}
\label{rk:X1Eff}
Let $X_1$ be $\mathbb{P}^2$ blow up at one point. The effective monoid of $X_1$ is generated by $e_1$ and $l-e_1$. 
\end{rem}

\begin{lem}\label{lem:AG54-73-1}
Suppose $1 \leq n \leq 8$, and $Q$ is a torsion quotient sheaf of $T_{X_n}$  on $X_n$ satisfying  $-Kc_1(Q)=1$.  Then $\ch_2(Q) \geq -\tfrac{1}{2}$.
\end{lem}


\begin{proof}
It is enough to consider the case when $Q$ is pure 1-dimensional quotient sheaf of $T_{X_n}$. Let $C$ denote the support of $Q$; then $C$ is an effective divisor on $X_n$ satisfying $-Kc_1(Q)=1$. 

\noindent
\underline{The case  $2\leq n \leq 8$}: Since $\Eff (X_n)$ is generated by exceptional curves and every exceptional curve has intersection number $1$ with $-K$, $C$ itself is linearly equivalent to an exceptional curve.  By Lemma \ref{lem:negcurverig}, $C$ itself is an exceptional curve.  

Now let $i$ denote the closed immersion $C \hookrightarrow X_n$, and write $Q = i_\ast \OO_C(m)$ for some integer $m$.  By Grothendieck-Riemann-Roch, we have 
\[
[X_n] \cap\left( \ch (Q )\td (T_{X_n}) \right)= i_\ast \left( [C] \cap (\ch(\OO_C(m))\td (T_C)) \right)
\]
where $\td (T_{X_n}) = \left( 1, -\tfrac{K}{2}, 1\right)$ and $\td (T_C)=(1,1)$.  Comparing the 0-cycle components on both sides above, we obtain $\ch_2 (Q)-\tfrac{KC}{2} = m+1$ and hence $\ch_2(Q)=m+\tfrac{1}{2}$.

Now the surjective morphism of sheaves $T_{X_n} \twoheadrightarrow Q$ factors as 
\[
T_{X_n} \to T_{X_n}|_C \twoheadrightarrow Q=i_\ast \OO_C(m)
\]
where $T_{X_n}|_C \cong i_\ast (\OO_C(2)\oplus \OO_C(-1))$ from \ref{para:Xntangbdlresexccurve}.  Therefore, by slope stability on $C$, we must have $m \geq -1$, and so $\ch_2(Q)\geq -\tfrac{1}{2}$.

\noindent
\underline{The case $n=1$}:  Remark \ref{rk:X1Eff} implies that an effective divisor $D$ satisfying $-KD=1$ must be linearly equivalent to $e_1$ or $l-e_1$. Hence $C$ is linearly equivalent to $e_1$ or $l-e_1$. Since every element in $|l-e_1|$ is the strict transform of a line passing through the blow up point, hence is smooth and isomorphic to $\mathbb{P}^1$. We have $C\simeq \mathbb{P}^1$, then the same computation above implies that $\ch_2(Q)\geq -\frac{1}{2}$.
\end{proof}


\begin{lem}\label{lem:AG54-112-1}
Let $X$ be a smooth projective surface and $\omega$ an ample class on $X$.  Suppose $E$ is a torsion-free sheaf in $\mathcal{T}_\omega$ and $0 \to A \to E \to Q \to 0$ is a short exact sequence in $\mathcal{B}_\omega$ corresponding to the outermost mini-wall  as in \ref{para:outermostwall-existence}.  Then $\omega c_1(Q^{-1})<0$ if $Q^{-1} \neq 0$.  That is, if $\omega c_1(Q^{-1})=0$ then $Q^{-1}=0$, in which case the short exact sequence in $\Bc_\omega$ above is also a short exact sequence in $\Coh (X)$.
\end{lem}

\begin{proof}
A priori, we know $\omega c_1(Q^{-1})\leq 0$ since $Q^{-1}$ lies in the torsion-free class $\Fc_\omega$ in $\Coh (X)$.  So it suffices to show $\omega c_1(Q^{-1})=0$ cannot occur when $Q^{-1}$ is nonzero.

Suppose $\omega c_1(Q^{-1})=0$ and $Q^{-1}\neq 0$.  Then $c_1(Q^{-1})^2 \leq 0$ by the Hodge index theorem.  Also, $Q^{-1}$ must have exactly one $\mu_\omega$-HN factor since all its HN factors have $\mu_\omega \leq 0$.  That is, $Q^{-1}$ is $\mu_\omega$-semistable, and the Bogomolov-Gieseker inequality gives $\ch_2(Q^{-1})\leq 0$.  As a result, in the long exact sequence of cohomology in $\Coh (X)$
\[
0 \to Q^{-1} \to A \overset{\beta}{\to} E \to Q^0 \to 0
\]
we have
\[
\ch_2(A) \leq \ch_2(\image \beta), \,\,\, \ch_0(A) > \ch_0(\image \beta), \,\,\, -Kc_1(A) = -Kc_1(\image \beta).
\]
These together imply that $\rho_{t\omega} (E)=\rho_{t\omega} (A) < \rho_{t\omega} (\image \beta)$ where $t>0$ is the outermost mini-wall.  Since $\image \beta$ is a $\Bc_\omega$-subobject of $E$, this contradicts the assumption that $A$ is a destabilizing subobject of $E$ at the outermost mini-wall in the family $\{\sigma_{s\omega}\}_{s>0}$.
\end{proof}

\begin{thm}\label{thm:TXn-stab-78}
For $n=7, 8$, the tangent bundle $T_{X_n}$ is $\sigma_{-tK_{X_n}}$-stable for all $t >0$.
\end{thm}


\begin{proof}
Let us write $X=X_n$ and $K=K_{X_n}$.  Recall that the tangent bundle $T_X$ is $\mu_{-K}$-stable for $n=7,8$, and so is $\sigma_{-tK}$-semistable for $t \gg 0$ by \cite[Lemma 4.8]{LM2}.

Suppose there exists some $t>0$ such that $T_X$ is strictly $\sigma_{-tK}$-semistable,  with an associated destabilizing exact sequence in $\mathcal{B}_{-K}$
\begin{equation}\label{eq:AG54-77-1}
0 \to A \to T_X \to Q \to 0.
\end{equation}
Since $\phi_{t}(T_X) = \phi_{t}(Q)$, we must have $-Kc_1(Q)>0$; note that we also have $-Kc_1(A)>0$ since $A \in \Tc_{-K}$ and $A$ must have positive rank.  When $n=8$, we have $K^2=1$, and so $K^2=-Kc_1(T_X) =-Kc_1(Q) -Kc_1(A)\geq 2$, a contradiction.  Hence $T_X$ is never strictly $\sigma_{-sK}$-semistable for any $s>0$ and so is $\sigma_{-sK}$-stable for all $s>0$ when $n=8$.  

From now on, suppose $n=7$, i.e.\ $K^2=2$.  Then $-Kc_1(A)=1=-Kc_1(Q)$, and so $A$ must be a $\mu_{-K}$-semistable sheaf (since each $\mu_{-K}$-HN factor of $A$ must have $-Kc_1>0$).  Consider the long exact sequence of cohomology of \eqref{eq:AG54-77-1}
\[
0 \to Q^{-1} \to A \overset{\beta}{\to} T_X \to Q^0 \to 0.
\]
If $\ch_0(\image \beta)=1$, then $\image \beta$ is a rank-1 subsheaf of $T_X$ with $\mu_{-K}(\image \beta) \geq 1$ while $\mu_{-K}(T_X)=1$, contradicting the $\mu_{-K}$-stability of $T_X$.  Hence $\ch_0(\image \beta)=2$, which implies $\ch_0(A) \geq 2$ and $Q^0$ is a torsion sheaf.  Also, since $-Kc_1(\image \beta) >0$, we have $-Kc_1(Q^0)=1$ or $0$.  We now divide into these two cases.


\noindent
\textbf{Case 1: $-Kc_1(Q^0)=0$.} In this case, $Q^0$ is supported in dimension $0$ and we have $-Kc_1(\image \beta)=2$ and hence $-Kc_1(Q^{-1})=-1$.  This means that $Q^{-1}$ has exactly one $\mu_{-K}$-HN factor $R_1$ with $-Kc_1=-1$, and possibly a $\mu_{-K}$-HN factor $R_0$ with $-Kc_1=0$.  Thus
\begin{align*}
    \ch_2(Q^{-1}) &= \ch_2(R_0) + \ch_2(R_1) \\
    &\leq \frac{c_1(R_0)^2}{2\ch_0(R_0)} + \frac{c_1(R_1)^2}{2\ch_0(R_1)} \text{ by the Bogomolov-Gieseker inequality} \\
    & \leq \frac{(Kc_1(R_0))^2}{2K^2\ch_0(R_0)} +  \frac{(Kc_1(R_1))^2}{2K^2\ch_0(R_1)} \text{ by the Hodge Index Theorem} \\
    &\leq \frac{1}{4}.
\end{align*}
Since $\ch_2(Q^{-1}) \in \tfrac{1}{2}\mathbb{Z}$, we obtain $\ch_2(Q^{-1})\leq 0$.  It follows that $\ch_2(Q) = \ch_2(Q^0)-\ch_2(Q^{-1}) \geq 0$, and so $\ch_2(A) = \ch_2(T_X)-\ch_2(Q) \leq -9$. 

 Now consider the wall equation $\rho_t(A) = \rho_t(T_X)$, i.e.\ 
\begin{align*}
    \frac{\ch_2(T_X) - t^2K^2}{tK^2} &= \frac{\ch_2(A) - \tfrac{t^2K^2}{2}\ch_0(A)}{t} \\
    \Rightarrow \frac{-9-2t^2}{2t} &\leq \frac{-9-\ch_0(A)t^2}{t} \\
    \Leftrightarrow  (\ch_0(A)-1)t^2 \leq -\tfrac{9}{2}
\end{align*}
which is impossible, and so this case cannot occur.

\noindent
\textbf{Case 2: $-Kc_1(Q^0)=1$.}  In this case, $-Kc_1(\image \beta)=1$ and $-Kc_1(Q^{-1})=0$.  The argument in the proof of Lemma \ref{lem:AG54-112-1} now shows that if $Q^{-1}$ is nonzero, $T_X$ would be destabilized by the subobject $\image \beta$ in $\mathcal{B}_{-K}$ with respect to $\sigma_{-tK}$, contradicting the assumption that $T_X$ is strictly $\sigma_{-tK}$-semistable.  Hence $Q^{-1}=0$.

On the other hand, $Q^0$ is a 1-dimensional sheaf with a pure 1-dimensional quotient sheaf $Q'$ such that $-Kc_1(Q')=-Kc_1(Q^0)=1$.  By Lemma \ref{lem:AG54-73-1}, we have $\ch_2(Q)=\ch_2(Q^0)\geq \ch_2(Q') \geq -\tfrac{1}{2}$, implying $\ch_2(A) = \ch_2(T_X)-\ch_2(Q) \leq -8.5$.  The wall equation $\rho_t(A) = \rho_t(T_X)$ now implies
\begin{align*}
    \frac{-9-2t^2}{2t} &\leq \frac{-8.5-\ch_0(A)t^2}{t} \\
    \Rightarrow  (\ch_0(A)-1)t^2 &\leq -4
\end{align*}
which is again impossible.

Hence $T_X$ is $\sigma_{-sK}$-stable for all $s>0$ when $n=7$.
\end{proof}

For del Pezzo surfaces $X_n$ where $n \geq 3$, it will be helpful to have an upper bound for the $\ch_2$ of semistable sheaves that is sharper than the usual Bogomolov-Gieseker inequality.  The principle behind the following result appears to be known to experts: the case  for K3 surfaces was  stated  in \cite[Section 6]{ABL}, while the same ideas appeared in the proof of \cite[Proposition 7.2]{levine2019brill}.

\begin{lem}\label{lem:AG54-108-1}
Suppose $(X,D)$ is a pair in either of these two cases:
\begin{itemize}
    \item[(1)] $X$ is a K3 surface and $D$ is an ample divisor on $X$.
    \item[(2)] $X$ is a del Pezzo surface and $D$ is an ample divisor on $X$ such that $DK_X <0$ (e.g.\ $D=-K_X$).
\end{itemize}

Then for any  $\mu_D$-stable torsion-free sheaf $A$ on $X$, we have
\[
\ch_2(A) \leq \frac{\ch_1(A)^2}{2\ch_0(A)} + \delta_X \left( \frac{1}{\ch_0(A)} - \ch_0(A) \right)
\]
where
\[
\delta_X = \begin{cases} 1 &\text{ in case (1)} \\
\tfrac{1}{2} &\text{ in case (2)}
\end{cases}.
\]
\end{lem}

\begin{proof}
From the short exact sequence of sheaves
\[
0 \to A \to A^{\ast\ast} \to T \to 0
\]
where $T$ is some 0-dimensional sheaf, we know  $A$ and  $A^{\ast\ast}$ have the same $\ch_0, \ch_1$ while $\ch_2(A)\leq \ch_2(A^{\ast\ast})$.  Therefore, it suffices for us to consider the case where $A$ is locally free.  From Hirzebruch-Riemann-Roch (e.g.\ see \cite[p.132]{huybrechts2006fourier}), we have
\[
\chi (A, A) = \int_X \ch(A) \ch(A^\ast) \mathrm{td}(T_X). 
\]
Recall that 
\begin{align*}
    \mathrm{td}_0(T_X)&=1 \\
    \mathrm{td}_1(T_X) &= \tfrac{-K}{2} \\
    \mathrm{td}_2(T_X) &= \tfrac{1}{12} ( K^2 + c_2(T_X) ) \\
  &= 1+p_a    \text{\quad from \cite[Appendix A, Example 4.1.2]{Harts}}
\end{align*}
where $p_a$ denotes the arithmetic genus of $X$.
Therefore,
\begin{align}
       \chi (A, A) &= \int_X \ch(A) \ch(A^\ast) (1, \tfrac{-K}{2}, 1 + p_a ) \notag\\
    &= \int_X \left( \ch_0(A)^2, 0, 2\ch_0(A)\ch_2(A) - \ch_1(A)^2 \right) \left( 1, \tfrac{-K}{2}, 1+p_a \right) \notag\\
    &= \ch_0(A)^2 (1+p_a) + 2\ch_0(A)\ch_2(A) - \ch_1(A)^2 \notag\\
    &= 2\ch_0(A) \left( \ch_2(A) - \frac{\ch_1(A)^2}{2\ch_0(A)} + \frac{\ch_0(A)}{2} ( 1+p_a ) \right). \label{eq:AG54-108-1}
\end{align}

Now, 
\begin{align*}
       \chi (A, A) &= \hom (A, A) - \mathrm{ext}^1(A,A) +\mathrm{ext}^2(A,A).
\end{align*}
Since $A$ is a $\mu_D$-stable torsion-free sheaf, we have $\hom (A, A)=1$.  By Serre duality, 
\[
\Ext^2(A,A) \cong \Hom (A, A \otimes \omega_X).
\]

\noindent
\textbf{Case (1):} In this case, we have $\hom (A, A \otimes \omega_X) = \hom (A,A)=1$ so $\chi (A,A) \leq 2$.

\noindent
\textbf{Case (2):} In this case, observe 
\[
Dc_1(A \otimes \omega_X) = D (c_1(A)+\ch_0(A)K ) = Dc_1(A) +\ch_0(A)DK < Dc_1(A).
\]
It follows that $A, A \otimes \omega_X$ are two $\mu_D$-stable sheaves with $\mu_D (A) > \mu_D (A \otimes \omega_X)$, and so  $\Hom (A, A \otimes \omega_X)=0$ and hence $\mathrm{ext}^2(A,A)=0$.  This gives $\chi (A,A) \leq 1$.

In summary, we have 
\[
\chi (A,A) \leq  \delta_X' \text{\quad where \quad} \delta_X' := \begin{cases}
    2 &\text{ in case (1)} \\
    1 &\text{ in case (2)}
    \end{cases}.
\]
Along with  \eqref{eq:AG54-108-1}, this gives us 
\[
\ch_2(A) \leq \frac{c_1(A)^2}{2\ch_0(A)} + \frac{\delta_X'}{2}\cdot \frac{1}{\ch_0(A)} - \frac{1+p_a}{2}\cdot \ch_0(A).
\]
Since $p_a=1$ when $X$ is K3 and $p_a=0$ when $X$ is del Pezzo, the result follows.
\end{proof}


We now turn our attention to the del Pezzo surface $X_n$ where $K^2=3$, i.e.\ $n=6$.  In this case, the wall equation $\rho_t(A)=\rho_t(T_X)$ from a destabilizing sequence $0 \to A \to T_X \to Q \to 0$ in $\Bc_{-K}$ in the sense of \ref{para:outermostwall-existence} reads
\begin{align}
    \rho_t (A)=\rho_t(T_X) & \Leftrightarrow \frac{\ch_2(A) - \tfrac{\omega^2}{2}\ch_0(A)}{\omega c_1(A)} = \frac{\ch_2(T_X)-\tfrac{\omega^2}{2}\ch_0(T_X)}{\omega c_1(T_X)} \notag \\
    &\Leftrightarrow \frac{\ch_2(A) -\tfrac{3}{2}\ch_0(A)t^2}{(-Kc_1(A))t} = \frac{-\tfrac{15}{2} - 3t^2}{3t}.  \label{eq:AG54-104-1}
\end{align}



\begin{lem}\label{lem:AG54-105-1}
Let $X=X_n$ where $0 \leq n \leq 8$.  Suppose $T_X$ is not strictly $\sigma_{-tK}$-semistable for all $t>0$ and $0 \to A \to T_{X} \to Q \to 0$ is a destabilizing sequence in $\Bc_{-K}$ corresponding to the outermost wall as in \ref{para:outermostwall-existence} that is simultaneously a short exact sequence in $\Coh (X)$.  Then 

(1) if $A$ is $\mu_{-K}$-semistable, then it is locally free; 

(2) when $n=4$, $5$, or $6$, if $\ch_0(A)=1$ then $Q$ is a rank-one torsion-free sheaf.
\end{lem}


Recall that a subsheaf $F$ of a torsion-free sheaf $E$ is called \emph{saturated} if the quotient sheaf $E/F$ is torsion-free (used in the proof of part (2) below).

\begin{proof}
For (1), since $A$ is a torsion-free sheaf, we have the commutative diagram
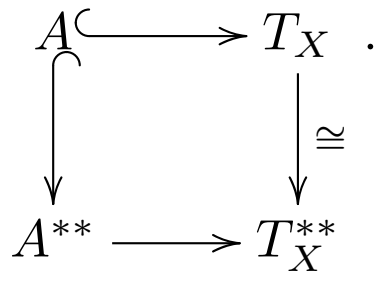
\[
\xymatrix{
  A \ar@{^{(}->}[r] \ar@{^{(}->}[d] & T_X \ar[d]^\cong \\
  A^{\ast\ast} \ar[r] & T_X^{\ast\ast}
}.
\]
Hence the lower horizontal map $A^{\ast\ast} \to T_X^{\ast\ast}$ is injective.  Let us denote the composite map 
 $A^{\ast\ast} \to T_X^{\ast\ast} \overset{\thicksim}{\to} T_X$ by $j$, and write $Q' = \cokernel j$.    Applying the octahedral axiom to the composition $A \to A^{\ast\ast} \overset{j}{\to} T_X$ gives the diagram
 
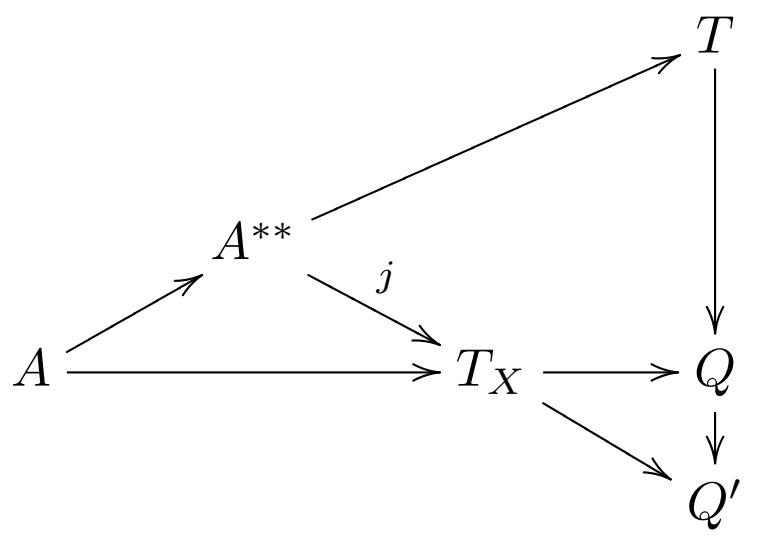
\[
 \xymatrix@R=1em{
 & & & T \ar[ddd] & \\
 & & & \\
 & A^{\ast\ast} \ar[uurr] \ar[dr]^j & &  \\
 A \ar[ur] \ar[rr] & & T_X \ar[r] \ar[dr] & Q \ar[d] \\
  & & & Q' 
 }
 \]
 where $T$ denotes the cokernel of $A \to A^{\ast\ast}$. 

 Now, since $A$ is $\mu_{-K}$-semistable torsion-free and $-Kc_1(A)>0$, the same hold for $A^{\ast\ast}$ and so $A^{\ast\ast} \in \Tc_{-K}$.  Also, since $Q'$ is a quotient sheaf of $Q$ from the diagram above, it also lies in $\Tc_{-K}$.  Hence the exact triangle $A^{\ast\ast} \to T_X \to Q' \to A^{\ast\ast}[1]$ corresponds to a short exact sequence 
 \[
 0 \to A^{\ast\ast} \to T_X \to Q' \to 0
 \]
 in $\Coh (X)$ that is also a short exact sequence in $\Bc_{-K}$.  If $A$ is not locally free, i.e.\ $T \neq 0$, then  $\rho_t (A^{\ast\ast}) > \rho_t(A)=\rho_t(T_X)$,  contradicting the assumption that $A$ is a destabilizing subobject of $T_X$ at the outermost mini-wall.     Hence $A$ must be locally free.

 Next we consider (2). Suppose $Q$ has a nonzero subsheaf $R$ supported in dimension 1.  We divide into three cases:

\begin{itemize}
    \item[(i)] $n=6$ (i.e.\ $K^2=3$): since $-Kc(A)>0$, we have
    \[
    3=K^2=-Kc_1(T_X)>-Kc_1(Q)>-Kc_1(Q/R)>0
    \]
    and so $-Kc_1(Q/R)$ is necessarily equal to $1$. Then $\mu_{-K}(T_X)=\tfrac{3}{2} > 1=\mu_{-K}(Q/R)$ which contradicts the $\mu_{-K}$-stability of $T_X$. 
    \item[(ii)] $n=5$ (i.e.\ $K^2=4$): As in case (i), we have 
    \[
    4=K^2=-Kc_1(T_X)>-Kc_1(Q)>-Kc_1(Q/R)>0
    \]
    and so $-Kc_1(Q/R)=2$ or $1$. We now have $\mu_{-K}(T_X)=2 \geq \mu_{-K}(Q/R)$ which again contradicts the $\mu_{-K}$-stability of $T_X$.  
    \item[(iii)] $n=4$ (i.e.\ $K^2=5$): the same argument as in the previous two cases gives $-Kc_1(Q/R)\leq 3$ which implies $-Kc_1(Q^*)\geq -3$.  Dualizing the short exact sequence of sheaves $0 \to A \to T_X \to Q \to 0$ gives the long exact sequence
    \[
    0 \to Q^\ast \to \Omega_X \to A^\ast \to \cdots
    \]
    from which we see $Q^\ast$ is a saturated invertible subsheaf of $\Omega_X$ of degree at least $-3$ with respect to $-K$, contradicting Proposition \ref{prop:K25678} (i) below.
\end{itemize}
 Hence any torsion subsheaf of $Q$ must be supported in dimension 0.  
 

 Suppose $Q$ has a nonzero subsheaf $G$ supported in dimension 0.  Then $Q/G$ is a quotient of $Q$ in $\Bc_{-K}$ with $\rho_t(T_X)=\rho_t(Q) > \rho_t(Q/G)$ (since $Q, Q/G$ have the same $\ch_0,  \ch_1$ while $\ch_2(Q)>\ch_2(Q/G)$), contradicting the assumption that $Q$ is a quotient object of $T_X$ at the outermost mini-wall.  Hence $Q$ must be a torsion-free sheaf.
\end{proof}

The following result of Jie Liu allows us to rule out a particular type of destabilizing subobjects on $X_6$.  Note that a cubic surface in $\mathbb{P}^3$ is isomorphic to $\mathbb{P}^2$ blown up at 6 points in general position, which is $X_6$ in our notation \cite[3.4.2]{liu2018geometry}.

\begin{prop}\cite[3.4.3 Proposition]{liu2018geometry}\label{pro:JieLiu-1}
Let $X$ be a cubic surface and $L \subset \Omega_X$ a saturated invertible subsheaf.  Then $-Kc_1(L)\leq -3$.
\end{prop}


\begin{lem}\label{lem:AG54-104-1}
Let $X=X_n$ where $n=6$.  Then in a destabilizing sequence $0 \to A \to T_{X} \to Q \to 0$ in $\Bc_{-K}$ corresponding to the outermost wall as in \ref{para:outermostwall-existence}, we must have $-Kc_1(A)=1$.
\end{lem}


\begin{proof}
Since $A$ is a nonzero torsion-free sheaf, we have $-Kc_1(A)>0$.  On the other hand, since  $\phi_t(A)=\phi_t (T_X)=\phi_t(Q)$, we cannot have $-Kc_1(Q)=0$.  Hence $-Kc_1(A)=1$ or $2$.  Consider the long exact sequence of sheaves 
\begin{equation}\label{eq:AG54-101-1}
    0 \to Q^{-1} \to A \overset{\beta}{\to} T_X \to Q^0 \to 0
\end{equation}
where $Q^i := H^i(Q)$.

Suppose $-Kc_1(A)=2$. Then $-Kc_1(\image \beta) \geq 2$ since $-Kc_1(Q^{-1}) \leq 0$.  If $\ch_0(\image \beta)=1$, then $\mu_{-K}(\image \beta) \geq 2 >\tfrac{3}{2}=\mu_{-K}(T_X)$ and, by the $\mu_{-K}$-stability of $\image \beta$ and $T_X$, we would have $\Hom (\image \beta,T_X)=0$ which is a contradiction.  Hence $\ch_0(\image \beta)=2$ which implies  $\ch_0(A) \geq 2$ and that $Q^0$ is a torsion sheaf.  Since $-Kc_1(Q)=1$, either $-Kc_1(Q^0)=1$ or $-Kc_1(Q^0)=0$.
\begin{itemize}
    \item If $-Kc_1(Q^0)=1$, then we have $\ch_2(Q^0) \geq -\tfrac{1}{2}$ by Lemma \ref{lem:AG54-73-1}.
    \item If $-Kc_1(Q^0)=0$, then $\ch_2(Q^0) \geq 0$.
\end{itemize}
Therefore, $\ch_2(Q^0) \geq -\tfrac{1}{2}$ overall which gives $\ch_2(\image \beta) = \ch_2(T_X)-\ch_2(Q^0) \leq -\tfrac{15}{2}+\tfrac{1}{2}=-7$.  The wall equation \eqref{eq:AG54-104-1} now becomes
\begin{align}
\frac{-\tfrac{15}{2}-3t^2}{3t} &= \frac{\ch_2(A)-\tfrac{3}{2}t^2\ch_0(A)}{2t}  \notag \\
\Leftrightarrow \text{\quad} ( \tfrac{3}{4}\ch_0(A)-1)t^2 &= \tfrac{1}{2}\ch_2(A)+\tfrac{5}{2}  \label{eq:AG54-101-2}
\end{align}

Now recall $-Kc_1(\image \beta) \geq 2$ from above.  If $-Kc_1(\image \beta)=2$, then $-Kc_1(Q^{-1})=0$ and $Q^{-1}=0$ by Lemma \ref{lem:AG54-112-1}.  Then $A=\image \beta$ and $\ch_2(A)=\ch_2(\image \beta) \leq -7$. 
In this case, the left-hand side of \eqref{eq:AG54-101-2} is positive (since $\ch_0(A) \geq 2$) while the right-hand side of \ref{eq:AG54-101-2} is negative, which is impossible.  So the case $-Kc_1(\image \beta)=2$ cannot occur.

From now on, suppose $-Kc_1(\image \beta) \geq 3$.  Then $-Kc_1(\image \beta)=3$ since $\image \beta$ is a subobject of $T_X$ in $\Bc_{-K}$ and $-Kc_1(T_X)=K^2=3$.  This means $-Kc_1(Q^{-1})=-1$, and so $Q^{-1}$ has exactly one $\mu_{-K}$-HN factor with $\mu_{-K}<0$, say $R_1$.  Then $-Kc_1(R_1)=-1$ and $c_1(R_1)^2 \leq \tfrac{(-Kc_1(R_1))^2}{K^2}= \tfrac{1}{3}$ by the Hodge index theorem.  By Lemma \ref{lem:nKDDsmod2}, we have $c_1(R_1)^2 \leq -1$ and so the Bogomolov-Gieseker inequality gives $\ch_2(R_1) \leq \tfrac{c_1(R_2)^2}{2\ch_0(R_1)}<0$, implying $\ch_2(R_1) \leq -\tfrac{1}{2}$.  Note that $Q^{-1}$ may also have one HN factor with $\mu_{-K}=0$, say $R_0$; the Hodge index theorem and the Bogomolov-Gieseker inequality together give $\ch_2(R_0) \leq 0$.  Hence $\ch_2(Q^{-1}) \leq - \tfrac{1}{2}$ overall.  Since $-Kc_1(Q^0)=0$, the sheaf $Q^0$ must be supported in dimension 0, and so  $\ch_2(Q) =\ch_2(Q^0)-\ch_2(Q^{-1}) \geq  \tfrac{1}{2}$, which in turn gives $\ch_2(A)=\ch_2(T_X)-\ch_2(Q) \leq -\tfrac{15}{2}-\tfrac{1}{2}=-8$.  Now, as in the previous paragraph, the left-hand side of the wall equation \eqref{eq:AG54-101-2} is positive while  the right-hand side is negative, again giving us a contradiction. 

Hence $-Kc_1(A)$ must be $1$.
\end{proof}

\begin{cor}\label{cor:AG54-104-1}
Let $X=X_n$ where $n=6$.  Then we cannot have a destabilizing sequence $0 \to A \to T_{X} \to Q \to 0$ in $\Bc_{-K_X}$ corresponding to the outermost wall as in \ref{para:outermostwall-existence} where $H^{-1}(Q)=0$ and $\ch_0(A)=1$.
\end{cor}

\begin{proof}
Suppose 
\begin{equation}\label{eq:TXdestabseq}
0 \to A \to T_{X} \to Q \to 0
\end{equation}
is a destabilizing sequence in $\Bc_{-K}$ corresponding to the outermost wall as in \ref{para:outermostwall-existence} such that $H^{-1}(Q)=0$ and $\ch_0(A)=1$.  Then \eqref{eq:TXdestabseq} is also a short exact sequence in $\Coh (X)$ and, by Lemma \ref{lem:AG54-105-1}, $A$ must be a line bundle while $Q$ must be a rank-1 torsion-free sheaf.  By Lemma \ref{lem:AG54-104-1}, we also know $-Kc_1(A)=1$ and $-Kc_1(Q)=2$.

Applying $\HHom(-, \OO_X)$ to \eqref{eq:TXdestabseq} gives us the long exact sequence of sheaves
\[
0 \to Q^\ast \to \Omega_X \overset{k}{\to} A^\ast \to \EExt^1(Q,\OO_X) \to 0.
\]
Since $Q^\ast$ is  locally free, by Proposition \ref{pro:JieLiu-1} we have $-Kc_1(Q^\ast)\leq -3$, i.e.\ $-Kc_1(Q) \geq 3$.  It follows that  $-Kc_1(A)=0$ and hence $\phi_t(T_X)=\phi_t(A)=1$ in the destabilizing sequence \eqref{eq:TXdestabseq}, which is impossible.
\end{proof}


We are now ready to prove the following result on the Bridgeland stability of $T_X$ when $X=X_6$.

\begin{prop}\label{pro:AG54-104-1}
Let $X=X_n$ where $n=6$, i.e.\ $K^2=3$.  Suppose $0 \to A \to T_{X} \to Q \to 0$ is a destabilizing sequence in $\Bc_{-K_X}$ corresponding to the outermost wall $t$ as in \ref{para:outermostwall-existence}.  Then 
\begin{itemize}
    \item[(i)] $H^{-1}(Q)=0$, and $A$ is a rank-two  $\mu_{-K}$-stable locally free subsheaf of $T_X$ with $-Kc_1(A)=1$, while $Q$ is a pure 1-dimensional Gieseker semistable sheaf with $-Kc_1(Q)=2$.
    \item[(ii)] The triple $(c_1(A)^2, \ch_2(A), t)$ must be one of the following six possibilities:
    \[
    (-1, -1, \tfrac{\sqrt{3}}{2}), \,\, (-1, -\tfrac{3}{2}, \tfrac{\sqrt{2}}{2}), \,\, (-1, -2, \tfrac{1}{2}), \,\, (-3, -\tfrac{3}{2}, \tfrac{\sqrt{2}}{2}), \,\, (-3, -2, \tfrac{1}{2}), \,\, (-5, -2, \tfrac{1}{2}).
    \]
\end{itemize}
As a result, if $t>\tfrac{\sqrt{3}}{2}$ then $T_X$ is $\sigma_{-tK}$-semistable.
\end{prop}

\begin{proof}
Suppose  $0 \to A \to T_{X} \to Q \to 0$ 
is a destabilizing sequence in $\Bc_{-K_X}$ as described.  Then $-Kc_1(A)=1$ by Lemma \ref{lem:AG54-104-1}.  Since $A \in \Tc_{-K}$ is a torsion-free sheaf, this means $A$ has exactly one $\mu_{-K}$-HN factor, and so is $\mu_{-K}$-semistable; that $-Kc_1(A)=1$ also implies $A$ cannot have more than one Jordan-H\"{o}lder factor, and so is $\mu_{-K}$-stable.   Consider the long exact sequence of sheaves 
\begin{equation}
    0 \to Q^{-1} \to A \overset{\beta}{\to} T_X \to Q^0 \to 0
\end{equation}
where $Q^i := H^i(Q)$.

Suppose $\ch_0(\image \beta)=1$.  By the $\mu_{-K}$-stability of $T_X$, we have $-Kc_1(\image \beta)=1$ and so $-Kc_1(Q^{-1})=0$.  By Lemma \ref{lem:AG54-112-1}, this means $Q^{-1}=0$ and 
$A=\image \beta$ is a rank-1 subsheaf of $T_X$; however, such a destabilizing sequence cannot exist by Corollary \ref{cor:AG54-104-1}.  Hence the case $\ch_0(\image \beta)=1$ cannot occur and we have $\ch_0(\image \beta)=2$.


Since $\image \beta$ is a $\Bc_{-K}$-suboject of $T_X$, we have $1 \leq -Kc_1(\image \beta)\leq 3$.  Our next step is to  rule out the cases $-Kc_1(\image \beta)=2,3$.

Suppose $-Kc_1(\image \beta)=3$.  Then $-Kc_1(Q^0)=0$ and $Q^0$ is supported in dimension 0.  It follows that $\ch_2(Q^0)\geq 0$ and $\ch_2(\image \beta)\leq \ch_2(T_X)=-\tfrac{15}{2}$.  On the other hand, we have $-Kc_1(Q^{-1})=-Kc_1(A)+Kc_1(\image \beta)=-2$.  Now, let $R_0$ denote the $\mu_{-K}$-HN factor of $Q^{-1}$ with $\mu_{-K}=0$ if it exists; note that $\ch_2(R_0) \leq 0$ and $R_0$.  In addition to $R_0$, $Q^{-1}$ can have either one or two $\mu_{-K}$-HN factors with $\mu_{-K}<0$.
\begin{itemize}
    \item If $Q^{-1}$ has  one HN factor with $\mu_{-K}<0$, say $R_1$, then $-Kc_1(R_1)=-2$ and $c_1(R_1)^2 \leq \tfrac{4}{3}$ by the Hodge index theorem.  By Lemma \ref{lem:nKDDsmod2}, it follows that $c_1(R_1)^2 \leq 0$.  The Bogomolov-Gieseker inequality then gives $\ch_2(R_1) \leq 0$.  Thus $\ch_2(Q^{-1})\leq 0$ overall.
    \item If $Q^{-1}$ has two HN factors with $\mu_{-K} < 0$, say $R_1$ and $R_2$, then we must have $-Kc_1(R_1)=-Kc_1(R_2)=-1$.  For each $i$, we have $c_1(R_i)^2\leq \tfrac{1}{3}$ by the Hodge index theorem; Lemma \ref{lem:nKDDsmod2} then gives $c_1(R_i)^2 \leq -1$, and the Bogomolov-Gieseker inequality gives $\ch_2(R_i)  <0$ and hence $\ch_2(R_i) \leq -\tfrac{1}{2}$.  Hence $\ch_2(Q^{-1})\leq -1$.
\end{itemize}
Regardless of which of these two cases $Q^{-1}$ is in, we have $\ch_2(Q^{-1})\leq 0$, and so $\ch_2(Q) \geq 0$ and hence $\ch_2(A) \leq \ch_2(T_X)=-\tfrac{15}{2}$.  The wall equation \eqref{eq:AG54-104-1} now becomes
\begin{equation}\label{eq:AG54-104-3}
\left( \tfrac{3}{2}\ch_0(A)-1\right)t^2 = \ch_2(A) + \tfrac{5}{2}.
\end{equation}
In this equation, the left-hand side is positive since $\ch_0(A) \geq \ch_0(\image \beta)=2$; on the other hand, the right-hand side is negative since $\ch_2(A) \leq -\tfrac{15}{2}$.  Hence the case $-Kc_1(\image \beta)=3$ cannot occur.

Now suppose $-Kc_1(\image \beta)=2$.  Then $Q^0$ is a torsion sheaf with $-Kc_1(Q^0)=1$.  Lemma \ref{lem:AG54-73-1} then gives us $\ch_2(Q^0) \geq -\tfrac{1}{2}$.  Consequently, $\ch_2(\image \beta)\leq \ch_2(T_X)+\tfrac{1}{2}=-7$.  We can now analyze the HN filtration of $Q^{-1}$ as we did in the previous paragraph.  As before, $Q^{-1}$ may have an HN factor $R_0$ with $\mu_{-K}=0$ (i.e.\ $R_0$ could be zero).  On the other hand, since $-Kc_1(Q^{-1})=-Kc_1(A)+Kc_1(\image \beta)=-1$, $Q^{-1}$ must have exactly one HN factor $R_1$ with $\mu_{-K}<0$ and $-Kc_1(R_1)=-1$.  Now $c_1(R_1)^2 \leq \tfrac{1}{3}$ by the Hodge index theorem, and so by Lemma \ref{lem:nKDDsmod2} we have $c_1(R_1)^2 \leq -1$.  The Bogomolov-Gieseker inequality then gives $\ch_2(R_1) <0$ and thus $\ch_2(R_1) \leq -\tfrac{1}{2}$, which in turn gives $\ch_2(Q^{-1}) \leq -\tfrac{1}{2}$; subsequently we have $\ch_2(A)=\ch_2(Q^{-1}) + \ch_2(\image \beta) \leq -\tfrac{15}{2}$.  In the wall equation \eqref{eq:AG54-104-3}, the left-hand side is  now positive while the right-hand side is negative, and so the case $-Kc_1(\image \beta)=2$ cannot occur, either.  We can now conclude $-Kc_1(\image \beta)=1$.

Since $-Kc_1(A)>0$, this means $-Kc_1(Q^{-1})=0$.  By Lemma \ref{lem:AG54-112-1}, $Q^{-1}$ must be zero.
Hence  $A$ is a rank-two subsheaf of $T_X$.  As a consequence, $Q=Q^0$ is a torsion sheaf with $-Kc_1(Q)=2$.  Moreover, $Q$ must be a pure 1-dimension sheaf for, if not, then $Q$ has a pure 1-dimensional quotient sheaf $Q'$ with $\ch_1(Q)=\ch_1(Q')$ and $\ch_2(Q) > \ch_2(Q')$.  This implies $\rho_t(T_X)=\rho_t(Q) > \rho_t(Q')$, contradicting the assumption that $Q$ is a destabilizing quotient of $T_X$ at the outermost mini-wall.  Similarly, $Q$ must be Gieseker semistable, for any quotient sheaf $Q'$ of $Q$ that destabilizes $Q$ with respect to Gieseker stability would satisfy $\rho_t(Q)> \rho_t(Q')$.

Also, since $A$ is a $\mu_{-K}$-stable, by Lemma \ref{lem:AG54-105-1} it is locally free.  This completes the proof of (i).

To see (ii), note that $\rho_t(T_X)=\rho_t(Q)$ by assumption which gives
\[
\frac{-\tfrac{15}{2}-3t^2}{3t}=\frac{\ch_2(Q)}{2t} \Leftrightarrow t^2 = -\tfrac{5}{2}-\tfrac{1}{2}\ch_2(Q)
\]
and so $-\tfrac{5}{2}-\tfrac{1}{2}\ch_2(Q)>0$, i.e.\ $\ch_2(Q)<-5$.  This in turn means $\ch_2(A)>\ch_2(T_X)+5=-\tfrac{5}{2}$, and so $\ch_2(A) \geq -2$.  Along with Lemma \ref{lem:AG54-108-1}, we now have
\[
-2 \leq \ch_2(A) \leq \tfrac{c_1(A)^2}{4}-\tfrac{3}{4}
\]
while the Hodge index theorem and Lemma \ref{lem:nKDDsmod2} together give $c_1(A)^2 \leq -1$.   
Since $c_1(A)^2 \leq -1$ and $c_1(A)^2 \equiv -Kc_1(A) \equiv 1$ modulo 2, $c_1(A)^2$ can only be $-1, -3, -5$.  As a result,
\[
\ch_2(A) = \begin{cases}
    -2 &\text{ when } c_1(A)^2=-5 \\
    -2 \text{ or } -\tfrac{3}{2} &\text{ when }c_1(A)^2=-3  \\
    -2, -\tfrac{3}{2}, \text{ or } -1 &\text{ when }c_1(A)^2=-1
\end{cases}.
\]
For each of these cases, we can use the wall equation $\rho_t(A)=\rho_t(T_X)$ which gives
\[
\frac{-\tfrac{15}{2}-3t^2}{3t}=\frac{\ch_2(A)-3t^2}{t} \Leftrightarrow 2t^2 - \tfrac{5}{2}=\ch_2(A)
\]
to complete the proof of (ii).
\end{proof}

%

\begin{lem}\label{lem:ch2bd-AG54-134-1}
Let $X$ be a smooth projective surface, and $\omega$ an (integral) ample divisor on $X$.  Suppose $E$ is a nonzero torsion-free sheaf in either $\Tc_\omega$ or $\Fc_\omega$.  Then 
\begin{equation}\label{eq:ch2bd-AG54-134-1}
\ch_2(E) \leq \frac{ |\omega c_1(E)|^3}{2\omega^2}.
\end{equation}

\end{lem}


\begin{proof}
Suppose $R_1, \cdots, R_m$ are the $\mu_\omega$-HN factors of $E$ satisfying $\omega c_1 \neq 0$.  In  case  $E \in \Fc_\omega$, let us write $R_0$ (possibly zero) to denote the HN factor  with $\omega c_1=0$; in case of $E \in \Tc_\omega$, let us set $R_0=0$.  The Hodge index theorem and the BG inequality togther give $\ch_2(R_0) \leq 0$.  On the other hand, for each $i \neq 0$ we have by the Hodge index theorem
\[
c_1(R_i)^2 \leq \frac{(\omega c_1(R_i))^2}{\omega^2} \leq \frac{(\omega c_1(E))^2}{\omega^2}.
\]
Then, using the Bogomolov-Gieseker inequality  as well, we obtain  
\begin{align*}
\ch_2(E) = \sum_{i=0}^m \ch_2(R_i) &\leq \sum_{i=1}^m \ch_2(R_i) \\
&\leq \sum_{i=1}^m \frac{ c_1(R_i)^2}{2\ch_0(R_i)} \\
&\leq \sum_{i=1}^m \frac{(\omega c_1(E))^2}{2\omega^2 \ch_0(R_i)}\text{\quad from above}\\
&\leq  \sum_{i=1}^m \frac{(\omega c_1(E))^2}{2\omega^2} \\\
&\leq \frac{|\omega c_1(E)|^3}{2\omega^2}
\end{align*}
where, in the last inequality, we used the fact that $m \leq |\omega c_1(E)|$ since $\omega$ is an integral divisor.
\end{proof}

\begin{prop}\label{pro:AG54-133-1}
Let $X=X_n$ where $n=5$, i.e.\ $K^2=4$.  Suppose $0 \to A \to T_{X} \to Q \to 0$ is a destabilizing sequence in $\Bc_{-K}$ corresponding to the outermost wall as in \ref{para:outermostwall-existence} for some $t>0$.  Then  $t \leq 1$, i.e.\ if $t>1$ then $T_X$ is $\sigma_{-tK}$-semistable.
\end{prop}

\begin{proof}
As usual, in the long exact sequence of sheaves 
\begin{equation}
    0 \to Q^{-1} \to A \overset{\beta}{\to} T_X \to Q^0 \to 0
\end{equation}
where $Q^i := H^i(Q)$ we have $\ch_0(\image \beta)=1$ or $2$.  The wall equation $\rho_t(A) = \rho_t(T_X)$ when $K^2=4$ reads
\begin{align}
    \frac{\ch_2(T_X) - t^2K^2}{tK^2} &= \frac{\ch_2(A)-\tfrac{t^2K^2}{2}\ch_0(A)}{-Kc_1(A)t}\notag\\
 \Leftrightarrow   \frac{-6-4t^2}{4t} &= \frac{\ch_2(A) - 2\ch_0(A) t^2}{-Kc_1(A)t} \notag \\
  \Leftrightarrow  \left( \tfrac{2}{-Kc_1(A)}\ch_0(A) - 1 \right) t^2 &= \tfrac{\ch_2(A)}{-Kc_1(A)} + \tfrac{3}{2}. \label{eq:AG54-86-wq}
\end{align}

\noindent
\textbf{If $\ch_0(\image \beta)=1$:} The slope stability of $T_X$ gives $-Kc_1(\image \beta) < \mu_{-K}(T_X)=2$ and so $-Kc_1(\image \beta)=1$.  Since $-Kc_1(A)>0$, this implies $-Kc_1(Q^{-1})=0$ and so $Q^{-1}=0$ by Lemma \ref{lem:AG54-112-1} and consequently $A=\image \beta$.  By the Hodge index theorem, we have $c_1(A)^2 \leq \tfrac{1}{4}$; together with Lemma \ref{lem:nKDDsmod2}, we obtain $c_1(A)^2 \leq -1$.  The Bogomolov-Gieseker inequality then gives $\ch_2(A) \leq -\tfrac{1}{2}$.  

The wall equation \eqref{eq:AG54-86-wq} now gives
\[
    t^2 = \ch_2(A) + \tfrac{3}{2} \leq 1
\]
and so $t>1$ would ensure the $\sigma_{-tK}$-semistability of $T_X$.

\noindent
\textbf{If $\ch_0(\image \beta)=2$:} In this case, $Q^0$ is a torsion sheaf with $0 \leq -Kc_1(Q^0)\leq 3$, and so $1 \leq -Kc_1(\image \beta) \leq 4$.  We divide into cases (i) through (iv) according to the value of $-Kc_1(\image \beta)$.

\noindent
Case (i): $-Kc_1(\image \beta)=1$.  In this case, that $-Kc_1(A)>0$ forces $-Kc_1(Q^{-1})=0$, and so $Q^{-1}=0$ by Lemma \ref{lem:AG54-112-1}.  Hence $A=\image \beta$.  Even though $A$ has rank 2, that $-Kc_1(A)=1$ and $A \in \Tc_{-K}$ means $A$ cannot have more than one $\mu_{-K}$-HN factor, i.e.\ $A$ is $\mu_{-K}$-semistable (and in fact, $\mu_{-K}$-stable).  The Hodge index theorem together with Lemma \ref{lem:nKDDsmod2} together give $c_1(A)^2 \leq -1$; the Bogomolov-Gieseker inequality subsequently gives $\ch_2(A) \leq -\tfrac{1}{4}$ and thus $\ch_2(A) \leq -\tfrac{1}{2}$ since $\ch_2$ is a half-integer.  The wall equation \eqref{eq:AG54-86-wq} now gives 
\[
 3t^2 = \ch_2(A) + \tfrac{3}{2} \leq 1,
\]
i.e.\ $t^2 \leq \tfrac{1}{3}$.

\noindent
Case (ii): $-Kc_1(\image \beta)=2$.  In this case, $\image \beta$ has rank 2 and $-Kc_1=2$; if it's not $\mu_{-K}$-semistable, it would have exactly two HN factors each of which has rank 1 and $-Kc_1=1$, which is impossible.  Hence $\image \beta$ is $\mu_{-K}$-semistable.  The Hodge index theorem along with Lemma \ref{lem:nKDDsmod2} gives us $c_1(\image \beta)^2 \leq 0$, and so by the Bogomolov-Gieseker inequality we have $\ch_2(\image \beta) \leq 0$.  Also, since $-Kc_1(A) >0$, we have $-Kc_1(Q^{-1})=0$ or $-1$.

\underline{If $-Kc_1(Q^{-1})=0$}: then $Q^{-1}=0$ by Lemma \ref{lem:AG54-112-1}, in which case $A=\image \beta$ has rank 2, is $\mu_{-K}$-semistable, and has $-Kc_1=2$ and $c_1(A)^2 \leq 0$.  We further divide into two cases:
\begin{itemize}
    \item Case (a): $A$ is $\mu_{-K}$-stable.  By the inequality in Lemma \ref{lem:AG54-108-1}, we have 
    \[
    \ch_2(A) \leq 0+\tfrac{1}{2} \left( \tfrac{1}{2} -2 \right)=-\tfrac{3}{4}
    \]
    and hence $\ch_2(A) \leq -1$.
    \item Case (b): $A$ is strictly $\mu_{-K}$-semistable.  Then $A$ must have exactly two JH factors, say $R_1, R_2$, each of which having rank 1 and $-Kc_1=1$.   For each $i$, the Hodge index theorem and Lemma \ref{lem:nKDDsmod2} give $c_1(R_i)^2 \leq -1$.  Applying the Bogomolov-Gieseker inequality to each $R_i$ now gives  $\ch_2(R_i) \leq -\tfrac{1}{2}$, and so we obtain $\ch_2(A) \leq -1$.
\end{itemize}
Overall, we have $\ch_2(A) \leq -1$ and the wall equation \eqref{eq:AG54-86-wq} gives
\[
t^2 = \tfrac{\ch_2(A)}{2}+\tfrac{3}{2} \leq 1.
\]

\underline{If $-Kc_1(Q^{-1})=-1$}: then $Q^{-1}$ has exactly one HN factor $R_1$ with $-Kc_1<0$, and we necessarily have $-Kc_1(R_1)=-1$.  Note that $Q^{-1}$ might also have a nonzero HN factor $R_0$ with $-Kc_1=0$.  In any case, the Hodge index theorem and Lemma \ref{lem:nKDDsmod2} together give $c_1(R_1)^2 \leq -1$, and so from the Bogomolov-Gieseker inequality we have $\ch_2(R_1) \leq -\tfrac{1}{2\ch_0(R_1)}$ and hence $\ch_2(R_1) \leq -\tfrac{1}{2}$.  Since $R_0$ is necessarily $\mu_{-K}$-semistable, the Hodge index theorem and the Bogomolov-Gieseker inequality together give $\ch_2(R_0)\leq 0$.  As a result, we have $\ch_2(Q^{-1}) = \ch_2(R_0)+\ch_2(R_1) \leq -\tfrac{1}{2}$ and hence $\ch_2(A) =\ch_2(Q^{-1})+\ch_2(\image \beta)\leq -\tfrac{1}{2}$.  The wall equation \eqref{eq:AG54-86-wq} now gives
\[
(2\ch_0(A)-1)t^2 = \ch_2(A) + \tfrac{3}{2} \leq 1 
\]
Since $\ch_0(A) = \ch_0(Q^{-1})+\ch_0(\image \beta) \geq 3$, it follows that $t^2 \leq \tfrac{1}{5}$.

Case (iii): $-Kc_1(\image \beta)=3$.  We will show that this case is impossible.  Note that $-Kc_1(Q^0)=1$ and so $\ch_2(Q^0) \geq -\frac{1}{2}$ by Lemma \ref{lem:AG54-73-1}.  Thus $\ch_2(\image \beta) \leq \ch_2(T_X)+\tfrac{1}{2}=-\tfrac{11}{2}$.  On the other hand, since $-Kc_1(A)>0$ we must have $0 \geq -Kc_1(Q^{-1}) \geq -2$.  By Lemma \ref{lem:ch2bd-AG54-134-1}, 
\[
\ch_2(Q^{-1}) \leq \tfrac{|-Kc_1(Q^{-1})|^3}{2K^2} \leq 1
\]
and so $\ch_2(A) = \ch_2(Q^{-1})+\ch_2(\image \beta) \leq -\tfrac{9}{2}$.  Now consider the wall equation \eqref{eq:AG54-86-wq}.  Since $1 \leq -Kc_1(A) \leq 3$, the right-hand side is 
\[
\tfrac{\ch_2(A)}{-Kc_1(A)} + \tfrac{3}{2} \leq \tfrac{-9/2}{3} + \tfrac{3}{2}=0
\]
while the coefficient of $t^2$ on the left-hand side is (noting $\ch_0(A) \geq \ch_0(\image \beta)=2$)
\[
\tfrac{2}{-Kc_1(A)}\ch_0(A) - 1 \geq \tfrac{2}{3}\cdot 2 -1
\]
which is positive, giving us a contradiction.

Case (iv): $-Kc_1(\image \beta)=4$.  We will also show that this case is impossible.  Since $-Kc_1(Q^0)=0$, the sheaf $Q^0$ must be supported in dimension 0 with $\ch_2 \geq 0$, and so $\ch_2(\image \beta) \leq \ch_2(T_X)=-6$.  Also note that $0 \geq -Kc_1(Q^{-1}) \geq -3$ in this case.

\underline{If $-Kc_1(Q^{-1})=0$}: then $-Kc_1(Q)=0$ and hence $\rho_\omega (Q)=\infty$, a contradiction.  

\underline{If $-Kc_1(Q^{-1})=-1$}: then $Q^{-1}$ has exactly one $\mu$-HN factor (say $R_1$) with $-Kc_1<0$.  Since the HN factor with $-Kc_1 =0$ - if it exists - has $\ch_2 \leq 0$, we obtain $\ch_2(Q^{-1}) \leq \ch_2(R_1)$.  Since $-Kc_1(R_1)=-1$, by the Hodge index theorem and Lemma \ref{lem:nKDDsmod2} together we get   $c_1(R_1)^2 \leq -1$.  The BG inequality then gives $\ch_2(R_1) \leq \frac{-1}{2\ch_0(R_1)} <0$, hence $\ch_2(Q^{-1}) \leq -\frac{1}{2}$.  Overall, we obtain $\ch_2(A) \leq -\tfrac{1}{2} -6 \leq -\frac{13}{2}$ while $-Kc_1(A)=3$.  The right-hand side of the wall equation \eqref{eq:AG54-86-wq} is now bounded from above by $-\tfrac{13}{6}+\tfrac{3}{2} = -\tfrac{2}{3}<0$ while the left-hand side is positive, which is impossible.

\underline{If $-Kc_1(Q^{-1})=-2$}: then we have the following two sub-cases
\begin{itemize}
    \item[(a)] $Q^{-1}$ has a single $\mu$-HN factor with nonzero $-Kc_1$ which has $-Kc_1=-2$;
    \item[(b)] $Q^{-1}$ has exactly two $\mu$-HN factors with nonzero $-Kc_1$, each of which has $-Kc_1=-1$.
\end{itemize}
In case (a), the single HN factor $R$ has even $c_1^2$ by Lemma \ref{lem:nKDDsmod2} while the Hodge index theorem gives $c_1(R)^2 \leq \frac{4}{K^2}=1$, and so $c_1(R)^2 \leq 0$.  The BG inequality then gives $\ch_2(Q^{-1})\leq \ch_2(R) \leq 0$.  In case (b), each of the two HN factors with nonzero $-Kc_1$, say $R_1$ and $R_2$, has odd $c_1^2$, while the Hodge index theorem gives $c_1^2(R_i) \leq \frac{1}{4}$ for each $i$.  Hence $c_1(R_i)^2 \leq -1$ for each $i$, and the BG inequality gives $\ch_2(R_i) \leq \frac{-1}{2\ch_0(R_i)}$, meaning we have  $\ch_2(R_i) \leq -\tfrac{1}{2}$ for each $i$.  This means $\ch_2(Q^{-1})\leq \ch_2(R_1)+\ch_2(R_2) \leq -1$.    In both case (a) and case (b), we have $\ch_2(Q^{-1}) \leq 0$.  This implies $\ch_2(A) = \ch_2(Q^{-1}) + \ch_2(\image \beta) \leq -6$ while $-Kc_1(A)=2$.  It is then easy to see that the right-hand side of the  wall equation \eqref{eq:AG54-86-wq}  is negative while the left-hand side is positive, making this case impossible.

\underline{If $-Kc_1(Q^{-1})=-3$}: then we have three cases below to consider depending on how $-Kc_1(Q^{-1})$ is partitioned.
\begin{itemize}
    \item[(a')] $Q^{-1}$ has a single HN factor with nonzero $-Kc_1$, and it has $-Kc_1=-3$;
    \item[(b')] $Q^{-1}$ has exactly two HN factors with nonzero $-Kc_1$, and they have  $-Kc_1$ equal to $-1, -2$;
    \item[(c')] $Q^{-1}$ has exactly three HN factors with nonzero $-Kc_1$, each  of which has $-Kc_1=-1$.
\end{itemize}
In case (a'), the unique HN factor $R$ with nonzero $-Kc_1$ must have odd $c_1^2$ by Lemma \ref{lem:nKDDsmod2}.  The Hodge index theorem gives $c_1(R)^2 \leq \tfrac{9}{4}$ and so $c_1(R)^2 \leq 1$.  Then the BG inequality gives $\ch_2(R) \leq \tfrac{1}{2}$ and so $\ch_2(Q^{-1}) \leq \tfrac{1}{2}$.  In case (b'), let us write $R_1, R_2$ to denote the two HN factors with nonzero $-Kc_1$, and suppose $-Kc_1(R_1)=-1, -Kc_1(R_2)=-2$.  The same calculation as in case (b) of the previous paragraph gives us $\ch_2(R_1) \leq -\tfrac{1}{2}$, while that in case (a) of the previous paragraph gives us $\ch_2(R_2)\leq 0$.  Overall, we have $\ch_2(Q^{-1})\leq -\tfrac{1}{2}$ in case (b'). 
In case (c'), the calculation in case (b) of the previous paragraph gives $\ch_2 \leq -\tfrac{1}{2}$ for each of the three HN factors with $-Kc_1=-1$, and so $\ch_2(Q^{-1}) \leq -\tfrac{3}{2}$ and we have $\ch_2(A) \leq -\tfrac{15}{2}$ overall.  
Putting cases (a'), (b'), and (c') together, we see that when $-Kc_1(Q^{-1})=-3$, we have $\ch_2(Q^{-1}) \leq \tfrac{1}{2}$ and hence $\ch_2(A) \leq -\tfrac{11}{2}$ while $-Kc_1(A)=1$.  In the wall equation \eqref{eq:AG54-86-wq}, the right-hand side is again negative while the left-hand side positive, and so this case cannot occur, either.

Putting all the cases above together, we see that when $n=5$, i.e.\ $K^2=4$, the tangent bundle $T_X$ of the surface $X=X_n$ is $\sigma_{-tK}$-semistable whenever $t>1$.
\end{proof}

The following result and its argument were communicated to us by Jie Liu (see also \cite[Proposition 6.3]{liu2018geometry}):

\begin{prop}[Liu]
\label{prop:K25678}
Let $X=X_n$ be a del Pezzo surface.  
\begin{itemize}
    \item[(i)] For $4 \leq K^2 \leq 8$,  we have
\[
\max{\{ -Kc_1(L) : L \text{ is a saturated invertible subsheaf of }\Omega_X \}} = -4.
\]
\item[(ii)] For $K^2=6, 7$, the sheaf $L=O_X(-2l +2e_i)$ is the unique saturated invertible subsheaf of $\Omega_X$ attaining the maximum value $-Kc_1(L)=-4$ in part (i).

\item[(iii)] For $K^2=5$, the sheaves $L=\cO_X(-4l+2\sum_{j} e_j)$ and $L'=\cO_X(-2l+2e_i)$ are the only saturated invertible subsheaves of $\Omega_X$ attaining the maximum value $-Kc_1(P)=-4$ in part (i).
\end{itemize}

\end{prop}

\begin{proof}
(i) Let $L$ be any rank-one subsheaf of $\Omega_X$.  

\underline{When $K^2=6, 7, 8$}: the $\mu_{-K}$-(semi)stability of $\Omega_X$ gives
\[
-Kc_1(L)=\mu_{-K}(L) \leq -4.
\]
On the other hand, we know  $L'=O_X(-2l+2e_i)$ is a saturated invertible subsheaf of $\Omega_X$ from \cite[Example 6.1(1)]{liu2021stability} and
\[
-Kc_1(L')=(3l-\sum_j e_j)(-2l+2e_i)=-4, 
\]
proving part (i) for $K^2=6, 7, 8$.

\underline{When $K^2=5$}: the $\mu_{-K}$-stability of $\Omega_X$ gives $-Kc_1(L)\leq -3$.  We will now show, that when $L$ is a saturated invertible subsheaf of $\Omega_X$, the equality $-Kc_1(L)=-3$ cannot be attained: Let us write 
\[
c_1(L) = -al - \sum_j b_j e_j
\]
for some integers $a, b_j$.  By the correspondence between rank-1 invertible subsheaves of $\Omega_X$ and global sections of twists of $\Omega_{\mathbb{P}^2}$ that vanish only in codimension 2 \cite[Section 1]{liu2021stability}, we know $\Omega_{\mathbb{P}^2}(a)$ must have a nonzero global section, i.e.\ there must be a nonzero morphism $\OO_X \to \Omega_{\mathbb{P}^2}(a)$ which in turn means the degree of $\Omega_{\mathbb{P}^2}(a)$ with respect to $-K$, i.e.  $-Kc_1(\Omega_{\mathbb{P}^2}(a))$, must be nonnegative; this implies   $a \geq 2$.  

If $a=2$, then $L$ must be either $O_X(-2l+2e_i)$ or $O_X(-2l)$ by \cite[Example 6.1(1)]{liu2021stability}; since $-K(-2l+2e_i)=-4$ from above and $-K(-2l)=-6$, this case cannot happen, so let us assume $a \geq 3$ from here on.

Now 
\[
-3 = -Kc_1(L)=-3a-\sum_j b_j, 
\]
i.e.\
\begin{equation}\label{eq:AG54-132-1}
3a+\sum_j b_j = 3.
\end{equation}
Since $L$ is not isomorphic to $O_X(-2l+2e_i)$ (their first Chern classes have different intersection numbers with $-K$), the line bundle $\omega_X \otimes L^\ast \otimes O_X(2l-2e_i)$ must have a nonzero global section by \cite[Lemma 3.4.2]{liu2018geometry} (see also \cite[Lemma 1]{fahlaoui1989stabilite}), and so there is an effective divisor $C_i$ on $X$ such that  $C_i \thicksim K -c_1(L)+(2l-2e_i)$.   Note that the divisor $-K-e_i$ is nef since its intersection with every generator of the effective monoid (described in Theorem \ref{thm:fanosurf-effconegenerators})  is nonnegative.  As a result,
\begin{align*}
0 &\leq (-K-e_i)C_i  \\
&= \left( 3l-2e_i - \sum_{j \neq i} e_j\right) \left( (a-1)l + (b_i-1)e_i + \sum_{j\neq i} (b_j+1)e_j \right) \\
&= 3(a-1)+2(b_i-1)+\sum_{j \neq i}(b_j+1) \\
&= 3a+ \left( \sum_j b_j\right) + b_i - 3 - 2 + (n-1) \text{ where $n=9-K^2=4$}\\
&= b_i+1.
\end{align*}
This shows that $b_i \geq -1$ for all $i$, and so from \eqref{eq:AG54-132-1} we have $3a=3-\sum_j b_j \leq 3+4=7$, i.e.\ $a \leq 2$, which is a contradiction.  Hence $-Kc_1(L)=-3$ cannot occur and statement (i) is proved.

\underline{When $K^2=4$}: the same argument as in the $K^2=5$ case show that 
\[
\max{\{ -Kc_1(L) : L \text{ is a saturated invertible subsheaf of }\Omega_X \}} \leq -3.
\]
Let $P$ be a saturated invertible subsheaf of $\Omega_X$ such that $-Kc_1(P)=-3$.  Let us write $c_1(P)=-al - \sum_j b_j e_j$ as in the previous case.  Then the same calculation  as in the $K^2=5$ case (taking the intersection of $-K-e_i$ with $K-c_1(P)+(2l-2e_i)$) shows that $b_i\geq -2$ for all $i$.  The same reasoning  as in the $K^2=5$ case also rules out the possibility $a=2$, while $a\geq 5$ cannot hold for otherwise $-Kc_1(P)=3a+\sum_j b_j \geq 15+4(-2)>3$, a contradiction.  Hence there are only two possibilities left to consider: (1) $a=3$, $\sum_j b_j=-6$, and (2) $a=4$, $\sum_jb_j=-9$.

Note that $l-e_i$ is also a nef divisor.  As in the previous case, there is an effective divisor on $X$ linearly equivalenct to $K-c_1(P)+(2l-2e_i)$, and so 
\[
\begin{split}
0 & \leq (l-e_i)(K-c_1(P)+(2l-2e_i))\\
& \leq (l-e_i)\left((a-1)l+(b_i-1)e_i+\sum_{j\neq i}(b_j+1)e_j\right)\\
&\leq (a-1)+(b_i-1) = a+b_i -2.
\end{split}
\]
This implies that $a+b_i\geq 2$ for all $i$.   In case (1), since  $\sum_{j=1}^5 b_j = -6$ and $b_j \geq -2$ for all $j$, we must have $b_j=-2$ for some $j$.  This means $a+b_j=1$, which contradicts $a+b_i \geq 2$ for all $i$.  So we must be in case (2).  In case (2), that $\sum_{j=1}^5 b_j=-9$ and $b_j \geq -2$ for all $j$ implies the $b_j's$ must be $-2, -2, -2, -2, -1$, and so 
\[
P\simeq \cO_X(-4l+e_i+2\sum_{j\neq i}e_j)
\]
and 
By Fahlaoui’s example, as recalled explicitly in \cite[Example 6.1(2)]{liu2021stability},  there is no subsheaf \(P\subset\Omega_X^1\) with  
\[
P\simeq\mathcal O_X\!\left(-4l+e_i+2\sum_{j\ne i}e_j\right),
\]
hence the actual bound for saturated subsheaf of $\Omega_X$ is $-4$.

(ii) From \cite[Example 6.1(1)]{liu2021stability}, we know $\OO_X (-2l+2e_i)$ is a saturated invertible subsheaf of $\Omega_X$.  Suppose $L$ is a saturated invertible subsheaf of $\Omega_X$ satisfying $-Kc_1(L)=-4$, and that $L$ is not isomorphic to $\OO_X(-2l+2e_i)$.  Let us write $c_1(L)=-al - \sum_j b_j e_j$ as in part (i), so that $-Kc_1(L)=-4$ becomes $3a+\sum_j b_j = 4$.  As in part (i) (see the argument in the $K^2=5$ case), we can assume $a \geq 3$, and there exists an effective divisor $C_i$ such that $C_i \thicksim K-c_1(L)+(2l-2e_i)$.  Since $-K-e_i$ is a nef divisor, we have
\begin{align*}
0 &\leq (-K-e_i)C_i \\
&= 3a+ \left( \sum_j b_j\right) + b_i - 5 + (n-1)  \\
&= b_i-1+(n-1)
\end{align*}
and so $b_i \geq 1-(n-1)$ for all $i$.

If $K^2=6$, then $n=3$ and $b_i \geq -1$.  Then $a=\tfrac{1}{3}\left( 4-\sum_j b_j\right) \leq \tfrac{7}{3}$ which is impossible.  If $K^2=7$, then $n=2$ and $b_i \geq 0$.  Then $a=\tfrac{1}{3}\left( 4-\sum_j b_j\right) \leq \tfrac{4}{3}$, which is also impossible.  So in either case, $L$ must be isomorphic to $\OO_X(-2l+2e_i)$.

(iii) Consider $K^2=5$.
Let $L$ be a saturated invertible subsheaf of $\Omega_X$ such that $-Kc_1(L)=-4$. 
As before, let us write $c_1(L)=-al-\sum_jb_je_j$.  

Suppose  $L\not\cong \cO_X(-2l+2e_i)$.  Then as in the $K^2=5$ case of part (i), we must have $a \geq 3$.  Calculating the intersection number of $(-K-e_i)$ with $K -c_1(L)+(2l-2e_i)$ as in the $K^2=5$ case in (i)  gives $b_i\geq -2$ for all $i$.  On the other hand, we must have $a < 5$, or else $-4=-Kc_1(L)= -3a-\sum_j b_j \leq (-3)5 +4\cdot 2 = -7$, which is impossible.  This leaves us with two cases:  (iii-1) $a=3$, $\sum_j b_j=-5$; and (iii-2) $a=4$, $\sum_jb_j=-8$. 

Calculating the intersection number of $(l-e_i)$ with $K-c_1(L)+(2l-e_i)$ as in the $K^2=4$ case of (i) now gives  $a+b_i\geq 2$ for all $i$, and so  case (iii-1) cannot happen.  That is, if $L$ is not isomorphic to $\OO_X (-2l+2e_i)$ then it is isomorphic to $\cO_X(-4l+2\sum_je_j)$.  From \cite[Example 6.1(1), (2)]{liu2021stability}, we know both $\OO_X (-2l+2e_i)$ and  $\cO_X(-4l+2\sum_je_j)$ are saturated invertible subsheaves of $\Omega_X$, and so we are done.
\end{proof}

\paragraph \label{para:wallcalc-01}
Let $X$ be a smooth projective surface, $H$ an ample class on $X$, and $E$ a nonzero torsion-free sheaf on $X$ lying in $\Bc_H$.  Let $\omega = tH$ where $t \in \mathbb{R}_{>0}$.  Suppose $0 \to A \to E \to Q \to 0$ is a short exact sequence in $\Bc_\omega$ where $A, Q \neq 0$.  Then
\begin{align}
 \rho_\omega (A) - \rho_\omega (E) &= \frac{\ch_2(A)-\tfrac{\omega^2}{2}\ch_0(A)}{\omega \ch_1(A)} - \frac{\ch_2(E) - \tfrac{\omega^2}{2}\ch_0(E)}{\omega \ch_1(E)} \notag \\
 &= \frac{tH^2}{2}\left( - \frac{\ch_0(A)}{H\ch_1(A)} + \frac{\ch_0(E)}{H\ch_1(E)}\right) + \frac{1}{t}\left( \frac{\ch_2(A)}{H\ch_1(A)} - \frac{\ch_2(E)}{H\ch_1(E)}\right)  \notag \\
  &= \frac{1}{t}\left( \frac{t^2H^2}{2} \left(  - \frac{1}{\mu_H (A)} + \frac{1}{\mu_H(E)}\right) + \left( \frac{\ch_2(A)}{H\ch_1(A)} - \frac{\ch_2(E)}{H\ch_1(E)}\right) \right) \label{eq:rhoAEdiff}
\end{align}
provided $H\ch_1(A), H\ch_1(E)$ are nonzero.  Assuming $E$ is $\mu_H$-semistable, from the long exact sequence in $\Coh (X)$
\begin{equation}\label{eq:AEQles}
0 \to Q^{-1} \to A \overset{\beta}{\to} E \to Q^0 \to 0
 \end{equation}
(where $Q^i$ denotes $H^i(Q)$) we have
\begin{align}
    0 < \mu_H (A) = \frac{Hc_1(A)}{\ch_0(A)} &= \frac{Hc_1(\image \beta)+Hc_1(Q^{-1})}{\ch_0(\image \beta) + \ch_0(Q^{-1})} \notag\\
    &\leq \mu_H(\image \beta) \notag\\
    &\leq \mu_H(E)  \label{eq:AEstringineq}
\end{align}
where the  last inequality follows from the $\mu_H$-semistability of $E$, and the second last inequality follows because $Hc_1(\image \beta), \ch_0(\image \beta)$ are both positive while $\ch_0(Q^{-1}) \geq 0$ and $Hc_1(Q^{-1}) \leq 0$.  As a result, we have
\[
 - \frac{1}{\mu_H (A)} + \frac{1}{\mu_H(E)} \leq 0
\]
in \eqref{eq:rhoAEdiff}.  

\subparagraph \label{para:rhoAEdiff-sub1} Note that if $\frac{1}{\mu_H (A)} = \frac{1}{\mu_H(E)}$, then $\mu_H(A)=\mu_H(E)$ and the string of inequalities above becomes a string of equalities,  forcing $\ch_0(Q^{-1})=0$, i.e.\ $Q^{-1}$ must be 0 in which case $A$ is a subsheaf of $E$.  If $E$ happens to be $\mu_H$-stable, then the condition $\mu_H(A)=\mu_H(E)$ forces $\ch_0(A)=\ch_0(E)$ and subsequently $Hc_1(A)=Hc_1(E)$, meaning $Q$ must be torsion sheaf supported in dimension 0.

\subparagraph  \label{para:rhoAEdiff-sub2} Suppose we are in the case $- \frac{1}{\mu_H (A)} + \frac{1}{\mu_H(E)} < 0$.  Then
\begin{equation}\label{eq:rhoAEcomparison}
\rho_\omega (A) > (=) \,\rho_\omega (E) \, \Leftrightarrow \, t^2 < (=) \,  \frac{ \tfrac{\ch_2(A)}{Hc_1(A)} - \frac{\ch_2(E)}{Hc_1(E)}}{\tfrac{\ch_0(A)}{Hc_1(A)}-\tfrac{\ch_0(E)}{Hc_1(E)}} \cdot \frac{2}{H^2}.
\end{equation}
In particular, whenever 
\begin{equation}\label{eq:t-generalbound-E}
 t^2 \geq  \frac{ \tfrac{\ch_2(A)}{Hc_1(A)} - \frac{\ch_2(E)}{Hc_1(E)}}{\tfrac{\ch_0(A)}{Hc_1(A)}-\tfrac{\ch_0(E)}{Hc_1(E)}} \cdot \frac{2}{H^2}
\end{equation}
we have $\rho_\omega (A) \leq \rho_\omega (E)$.

The next proposition will be useful when it comes to proving the Bridgeland stability of $T_X$ on $X=X_n$ for $n \neq 1$. It shows that when we strengthen the assumption of \cite[Lemma 4.8]{LM2} from twisted Gieseker semistability to slope stability, we can improve the conclusion from Bridgeland semistability to Bridgeland stability.

\begin{prop}\label{pro:Establehasoutermostwall}
Let $X$ be a smooth projective surface, $H$ an ample class on $X$, and $E$ a nonzero torsion-free sheaf in $\mathcal{T}_H$ that is $\mu_H$-stable.  Then there exists $t_1>0$ such that $E$ is $\sigma_{tH}$-stable for all $t>t_1$.  In particular, if $E$ is not $\sigma_{tH}$-stable for all $t>0$, then there exists $t_0>0$ such that $E$ is $\sigma_{tH}$-stable for all $t>t_0$,  strictly $\sigma_{t_0H}$-semistable, and $\sigma_{tH}$-unstable for all $0<t<t_0$, and $t_0$ is the outermost mini-wall in the sense of Lemma \ref{lem:goingup-v3}.
\end{prop}

Since slope stability implies twisted Gieseker stability, Proposition \ref{pro:Establehasoutermostwall} is a corollary of the `going-up lemma' in \cite[Lemma 7.11]{CLSY} (up to a $\widetilde{\mathrm{GL}}^+(2,\mathbb{R})$-action on the central charge).  Nonetheless, we include a direct proof here that is somewhat different from the argument presented in \cite{CLSY} to make clear  the connection with the outer-most mini-wall as presented in \ref{para:outermostwall-existence}.

\begin{proof}
Consider the usual long exact sequence of sheaves \eqref{eq:AEQles}.  If $Q^{-1} \neq 0$, then $\ch_0(Q^{-1})>0$ while $Hc_1(Q^{-1}) \leq 0$, and from the string of inequalities \eqref{eq:AEstringineq} we see $0< \mu_H(A) < \mu_H(E)$; from \eqref{eq:rhoAEdiff}, it follows that  $\rho_{tH}(A) < \rho_{tH}(E)$ for $t \gg 0$.  So let us assume $Q^{-1}=0$, i.e.\ $A=\image \beta$ is a subsheaf of $E$ from now on.

Suppose $\ch_0(A)=\ch_0(E)$.  Then $Q$ is a torsion sheaf.  If $\dimension Q=1$ then $\mu_H(A)<\mu_H(E)$ and again $\rho_{tH}(A)< \rho_{tH}(E)$ for $t \gg 0$; if $\dimension Q=0$ then $\rho_{tH}(A)< \infty = \rho_{tH}(Q)$ for all $t > 0$.

Suppose $\ch_0(A) < \ch_0(E)$.  Then $\mu_H(A) < \mu_H(E)$ by the $\mu_H$-stability of $E$ and we get  $\rho_{tH}(A) < \rho_{tH}(E)$ for $t \gg 0$ \eqref{eq:rhoAEdiff}.

Now, from \cite[Proposition 4.1(b)]{bayer2019stability}, we know there is a polynomial stability condition $\wt{\sigma}$ on $D^b(X)$ such that $E$ is $\wt{\sigma}$-stable if and only if it is $\sigma_{tH}$-stable for $t \gg 0$.  From the proof of \cite[Proposition 4.1]{bayer2019stability}, we know $\Bc_H$ coincides with the heart $\wt{\mathcal{P}}(0,1]$ where $\wt{\mathcal{P}}$ denotes the slicing of the polynomial stability condition $\wt{\sigma}$.  Our calculations in the previous paragraphs thus show that $E$ is a $\wt{\sigma}$-stable object in $\mathcal{P}(0,1]$, and so $E$ is $\sigma_{tH}$-stable for $t \gg 0$ by \cite[Proposition 4.1(b)]{bayer2019stability}.  This proves the first assertion.

For the second assertion, note  that since for any nonzero proper subobject $A$ of $E$ in $\mathcal{B}_H$ we have $\rho_{tH}(A) < \rho_{tH}(E)$ for $t \gg 0$, it follows that  $A$ must be in either case (ii) or case (iii) in the proof of Lemma \ref{lem:goingup-v3}. Therefore, assuming  
 $E$ is not $\sigma_{tH}$-stable for all $t>0$, $E$ must have a subobject $A$ that falls in case (iii).    Hence $E$ is not $\sigma_{tH}$-semistable for all $t>0$, and so there exists an outer-most mini-wall $t_0$ as in  Lemma \ref{lem:goingup-v3}.  This means that $E$ is strictly $\sigma_{t_0H}$-semistable and $\sigma_{tH}$-unstable for all $0<t<t_0$.  Since every nonzero subobject $A$ of $E$ falls in either case (ii) or case (iii), that $E$ is $\sigma_{tH}$-semistable for every $t>t_0$ also means that $E$ must be $\sigma_{tH}$-stable for every $t>t_0$, thus proving the second assertion in the proposition.
\end{proof}

Note that in the setting of Proposition \ref{pro:Establehasoutermostwall}, the value $t_0$ is exactly the outermost wall in the sense of \ref{para:outermostwall-existence}.

\paragraph[Matlab scripts] \label{para:Fano-Matlabcode-outline} Now suppose $X$ is a del Pezzo surface of the form $X_n$, i.e.\ $X$ is $\mathbb{P}^2$ blown up at $n$ points where $1 \leq n \leq 8$.  Recall from \ref{para:Fanosurf-general} that $T_X$ is $\mu_{-K}$-stable except when $n=1$, in which case $T_X$ is strictly $\mu_{-K}$-semistable.  So far, our main approach to finding a bound $t_0>0$ such that $t\geq t_0$ ensures the $\sigma_{-tK}$-semistability of $T_X$ has been based on first finding an upper bound for the second Chern character $\ch_2$ of the destabilizing subobject $A$ at the outermost mini-wall.  This usually involves first finding upperbounds for $\ch_2(Q^{-1})$ and $\ch_2(\image \beta)$.  In the case where we do not know whether $Q^{-1}$ or $\image \beta$ are $\mu_{-K}$-semistable sheaves, we have resorted to considering all possible HN filtrations of these sheaves and bounding the $\ch_2$ of the corresponding HN factors (e.g.\ in the analysis of $Q^{-1}$ in the proof of Proposition \ref{pro:AG54-104-1}).  As $K^2$ (which equals $9-n$) increases, however, the range of possible values of $-Kc_1(\image \beta)$ as well as that of $-Kc_1(Q^{-1})$ both increase, leading to more cases to check.  

In order to avoid checking all these cases manually where the computations are mechanical and vastly similar, we handwrote Matlab code to automate some calculations for us.  In the paragraphs below, we give a high-level picture of the functions of the different scripts we have written  and how they fit together.

The main script \url{wall_est.m} is a function that takes a positive integer $3 \leq \kappa \leq 8$ as input.  The integer $\kappa$ denotes   $K^2_X$ of a del Pezzo $X$ of the form $X_n$ where $n=9-K^2_X$.  For each possibility of $\ch_0(\image \beta), -Kc_1(\image \beta)$, and $-Kc_1(Q^{-1})$,  the Matlab script \url{wall_est.m}  calculates a value $t'_0$ such that $t^2\geq (t'_0)^2$ implies $\rho_\omega (A) \leq \rho_\omega (T_X)$ provided $A$ is a destabilizing subobject corresponding to the outermost wall in the sense of \ref{para:outermostwall-existence}.  The script then gives as output the maximum among all the $t_0'$, which is then an \emph{upper bound} for the outermost wall for the $\sigma_\omega$-semistability of $T_X$ on the surface $X$.  

The function \url{wall_est.m} depends on other scripts described below.  The Matlab codes of all these scripts  are available at \cite{Lo-Matlab-code-ch2}.  We have attempted to include comments in the Matlab code that allows the reader to follow the logic of the various scripts, including how mathematical results such as the Hodge index theorem, the BG inequality, and Lemmas \ref{lem:nKDDsmod2}, \ref{lem:AG54-73-1} and  \ref{lem:AG54-108-1} are used in calculations.
\begin{itemize}
    \item \url{intpartitions.m}: This was a Matlab script written by Wong \cite{WongIntegerPartitions}  that generates all unordered positive-integer partitions of a given positive integer.
     \item \url{ch2_uppbd_singleHN_withoutrank.m}: Given a $\mu_{-K}$-semistable torsion-free sheaf on $X$ where the rank of $E$ is not known, this script computes an upper bound of $\ch_2(E)$ using the Hodge index theorem, the Bogomolov-Gieseker inequality, the integrality of $-Kc_1$,  the half integrality of $\ch_2$ for sheaves on a surface, and Lemma \ref{lem:nKDDsmod2}.  
      \item \url{ch2_uppbd_singleHN_withrank.m}: Given a $\mu_{-K}$-semistable torsion-free sheaf on $X$ where the rank of $E$ is known, this script computes an upper bound of $\ch_2(E)$ using the Hodge index theorem, the Bogomolov-Gieseker inequality, the integrality of $-Kc_1$,  the half integrality of $\ch_2$ for sheaves on a surface, and both of   Lemmas \ref{lem:nKDDsmod2} and  \ref{lem:AG54-108-1}   
    \item \url{ch2_uppbd_Qm1.m}: Given a nonzero torsion-free sheaf in the torsion-free class $\mathcal{F}_{-K}$, this script computes an upper bound of $\ch_2(E)$ by considering the values of $-Kc_1$ of all possible $\mu_{-K}$-HN filtrations of $E$. 
    \item \url{ch2_uppbd_imbeta_rk2.m}: This script computes an upper bound of $\ch_2(\image \beta))$ when $\image \beta$ is of rank 2.  It uses Lemma \ref{lem:AG54-73-1} as part of the calculation.  
\end{itemize}


\begin{figure*}[t]
\[
\xymatrix{
 & \fbox{\url{wall_est.m}} \ar[dl] \ar[d] \ar@/^2pc/[dd]  \\
 \fbox{\url{ch2_uppbd_Qm1.m}} \ar[dd] \ar@/_2pc/[ddr] & \fbox{\url{ch2_uppbd_imbeta_rk2.m}}\ar[d]  \\
 & \fbox{\url{ch2_uppbd_singleHN_withrank.m}} \\
 \fbox{\url{intpartitions.m}} & \fbox{\url{ch2_uppbd_singleHN_withoutrank.m}}
}
\]
\caption{Dependence relations among the MATLAB scripts.  There is an arrow from script A to script B if A directly calls B.}
\label{fig:matlab-script-dependencies}
\end{figure*}
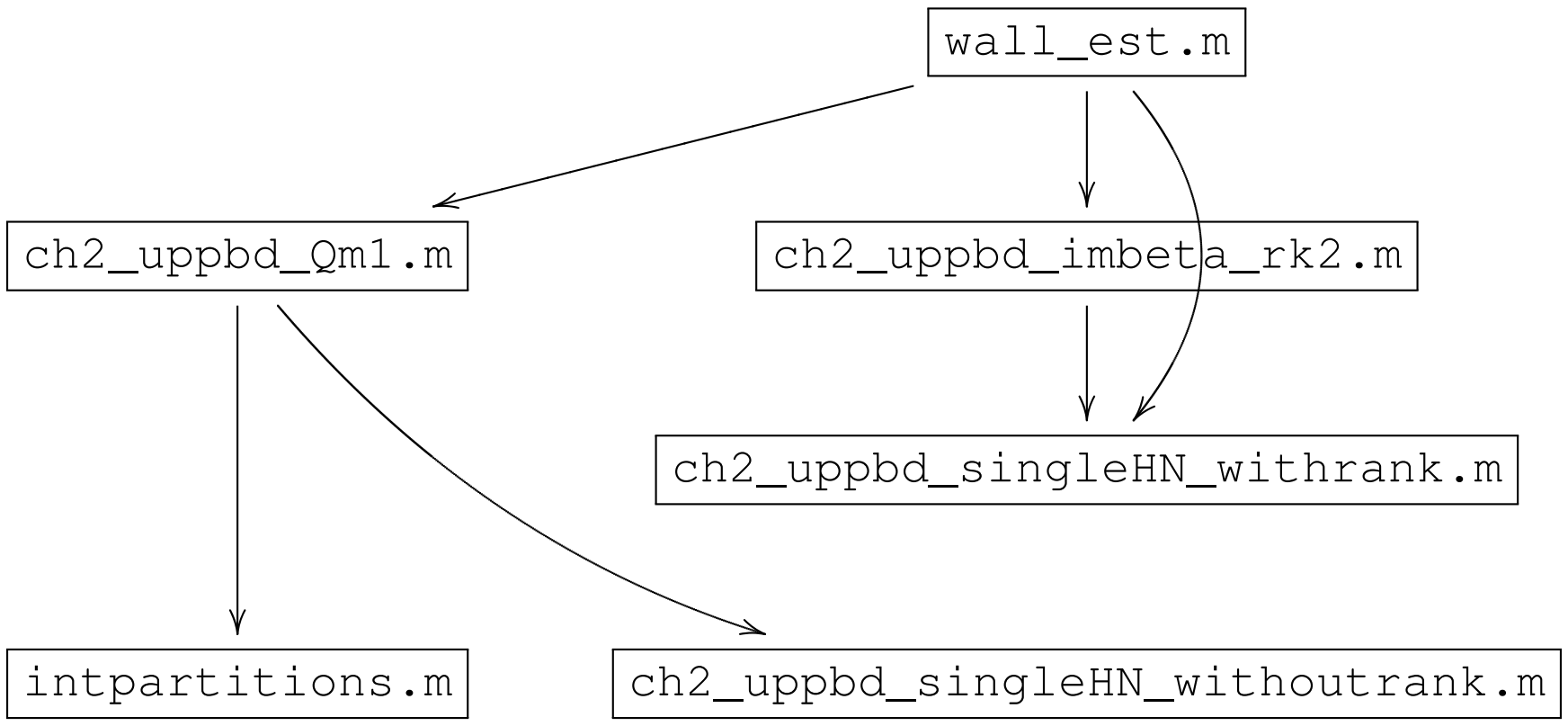

\paragraph Based on the calculation \eqref{eq:rhoAEdiff}, 
the wall equation $\rho_\omega (A) = \rho_\omega (T_X)$ can be rewritten as 
\begin{equation}\label{eq:walleq_C123}
C'_1 t^2 = C'_2 + C'_3
\end{equation}
where
\begin{align*}
    C'_1 &:= \left( \frac{(K^2/2)\ch_0(A)}{-Kc_1(A)} - \frac{(K^2/2)\ch_0(T_X)}{-Kc_1(T_X)}  \right) = \left( \frac{(K^2/2)\ch_0(A)}{-Kc_1(A)} - 1 \right), \\
    C'_2 &:= \frac{\ch_2(A)}{-Kc_1(A)}, \\
    C'_3 &:= -\frac{\ch_2(T_X)}{-Kc_1(T_X)}.
\end{align*}
Note that we can also write $C_1'$ as 
\[
C_1' = \frac{K^2}{2} \left( \frac{1}{\mu_{-K}(A)} - \frac{1}{\mu_{-K}(T_X)} \right). 
\]
By the $\mu_{-K}$-semistable stability of $T_X$ and \ref{para:wallcalc-01}, we have $0 < \mu_{-K}(A) \leq \mu_{-K}(T_X)$ and so $C_1' \geq 0$.  These observations will be used in the propositions below and in \url{wall_est.m}.


The script  \url{wall_est.m} now allows us to compute bounds of stability for $T_X$ on more del Pezzo surfaces.

\begin{prop}\label{prop:fano-four-wall-bound}
Let $X=X_n$ where $n=4$, i.e.\ $K^2=5$. Then 
 $T_X$ is $\sigma_{-tK}$-semistable for $t^2 \geq \tfrac{11}{10}$.
\end{prop}

\begin{proof}
Consider a destabilizing sequence
$0\to A\to T_X\to Q\to0$ at the outermost mini-wall $t$ as in
\ref{para:outermostwall-existence}  and the corresponding long exact sequence  in $\Coh (X)$
\[
0 \to Q^{-1} \to A \overset{\beta}{\to} E \to Q^0 \to 0.
 \]

From the 
 Matlab script \url{wall_est.m} we know  $t^2 \leq \tfrac{18}{5}$ and this upper bound comes from the case $(\ch_0(\image \beta), -Kc_1(\image \beta), -Kc_1(Q^{-1})) = (1, 2, 0)$.  Suppose we are in this case.  Then by Lemma \ref{lem:AG54-112-1}, $Q^{-1}=0$ and so $A = \image \beta$ is a rank-one subsheaf of $T_X$ with $-Kc_1(A)=2$.  Since $\operatorname{ch}_2(T_X)=-\tfrac92$,
the wall equation \eqref{eq:walleq_C123} can be rewritten as
\[
\tfrac14 t^2= \tfrac{\ch_2(A)}{2} + \tfrac{9}{10}. 
\]

The Hodge index theorem, the integrality of $-Kc_1(A)$ and the Bogomolov--Gieseker inequality
give $\operatorname{ch}_2(A)\leq0$.  Using the wall equation above, the  positivity of $t^2$ and
half-integrality of $\ch_2$ then give $-3/2\leq\operatorname{ch}_2(A)\leq0$.  By Lemma \ref{lem:AG54-105-1}(1), $A$ is a line bundle;   by Lemma \ref{lem:AG54-105-1}(2), $Q$ is a rank-one torsion-free sheaf. 
Dualizing the short exact sequence of sheaves $0 \to A \to T_X \to Q \to 0$ now gives the exact sequence of sheaves
\[
0 \to Q^\ast \to \Omega_X \to A^\ast
\]
meaning $Q^\ast$ is a saturated subsheaf of $\Omega_X$ with $-Kc_1=-3$, contradicting Proposition \ref{prop:K25678}(i).  Hence $T_X$ cannot fit in a destabilizing sequence with $(\ch_0(\image \beta), -Kc_1(\image \beta), -Kc_1(Q^{-1})) = (1, 2, 0)$.


From the output of \url{wall_est.m}, this means the outermost wall $t$ for $T_X$ now satisfies $t^2 \leq \frac{27}{20}$ and comes from the case where $(\ch_0(\image \beta), -Kc_1(\image \beta), -Kc_1(Q^{-1})) = (2, 3, 0)$.  Suppose we are in this case. Again by Lemma
\ref{lem:AG54-112-1}, $Q^{-1}=0$, so $A=\image\beta$ is a
rank-two subsheaf of $T_X$ with $-Kc_1(A)=3$. We claim that
$\ch_2(A)\leq-\tfrac12$.  To see this, we divide into two cases according to the $\mu_{-K}$-stability of $A$.

If $A$ is $\mu_{-K}$-stable, then the Hodge index theorem gives
\[
c_1(A)^2\leq \tfrac{(-Kc_1(A))^2}{K^2}=\tfrac95,
\]
so $c_1(A)^2\leq1$ by integrality. By Lemma \ref{lem:AG54-108-1}, we  have
\[
\ch_2(A)
\leq \frac{c_1(A)^2}{2\ch_0(A)}
     +\frac12\left(\frac{1}{\ch_0(A)}-\ch_0(A)\right)
\leq \frac14+\frac12\left(\frac12-2\right)
=-\frac12.
\]

If $A$ is not $\mu_{-K}$-stable, then it cannot be
$\mu_{-K}$-semistable since its rank and $-Kc_1$ are coprime.
Let $R_1,R_2$ be its $\mu_{-K}$-HN factors in decreasing order of $\mu_{-K}$, which must both
have rank one. Since $A\in\mathcal T_{-K}$, the degrees
$-Kc_1(R_1)$ and $-Kc_1(R_2)$ are positive integers; they sum
to $3$ and are strictly decreasing, so
\[
-Kc_1(R_1)=2,\qquad -Kc_1(R_2)=1.
\]
The Hodge index theorem gives
\[
c_1(R_1)^2\leq\frac45,\qquad
c_1(R_2)^2\leq\frac15.
\]
Hence $c_1(R_1)^2\leq0$. Also, $c_1(R_2)^2$ is odd by
Lemma \ref{lem:nKDDsmod2}, so $c_1(R_2)^2\leq-1$. The
Bogomolov--Gieseker inequality now gives
\[
\ch_2(R_1)\leq0,\qquad \ch_2(R_2)\leq-\frac12,
\]
and therefore
$\ch_2(A)=\ch_2(R_1)+\ch_2(R_2)\leq-\tfrac12$
in this case as well.

Thus regardless of the $\mu_{-K}$-stability of $A$ we have $\ch_2(A)\leq-\tfrac12$. Since
$\ch_0(A)=2$ and $-Kc_1(A)=3$, the wall equation
\eqref{eq:walleq_C123} gives
\[
\left(\frac53-1\right)t^2
=\frac{\ch_2(A)}3+\frac9{10}
\leq-\frac16+\frac9{10}
=\frac{11}{15},
\]
and so $t^2\leq\frac{11}{10}$.

All the remaining cases in the output of \url{wall_est.m}
either cannot occur or give upper bounds for $t^2$ of at most
$\frac35$. Hence the outermost mini-wall for $T_X$ satisfies
$t^2\leq\frac{11}{10}$. In particular, $T_X$ is
$\sigma_{-tK}$-semistable whenever $t \geq \sqrt{\frac{11}{10}}$.
\end{proof}

\begin{prop}\label{prop:fano-three-wall-bound}
Let $X=X_n$ where $n=3$, i.e.\  $K^2=6$. Then $T_X$ is
$\sigma_{-tK}$-semistable for $t \geq 1$.
\end{prop}

\begin{proof}
Consider a destabilizing sequence
$0\to A\to T_X\to Q\to0$ at the outermost mini-wall $t$ as in
\ref{para:outermostwall-existence}  and the corresponding long exact sequence in
$\Coh(X)$
\[
0\to Q^{-1}\to A\overset{\beta}{\to}T_X\to Q^0\to0
\]
where $Q^i := H^i(Q)$. From the output of \url{wall_est.m} for $K^2=6$, we know
$t^2\leq\tfrac54$ and this upper bound comes from the case
\[
(\ch_0(\image\beta),-Kc_1(\image\beta),-Kc_1(Q^{-1}))
=(2,4,0).
\]
(All the remaining cases either cannot occur or give upper
bounds for $t^2$ of at most $1$.) Suppose we are in this  case. By Lemma \ref{lem:AG54-112-1}, $Q^{-1}=0$, so
$A=\image\beta$ is a rank-two subsheaf of $T_X$ with
$-Kc_1(A)=4$. We claim that $\ch_2(A)\leq0$.

If $A$ is $\mu_{-K}$-stable, then the Hodge index theorem gives
\[
c_1(A)^2\leq\frac{(-Kc_1(A))^2}{K^2}
=\frac{16}{6}=\frac83,
\]
so $c_1(A)^2\leq2$ by integrality. By Lemma
\ref{lem:AG54-108-1}(2), we then have
\[
\ch_2(A)
\leq\frac{c_1(A)^2}{2\ch_0(A)}
  +\frac12\left(\frac{1}{\ch_0(A)}-\ch_0(A)\right)
\leq\frac24+\frac12\left(\frac12-2\right)
=-\frac14<0.
\]

If $A$ is not $\mu_{-K}$-stable, we first show that it must
be $\mu_{-K}$-semistable. Otherwise, its two $\mu_{-K}$-HN
factors have rank one and positive integral degrees since
$A\in\mathcal T_{-K}$. These degrees sum to $4$ and are
strictly decreasing, so they must be $3$ and $1$. The left
HN factor would then be a rank-one subsheaf of $T_X$ of
slope $3=\mu_{-K}(T_X)$, contradicting the
$\mu_{-K}$-stability of $T_X$.

Thus $A$ is strictly $\mu_{-K}$-semistable. Let $R_1,R_2$
be its $\mu_{-K}$-Jordan--H\"older factors. They have rank one
and slope $\mu_{-K}(A)=2$, so
\[
-Kc_1(R_1)=-Kc_1(R_2)=2.
\]
The Hodge index theorem gives $c_1(R_i)^2\leq4/6<1$,
hence $c_1(R_i)^2\leq0$ for each $i$. The
Bogomolov--Gieseker inequality now gives
\[
\ch_2(A)=\ch_2(R_1)+\ch_2(R_2)\leq0
\]
in this case as well.  That is, regardless of the $\mu_{-K}$-stability of $A$, we have $\ch_2(A) \leq 0$.

Since $\ch_2(T_X)=-3$, $\ch_0(A)=2$, and $-Kc_1(A)=4$,
the wall equation \eqref{eq:walleq_C123} gives
\[
\left(\frac32-1\right)t^2
=\frac{\ch_2(A)}4+\frac12\leq\frac12,
\]
and so $t^2\leq1$.  All the remaining cases in the output of \url{wall_est.m}
either cannot occur or give upper bounds for $t^2$ that are at most $1$. As a result, $T_X$ is
$\sigma_{-tK}$-semistable whenever $t \geq 1$.
\end{proof}

\paragraph[Regions of stability so far] For del Pezzo surfaces of the form $X=X_n$ where $0 \leq n \leq 8$, the current state of our knowledge of the regions of stability of $T_X$ on the ray $\{\sigma_{-tK}\}_{t>0}$ is as follows.
\begin{itemize}
    \item For $K^2=1, 2, 9$, we showed that the tangent bundle $T_X$ is $\sigma_{-tK}$-semistable for all $t>0$ (Theorems \ref{thm:AG54-53-1} and \ref{thm:TXn-stab-78}).
    \item For $K^2=3$, the tangent bundle $T_X$ is $\sigma_{-tK}$-semistable for all $t>\tfrac{\sqrt{3}}{2}$ (Proposition \ref{pro:AG54-104-1}).
    \item For $K^2=4$, the tangent bundle $T_X$ is $\sigma_{-tK}$-semistable for all $t>1$ (Proposition \ref{pro:AG54-133-1}).
    \item For $K^2=5$, additional analysis of the output of the Matlab script \url{wall_est.m} gave that $T_X$ is $\sigma_{-tK}$-semistable for $t> \sqrt{\frac{11}{10}}$ (Proposition \ref{prop:fano-four-wall-bound}). 
    \item For $K^2=6$, additional analysis of the Matlab output gave that $T_X$ is $\sigma_{-tK}$-semistable for $t> 1$ (Proposition \ref{prop:fano-three-wall-bound}). 
    \end{itemize}
These results, together with  upper bounds $t_0$ produced  by the Mathlab script \url{wall_est.m}, are  summarized in Table \ref{tab:wall-est-maximizers}. For each $1 \leq n \leq 6$, we also list in the table all possible cases of the triple
\[
 (\ch_0(\image \beta), -Kc_1(\image \beta), -Kc_1(Q^{-1})
\]
for which there might exist a destabilizing sequence 
\[
0 \to A\to T_X \to Q \to 0
\]
in $\mathcal{B}_{-K}$ with associated long exact sequence in $\Coh (X)$
\[
0 \to Q^{-1}\to A \overset{\beta}{\to} T_X \to Q^0 \to 0
\]
where $Q^i := H^i(Q)$.  These cases can be identified from the statement or the proof of the proposition listed under `Remark' in the cases $K^2=3, 4, 5$.  For the $K^2=6$ case, the triple $(2, 4, 0)$ can be identified from the proof of Proposition \ref{prop:fano-three-wall-bound} while the case $(1, 2, 0)$ is identified by  \url{wall_est.m}.  All the triples listed for $K^2=7, 8$ in Table \ref{tab:wall-est-maximizers} are identified by \url{wall_est.m}.

\begin{table}[htbp]
\centering
\small
\caption{For each value of $K^2$, we show the corresponding value of  $n=9-K^2$, an upper bound  $t_0^2$ and the corresponding $t_0$ such that $T_X$ is $\sigma_{-tK}$-semistable for all $t>t_0$.  We also list all possible cases of $(\ch_0(\image \beta), -Kc_1(\image \beta), -Kc_1(Q^{-1}))$ that give rise to the listed upper bound $t_0$.  The `Remark' column shows the result whose proof contains the argument identifying the listed case, unless that case is identified by the Matlab script.}
\label{tab:wall-est-maximizers}
\begin{tabular}{cccccccc}
\hline
$n$ & $K^2$ & upp.\ bound $t_0^2$ & upp.\ bound $t_0$ &
$\ch_0(\image \beta)$ & $-Kc_1(\image \beta)$ & $-Kc_1(Q^{-1})$ & Remark \\
\hline
 0&9 &0 &0 & & & & Theorem \ref{thm:AG54-53-1}\\ 
 1&8 & $\frac12=0.5$& $\frac{1}{\sqrt{2}} \approx 0.707$ &1 &3 &0 & \url{wall_est.m} \\ 
 & & & & 1& 4& $-1$& \url{wall_est.m}\\ 
 & & & & 1& 4& $-2$& \url{wall_est.m}\\ 
 & & & & 1& 4& $-3$& \url{wall_est.m}\\ 
 & & & & 2& 6& 0& \url{wall_est.m}\\ 
 & & & & 2& 7& 0& \url{wall_est.m}\\ 
 2& 7& $\frac{16}{7}\approx 2.286$ & $\frac4{\sqrt{7}} \approx 1.512$ &1 &3 &0  & \url{wall_est.m}\\ 
 3& 6& 1& 1& 1& 2& 0& \url{wall_est.m}\\ 
 & & & & 2& 4& 0& Proposition \ref{prop:fano-three-wall-bound}\\ 
 4& 5& $\frac{11}{10}=1.1$& $\sqrt{\frac{11}{10}}\approx 1.049$ & 2& 3& 0&Proposition \ref{prop:fano-four-wall-bound}\\ 
 5& 4& 1& 1& 1& 1& 0& Proposition \ref{pro:AG54-133-1} \\ 
 & & & & 2& 2& 0& Proposition \ref{pro:AG54-133-1} \\ 
 6& 3& $\frac34$ &  $\sqrt{\frac34}\approx 0.886$  &2 &1 &0 &  Proposition \ref{pro:AG54-104-1}\\ 
 7& 2& 0&0 & & & & Theorem \ref{thm:TXn-stab-78}\\ 
 8&1 &0 &0 & & & & Theorem \ref{thm:TXn-stab-78}\\ 
\hline
\end{tabular}
\end{table}

In the rest of this article, we focus on ruling out  walls listed in  Table \ref{tab:wall-est-maximizers} or identified by \url{wall_est.m}.

We end this section with another technical lemma that will be used in ruling out certain walls.

\begin{lem}\label{lem:AG55-59-1}
Suppose $X$ is a smooth projective surface, $\omega$ an ample class on $X$, $E$ a torsion-free sheaf in $\mathcal{T}_\omega$  and $0 \to A \to E \to Q \to 0$  a destabilizing sequence in $\mathcal{B}_\omega$ corresponding to the outermost mini-wall  as in \ref{para:outermostwall-existence}.  If $\ch_2(Q^{-1}), \ch_2(\image \beta) \leq 0$ then $Q^{-1}=0$.
\end{lem}

\begin{proof}
Suppose $Q^{-1} \neq 0$.  By Lemma \ref{lem:AG54-112-1}, we must have $\omega c_1(Q^{-1}) < 0$.  Since $A$ is a nonzero torsion-free sheaf with $\omega c_1>0$ and $\ch_2(A) = \ch_2(Q^{-1})+\ch_2(\image \beta) \leq 0$, we have 
\[
\rho_\omega (A) = \frac{\ch_2(A) - \tfrac{\omega^2}{2}\ch_0(A)}{\omega c_1(A)} < 0
\]
while 
\[
\rho_\omega (Q^{-1}[1]) = \frac{\ch_2(Q^{-1}[1])- \tfrac{\omega^2}{2}\ch_0(Q^{-1}[1])}{\omega c_1(Q^{-1}[1])} \geq 0.
\]
Now, the short exact sequence in $\Coh (X)$
\[
0 \to Q^{-1} \to A \to \image \beta \to 0
\]
gives an exact triangle in $D^b(X)$
\[
A \to \image \beta \to Q^{-1}[1] \to A[1].
\]
Since $\rho_\omega (A) < \rho_\omega (Q^{-1}[1])$ from above, we obtain $\rho_\omega (A) < \rho_\omega (\image \beta)$.  

Note also that the short exact sequence 
\[
0 \to \image \beta \to E \to Q^0 \to 0
\]
is also a short exact sequence in $\Bc_\omega$ and so $\image \beta$ is also a $\Bc_\omega$-subobject of $E$.  This contradicts the assumption that $A$ is a destabilizing subobject of $E$ at the outermost wall.  Hence $Q^{-1}$ must vanish.
\end{proof}

\section{Precise bounds for the outermost walls} 
In this section, we determine the sharp stability bound for the tangent bundle \(T_{X_n}\) for \(0\leq n\leq 8\). The result is summarized in the following Theorem.

\begin{thm}
\label{thm:mainthm}
Let $X=X_n$ where $0 \leq n \leq 8$, and let $t_0$ be given in Table \ref{tab:outermostwall}. Then $T_{X}$ is $\sigma_{-tK}$-semistable for $t >t_0$ and unstable for $t<t_0$. Furthermore, if $n\neq 1$, then $T_{X}$ is $\sigma_{-tK}$-stable for $t>t_0$, whereas if $n=1$, it is strictly $\sigma_{-tK}$-semistable for all $t\geq t_0$.

\begin{table}[htbp]
\small
\centering
\caption{The precise location of the outermost wall for $T_{X_n}$ along the family $\sigma_{-tK_{X_n}}$.}
\label{tab:outermostwall}
\begin{tabular}{cccl}
\hline
$n$ & $K^2$ & $t_0$ & Remark \\
\hline
0 & 9& 0& Theorem \ref{thm:AG54-53-1}\\ 
 1& 8& $\sqrt{\frac12}$ & Remark \ref{rem:n=2}\\ 
 2& 7& $\sqrt{\frac27}$ & Corollary \ref{cor:X2stable}\\ 
 3& 6& 0& Corollary \ref{cor:X3stable}\\ 
 4& 5& 0& Corollary \ref{cor:X4stable}\\ 
 5& 4& 0& Corollary \ref{cor:X5stable}\\ 
 6& 3& 0& Corollary \ref{cor:X6stable}\\ 
 7& 2& 0& Theorem \ref{thm:TXn-stab-78}\\ 
 8& 1& 0& Theorem  \ref{thm:TXn-stab-78}\\ 
\hline
\end{tabular}
\end{table}

\end{thm}
\begin{proof}
We prove the theorem for each $n$ in the rest of this section, with the corresponding results summarized in the ``Remark'' column. For $n=1$, we record some destabilizing short exact sequences at $t=t_0$  in Proposition \ref{prop:X1notstable} and Proposition \ref{prop:X1rank2wall}. For $n=2$, the destabilizing short exact sequences are presented in Proposition \ref{cor:X2stable}.
\end{proof}
Following the notation of the previous section, we denote the pullback of the hyperplane class by $l$ and the exceptional divisors by $e_i$.  In this section, we denote the canonical class of $X_n$ by $K_n$; when $n$ is understood, we simply write $X, K$ to denote $X_n, K_n$, respectively. The slope of $T_{X_n}$ is given by 
\begin{equation}
\label{eq:slopeTX}
\rho_{-tK_n}(T_{X_n})=\frac{\frac{3}{2}(1-n)}{t(9-n)}-t.
\end{equation}
The case $n=1$ is special, as $T_{X_1}$ can be $\sigma_{-tK_1}$-semistable but never stable. We treat this case first.

\begin{prop}
\label{prop:X1notstable}
(i) There exists a short exact sequence in $\mathcal{B}_{-K_1}$
\[
0\to \cO_{X_1}(l)\to T_{X_1}\to Q\to 0,
\]
where $\rho_{-tK_1}(\cO_{X_1}(l))> \rho_{-tK_1}(T_{X_1})$ when $t< \sqrt\frac12$.

(ii) There exists a short exact sequence in $\mathcal{B}_{-K_1}$
\[
0\to \cO_{X_1}(l+e_1)\to T_{X_1}\to Q\to 0,
\]
where $\rho_{-tK_1}(\cO_{X_1}(l+e_1))= \rho_{-tK_1}(T_{X_1})$ for all $t>0$.
\end{prop}
\begin{proof} 
We first verify the existence of the morphisms $\cO_{X_1}(l)\to T_{X_1}$ and $\cO_{X_1}(l+e_1)\to T_{X_1}$. Let $f:X_1\to X_0$ be the blow-up morphism. Recall that there exists a short exact sequence
\[
0\to T_{X_1}\to f^*T_{X_0}\to \cO_{e_1}(1)\to 0.
\]
Let $a, b\in \mathbb{Z}$, then
\[
\begin{split}
\Hom(\cO_{X_1}(al-be_1), f^*T_{X_0})&\simeq H^0(f^*T_{X_0}(-a)\otimes \cO_{X_1}(be_1))\\
&\simeq H^0(f_*(f^*T_{X_0}(-a)\otimes \cO_{X_1}(be_1)))\\
&\simeq H^0(T_{X_0}(-a)\otimes f_*\cO_{X_1}(be_1))
\end{split}
\]

Claim. If $a\leq 1$ and $b\geq 0$, then $\Hom(\cO_{X_1}(al-be_1), T_{X_1})\neq 0$. 

First, observe that $f_*\cO_{X_1}(be_1)\simeq \cO_{X_0}$ for all $b\geq 0$. When $b=0$, this follows from Zariski's Main Theorem. For $b>0$, the result follows by induction using the following short exact sequence: 
\[
0\to \cO_{X_1}((b-1)e_1)\to \cO_{X_1}(be_1)\to \cO_{e_1}(-b)\to 0.
\]
When $b\geq 0$, we have $\Hom(\cO_{X_1}(al-be_1), f^*T_{X_0})\simeq H^0(T_{X_0}(-a))\simeq \mathbb{C}^n$. When $a\leq 1$ we have $n\geq 3$. 
On the other hand, $\Hom(\cO_{X_1}(al-be_1), \cO_{e_1}(1))\simeq H^0(\cO_{e_1}(1-b))\simeq \mathbb{C}^m$.  When $b>1$, we have $m=0$; when $b=1$, we have $m=1$; and when $b=0$, we have $m=2$. In each case, $n-m>0$, which implies that $\Hom(\cO_{X_1}(al-be_1),T_{X_1})\neq 0$.

Now consider $a=1$, $b=-1$. Identity $X_1$ with $\mathbb{F}_1$, the first Hirzebruch surface.  We have the relative tangent sequence
\[
0\to T_{X_1/\mathbb{P}^1}\to T_{X_1}\to \pi^*T_{\mathbb{P}^1}\to 0.
\]
Then $T_{X_1/\mathbb{P}^1}\simeq \cO_{X_1}(-K_1-2F)$ where $F=l-e_1$.  Hence $T_{X_1/\mathbb{P}^1}\simeq \cO_{X_1}(l+e_1)$, which gives the required morphism. Both $\cO_{X_1}(l)$ and $\cO_{X_1}(l+e_1)$ lie in $\Bc_{-K_1}$ by computing the $\mu$-slopes. Since $T_{X_1}$ is $\mu$-semistable, $T_{X_1}$ and its quotient in $\Coh(X)$ also lies in $\Bc_{-K_1}$. Hence the short exact sequences in (i) and (ii) are both short exact sequences in $\Bc_{-K_1}$.

Now we compare the Bridgeland slope of the line bundles $\cO_{X_1}(al-be_1)$ and $T_{X_1}$, we have 
\[
\rho_{-tK_1}(\cO_{X_1}(al-be_1))=\frac{\frac{1}{2}(a^2-b^2)-\frac{1}{2}t^2K_1^2}{-tK_1\cdot(al-be_1)}=\frac{\frac{1}{2}(a^2-b^2)-4t^2}{t(3a-b)},
\]
\[
\rho_{-K_1}(T_{X_1})=-t.
\]
Then $\rho_{-tK_1}(\cO_{X_1}(al-be_1))> \rho_{-tK_1}(T_{X_1})$  iff 
\begin{equation}
\label{eq:X1slope}
t^2(4-3a+b)< \frac{1}{2}(a^2-b^2).
\end{equation}
When $a=1$, $b=-1$, $\rho_{-K_1}(\cO_{X_1}(l+e_1))=\rho_{-K_1}(T_{X_1})$ for all $t>0$. When $a=1$, $b=0$, the inequality \eqref{eq:X1slope} is equivalent to $t<\sqrt\frac12$. 
\end{proof}

There also exists a rank-two object that destabilizes $T_{X_1}$. 
\begin{prop}
\label{prop:X1rank2wall}
There exists a short exact sequence in $\mathcal{B}_{-K_1}$
\[
0\to M\to T_{X_1}\to \cO_{e_1}(-1)\to 0,
\]
where $\rho_{-tK_1}(M)> \rho_{-tK_1}(T_{X_1})$ when $t< \sqrt\frac12$.
\end{prop}
\begin{proof}
Recall that $T_{X_1}|_{e_1}\simeq \cO_{e_1}(-1)\oplus\cO_{e_1}(2)$. Then $M$ fits into the short exact sequence
\[
0\to T_X(-e_1)\to M\to O_{e_1}(2)\to 0.
\]
Since $T_{X_1}$ is $\mu_{-K_1}$-semistable, and $\mu_{-K_1}(T_{X_1}(-e_1))=\frac{6}{2}>0$, we have $T_{X_1}(-e_1)\in \Bc_{-K_1}$. This implies that $M\in \Bc_{-K_1}$.
Comparing the Bridgeland slopes of $M$ and $T_{X_1}$, we obtain the conclusion. 
\end{proof}
From the MATLAB output for $X_1$ in Table \ref{tab:wall-est-positive-n1}, we are left with one case to consider.
\begin{lem}\label{lem:Xeq1-specialcase}
Let $X=X_1$. Then the outermost wall in the sense of \ref{para:outermostwall-existence} exists in the family $\{ \sigma_{-tK}\}_{t>0}$.   A destabilizing sequence $0 \to A \to T_X \to Q \to 0$ in $\mathcal{B}_{-K}$ at the outermost wall, however, cannot have  
\begin{equation}\label{eq:Xeq1-eq1}
(\ch_0(\image \beta), -Kc_1(\image \beta), -Kc_1(Q^{-1})) = (1, 4, 0)
\end{equation}
where $\beta$ denotes the morphism $A \to T_X$ in the long exact sequence in $\Coh (X)$.
\end{lem}

\begin{proof}
By Proposition \ref{prop:X1notstable}(i), $T_X$ is not $\sigma_{-tK}$-semistable for all $t>0$.  Then the existence of the outermost wall follows from Lemmas \ref{lem:X1Txsslarget} and \ref{lem:goingup-v3}.

Now suppose $0 \to A \to T_X \to Q \to 0$ is a destabilizing sequence in $\mathcal{B}_{-K}$ at the outermost wall $t$, and that we are in the case  \eqref{eq:Xeq1-eq1}.  By Lemma \ref{lem:AG54-112-1}, we have $Q^{-1}=0$ and so $A$ is a rank-one subsheaf of $T_X$ with $-Kc_1(A)=4$.  Lemma \ref{lem:AG54-105-1} tells us that $A$ is locally free, and so dualizing gives the exact sequence of sheaves
\[
0 \to Q^\ast \to \Omega_X \overset{\alpha}{\to} A^\ast.
\]
Since $-Kc_1(A)=4$, it follows that $-Kc_1(A^\ast)=-4$ and $-Kc_1(\image \alpha) \leq -4$, which in turn implies $-Kc_1(Q^\ast) \geq -4$.  By the $\mu_{-K}$-semistability of $\Omega_X$, we get $-Kc_1(Q^\ast)=-4$ and $-Kc_1(\mathrm{coker}\, \alpha)=0$, which in turn implies $Q$ is torsion-free and so $Q^\ast$ is a line bundle.  Since $Q^\ast$ is also a saturated subsheaf of $\Omega_X$ of degree $-Kc_1(Q^\ast)=-4$, from \cite[Section 3]{fahlaoui1989stabilite} we know $Q^\ast \cong \OO_X (-2l+2e_1)$ and thus $c_1(Q)=2l-2e_1$, giving us $c_1(A)=-K-(2l-2e_1)=l+e_1$.  That is, $A \cong \OO_X (l+e_1)$.  This contradicts Proposition \ref{prop:X1notstable}(ii), however, since the condition of $t$ being an outermost wall in the sense of \ref{para:outermostwall-existence} requires that $\rho_{-tK}(A) > \rho_{-tK}(T_X)$ for values of $t$ close to 0.
\end{proof}
\begin{rem}
\label{rem:n=2}
Proposition \ref{prop:X1notstable} shows that $T_{X_1}$ is not $\sigma_{-tK_{X_1}}$-stable for any $t>0$. Combined with the MATLAB computations in Table \ref{tab:wall-est-maximizers} and Lemma \ref{lem:Xeq1-specialcase}, this implies that $T_{X_1}$ is strictly semistable for $t\geq \sqrt{\frac12}$. Finally, Propositions \ref{prop:X1notstable}(i) and \ref{prop:X1rank2wall} imply that it is unstable for $t<\sqrt{\frac12}$. 
\end{rem}

Since $T_{X_n}$ is stable for all $t>0$ when $n=7,8$, and we have already determined the sharp semistability region for $n=1$, we will consider only $2\leq n\leq 6$ for the rest of this section. Recall that in these cases, the effective monoid. $\operatorname{Eff}(X_n)$ is generated by classes of negative rational curves.
We begin with several technical lemmas analyzing potential destabilizing subsheaves of $T_{X_n}$ of rank two.

Let $C_1$ and $C_2$ be exceptional curves on $X_n$ such that $C_1\cdot C_2=1$. By the Grothendieck–Riemann–Roch, we have $\Ext^1(\cO_{C_2}(k_2), \cO_{C_1}(k_1))\simeq \mathbb{C}$, for any $k_1$, $k_2$. We denote the unique nontrivial extension by $E_{k_1,k_2}$, which fits into the short exact sequence 
\[
0\to \cO_{C_1}(k_1)\to E_{k_1, k_2}\to \cO_{C_2}(k_2)\to 0.
\]

\begin{lem}
\label{lem:norcrossquot}
Let $C_1$ and $C_2$ be exceptional curves on $X_n$ such that $C_1\cdot C_2=1$. Then 
\[
\Hom(T_{X_n}, E_{-2, -1})=0.
\]
\end{lem}
\begin{proof}
The argument does not depend on $n$ so we omit $n$ in the notation and write $X$ for $X_n$.
Let $D:=\text{supp}(E_{-2, -1})=C_1\cup C_2$. We show that there does not exist a map $T_{X}|_D\to E_{-2, -1}$ in $\Coh(X)$.
Recall that for any exceptional curve $C$, $T_{X}|_{C}\simeq \cO_{C}(2)\oplus \cO_C(-1)$,
hence $\Hom(T_{X}|_{C_1}, \cO_{C_1}(-2))=0$.
Since $\Hom (T_X|_D, \OO_{C_2}(-1))$ is one-dimensional, to prove the claim, we only need to show the connecting homomorphism 
\[
\Hom(T_{X}|_D, \cO_{C_2}(-1))\xrightarrow{g}\Ext^1(T_{X}|_D, \cO_{C_1}(-2))
\]
is not the zero map. This forces the only map $T_{X}\to E_{-2, -1}$ to be the zero map. 
Applying the connecting homomorphism to the morphism $T_{X}|_D\to T_{X}|_{C_2}$, we obtain a diagram 
\[
\begin{tikzcd}
& 0\ar{d}\\
&\Hom(T_{X}\otimes \cO_{C_1}(-1), \cO_{C_1}(-2))\ar{d}{l}\\
\Hom(T_{X}|_{C_2}, \cO_{C_2}(-1))\ar{r}{f}\ar{d}{e} & \Ext^1(T_{X}|_{C_2}, \cO_{C_1}(-2))\ar{d}{h}\\
\Hom(T_{X}|_{D}, \cO_{C_2}(-1))\ar{r}{g} & \Ext^1(T_{X}|_{D}, \cO_{C_1}(-2))
\end{tikzcd}
\]
The columns are obtained by applying $\Hom(-, \cO_{C_2}(-1))$ and $\Hom(-, \cO_{C_1}(-2))$ to the short exact sequence
\begin{equation}
\label{eq:caseiii}
0\to T_{X}\otimes \cO_{C_1}(-1)\to T_{X}|_D\to T_{X}|_{C_2}\to 0.
\end{equation}
Since $e$ is an isomorphism, it is enough to show the map $h\circ f$ is not the zero map. Note that we have $\Hom(T_{X}|_{C_2}, \cO_{C_2}(-1))\simeq \mathbb{C}$ and $\Hom(T_{X}\otimes \cO_{C_1}(-1), \cO_{C_1}(-2))\simeq \mathbb{C}$. We know that 
\[
\Ext^1(T_X|_{C_2}, \cO_{C_1}(-2))\simeq \Ext^1(\cO_{C_2}(2), \cO_{C_1}(-2))\oplus \Ext^1(\cO_{C_2}(-1), \cO_{C_1}(-2)).
\]
Tracing through the definition of the connecting homomorphism, we have 
\[
\text{im}\, f \subseteq \Ext^1(\cO_{C_2}(-1), \cO_{C_1}(-2))\subseteq \Ext^1(T_{X}|_{C_2}, \cO_{C_1}(-2)).
\]

On the other hand, using local coordinates of $\Spec \mathbb{C}[x, y]$, and set $C_1$ be defined by $y=0$ and $C_2$ defined by $x=0$. Applying $T_X\otimes-$ to $0\to \cO_{C_1}(-1)\to \cO_D\to \cO_{C_2}\to 0$ results in
\[
0\to \frac{\partial}{\partial x}(x)/(xy)\oplus \frac{\partial}{\partial y}(x)/(xy)\to \frac{\partial}{\partial x}\mathbb{C}[x, y]/(xy)\oplus \frac{\partial}{\partial y}\mathbb{C}[x, y]/(xy)\to \frac{\partial}{\partial x}\mathbb{C}[x, y]/(x)\oplus \frac{\partial}{\partial y}\mathbb{C}[x, y]/(x)\to 0
\]
where $\frac{\partial}{\partial x}$ is the local section of $\cO_{C_1}(1)$ and $\cO_{C_2}(-1)$, $\frac{\partial}{\partial y}$ is the local section of
$\cO_{C_1}(-2)$ and $\cO_{C_2}(2)$.

Hence
\[
\text{im}\, l\subseteq \Ext^1(\cO_{C_2}(2), \cO_{C_1}(-2))\subseteq \Ext^1(T_{X}|_{C_2}, \cO_{C_1}(-2)).
\]
Since $\ker h= \text{im}\, l$, we have $h\circ f\neq 0$.

\end{proof}

We fix some notation for the remainder of this section. Let $Q$ be a pure one-dimensional sheaf, and set
\[
D:=\operatorname{Supp}(Q),
\]
where $\operatorname{Supp}(Q)$ is the subscheme defined by the annihilator ideal sheaf $\operatorname{Ann}(Q)$.
Since $X_n$ is smooth, $D$ is an effective Cartier divisor, and any morphism $T_{X_n}\to Q$ factors through the restriction $T_{X_n}|_D\to Q$. Moreover, $c_1(Q)-c_1(\cO_D)$ is zero or is an effective divisor. When no confusion arises, we also use $D$ to denote its divisor class.

\begin{lem}
\label{lem:destabnonred}
Let $X=X_n$ where $2\leq n\leq 6$. Let $0\to A\to  T_{X}\to Q\to 0$ be a short exact sequence in $\Coh(X)$ such that $Q$ is a pure $1$-dimensional sheaf, $c_1(Q)=mC$ where $C$ is the class of an exceptional curve on $X$. Then $\rho_{-tK}(Q)\geq -\frac{1}{2t}$.
\end{lem}

\begin{proof}
We have the following short exact sequences 
\[
0\to Q_{i+1}\to Q_i\to Q_i|_{C}\to 0
\]
where $Q_0=Q$, $Q_{i+1}\simeq I_CQ_i$ and $Q_i\simeq I_C^iQ$. 
Since $D=bC$ for some $1\leq b\leq m$, there is a natural surjection $\cO_{D}\twoheadrightarrow \cO_C$. 
Hence, we obtain the following short exact sequences:
\[
0\to M_{i+1}\to M_i\to M_i|_{C}\to 0,
\]
where $M_0=\cO_D$, $M_{i+1}=I_CM_i$, and $M_i\simeq I_C^i\cO_D$. Since $T_X\otimes \cO_D\twoheadrightarrow Q$, we obtain 
\[
T_X\otimes M_i\twoheadrightarrow Q\otimes M_i\twoheadrightarrow I_C^iQ=Q_i.
\]
Then we have the following commutative diagram:
\[
\begin{tikzcd}
0\arrow[r] & T_X\otimes M_{i+1}\arrow[r]\arrow[d] & T_{X}\otimes M_i\arrow[r]\arrow[d] &T_X\otimes M_i|_{C}\arrow[r]\arrow[d] & 0\\
0\arrow[r] & Q_{i+1} \arrow[r]\arrow[d]&Q_i  \arrow[r]\arrow[d] &Q_i|_{C}\arrow[r]\arrow[d]  & 0\\
& 0 & 0&0 &
\end{tikzcd}
\]
Since $T_X(-iC)\twoheadrightarrow T_X\otimes M_i$ we have 
\[
\ch_2(Q_i|_C)\geq \text{min}\{\ch_2(\mathcal{V})|T_X(-iC)|_{C}\twoheadrightarrow \mathcal{V}\}.
\]
Note that $T_X(-iC)|_C\simeq \cO_C(i+2)\oplus \cO_C(i-1)$. Since $i\geq 0$, we have $\ch_2(Q_i|_C)\geq -\frac{1}{2}$, hence  $\rho_{-tK}(Q)\geq -\frac{1}{2t}$.  
\end{proof}

The following Proposition will be needed when we consider the destabilizing short exact sequence when $n=5$ and $n=6$. 

\begin{prop}
\label{prop:desesdeg2}
Let $X=X_n$ where $2\leq n\leq 6$. Let $0\to A\to T_{X}\to Q\to 0$ be a short exact sequence in $\Coh(X)\cap \Bc_{-K}$ such that $Q$ is a pure 1-dimensional sheaf. Assume that $-K\cdot c_1(Q)= 2$, then $\rho_{-tK}(Q)\geq -\frac{1}{2t}$.
\end{prop}
\begin{proof}

Since $-K\cdot c_1(Q)=2$, we have $c_1(Q)= C_1+C_2$ where $C_1$ and $C_2$ are the divisor classes of two exceptional curves. By Lemma \ref{lem:negcurverig}, $C_i$ is the only curve in the class, so we abuse notation and use $C_i$ to denote the subvariety if there is no confusion. 
There are three cases to consider based on the intersection number between $C_1$ and $C_2$, 

(i) $C_1= C_2$, 

(ii) $C_1\neq C_2$ and $C_1\cdot C_2=0$, 

(iii) $C_1\neq C_2$ and $C_1\cdot C_2=1$. 

 Case (i) following from Lemma \ref{lem:destabnonred}.
For case (ii), $C_1$ and $C_2$ must be disjoint and so $T_X|_D\simeq T_X|_{C_1}\oplus T_X|_{C_2}$. Hence the minimal slope of $Q$ is obtained at $Q\simeq \cO_{C_1}(-1)\oplus \cO_{C_2}(-1)$ or $Q\simeq \cO_{C_1}(-1)$ or $Q\simeq \cO_{C_2}(-1)$. In either case $\rho(Q)=-\frac{1}{2t}$.

For case (iii), we have $c_1(\cO_D)=c_1(Q)=C_1+C_2$. 
Consider the short exact sequences 
\[
\begin{split}
&0\to \cO_X\to \cO_X(C_1)\to \cO_{C_1}(-1)\to 0,\\
&0\to \cO_X(C_1)\to \cO_X(C_1+C_2)\to \cO_{C_2}\to 0.
\end{split}
\]
Then we have $H^0(\cO_X(C_1))\simeq \mathbb{C}$, $H^0(\cO_{C_2})\simeq \mathbb{C}$. Moreover, $H^1(\cO_X)=0$ implies that $H^1(\cO_X(C_1))=0$. Hence $H^0(\cO_X(C_1+C_2))\simeq \mathbb{C}^2$, i.e. the linear system $|C_1+C_2|$ is $1$-dimensional. By \cite[Theorem III.1(a)]{Har97}, this linear system is base point free. Then Bertini's Theorem implies that the general element is smooth. 
Furthermore, consider the short exact sequence 
\[
0\to \cO_X(-C_1-C_2)\to \cO_X\to \cO_D\to 0.
\]
Since $H^1(\cO_X(-C_1-C_2))=0$, we have $H^0(\cO_D)\simeq \mathbb{C}$, hence $D$ is connected. 
Since 
\[
K_D\simeq (D+K_X)D=-2,
\]
which implies that the geometric genus of the generic smooth element of $|C_1+C_2|$ is $0$, hence the generic smooth element of $|C_1+C_2|$ is isomorphic to $\mathbb{P}^1$. Since $p_a=0$ for any element of the linear system, the singular element of $|C_1+C_2|$ must be reducible, and hence $D=D'\cup D''$, since $-K\cdot D=2$, this forces $-K\cdot D'=1$, $-K\cdot D''=1$. Hence $D'$ and $D''$ are exceptional curves with $D'\cdot D''=1$.

We have $D\in |C_1+C_2|$. If $D$ is the generic smooth $\mathbb{P}^1$ in $|C_1+C_2|$, then 
\[
T_X|_D\simeq \cO_D(2)\oplus \cO_D.
\]

Thus $Q\simeq\mathcal O_D(k)$ for some $k\ge 0$, and $\ch_2(Q)=k$. This implies that $\rho_{-K}(Q)\geq 0$, with equality achieved at $Q\simeq \cO_{D}$. We next consider when $D=D'\cup D''$; without loss of generality, we may assume $D=C_1\cup C_2$. Tensoring the short exact sequence
\[
0\to \cO_{C_1}(-1)\to \cO_D\to \cO_{C_2}\to 0
\]
with $T_X$, we have the following diagram
\[
\begin{tikzcd}
0\arrow[r] & T_X|_{C_1}\otimes \cO_{C_1}(-1)\arrow[r]\arrow[d] & T_{X}|_D\arrow[r]\arrow[d] &T_X|_{C_2}\arrow[r]\arrow[d] & 0\\
0\arrow[r] & Q_1 \arrow[r]\arrow[d] &Q  \arrow[r]\arrow[d] &Q|_{C_2}\arrow[r]\arrow[d] & 0\\
& 0 & 0&0 &
\end{tikzcd}
\]
The left column is the composition of 
\[
T_X|_{C_1}\otimes \cO_{C_1}\twoheadrightarrow Q\otimes \cO_{C_1}(-1)\twoheadrightarrow \cO_{C_1}(-1)Q\simeq Q_1,
\]
hence it is surjective. 
Here $T_X|_{C_2}\simeq \cO_{C_2}(2)\oplus \cO_{C_2}(-1)$ and $T_X|_{C_1}\otimes \cO_{C_1}(-1)\simeq \cO_{C_1}(1)\oplus \cO_{C_1}(-2)$. Since $Q$ is pure, $Q_1$ is pure and $Q_1\simeq \cO_{C_1}(k_1)$ where $k_1=-2$ or $k_1\geq 1$. Moreover, $Q|_{C_2}\simeq P\oplus \cO_{C_2}(k_2)$, where $k_2=-1$ or $k_2\geq 2$, and $P$ is a sheaf supported in dimension $0$ if nonzero. The only possibility for $\rho_{-tK}(Q)<-\frac{1}{2t}$  happens when $Q_1\simeq \cO_{C_1}(-2)$, and $Q|_{C_2}\simeq \cO_{C_2}(-1)$. By Lemma \ref{lem:norcrossquot}, $\Hom(T_X|_D\to E_{-2, -1})=0$, hence this cannot happen and $\rho_{-tK}(Q)\geq -\frac{1}{2t}$.
\end{proof}

We also analyze the cases when the degree of $c_1(Q)$ equals $3$ or $4$. The following Proposition will be used when we consider the destabilizing short exact sequence for $n=5$.  

\begin{prop}
\label{prop:desesdeg3}
Let $X=X_n$, where $2\leq n\leq 6$.
Let $0\to A\to T_{X}\to Q\to 0$ be a short exact sequence in $\Coh(X)\cap \mathcal{B}_{-K}$, where $Q$ is a pure $1$-dimensional sheaf. Assume that $-K\cdot c_1(Q)=3$, then $\rho_{-tK}(Q)\geq -\frac{1}{2t}$. 
\end{prop}
\begin{proof}
 We can assume $D$ is connected, for if not, we can analyze each component separately and reduced to cases where $c_1(Q)$ has smaller degree. 

We have $c_1(Q)= C_1+C_2+C_3$ where $C_i$'s are the divisor classes of exceptional curves. There are four possible cases:

(i) $C_1=C_2=C_3$; 

(ii) $C_1= C_2$, $C_1\neq C_3$ and $C_1\cdot C_3=1$; 

(iii) $C_i\neq C_j$ for $i\neq j$ and $C_1\cdot C_3=C_2\cdot C_3=1$, $C_1\cdot C_2=0$.

(iv) $C_i\neq C_j$ for $i\neq j$ and $C_i\cdot C_j=1$ for all $i\neq j$.

Case (i) follows from Lemma \ref{lem:destabnonred}. We consider case (ii). In this case, we have $D= 2C_1+C_3$ or $D= C_1+C_3$, i.e. $D$ has two irreducible components with multiplicity given by the coefficients.  Define $D_0$ to be the reduced support with class $C_1+C_3$. For 
$D=2C_1+C_3$, consider the following diagram:
\[
\begin{tikzcd}
&&&0\arrow[d]&\\
&0\arrow[d]&&\cO_{C_1}(-C_1-C_3)\arrow[d]&\\
0\arrow[r]&\cO_X(-2C_1-C_3)\arrow[r]\arrow[d]&\cO_X\arrow[r]\arrow[d, "\simeq"]&\cO_D \arrow[r]\arrow[d]&0\\
0\arrow[r]&\cO_X(-C_1-C_3)\arrow[r]\arrow[d]&\cO_X\arrow[r]&\cO_{D_0}\arrow[r]\arrow[d]&0\\
&\cO_{C_1}(-C_1-C_3)\arrow[d]&&0&\\
&0&&&
\end{tikzcd}
\]

Consider the following short exact sequence
\[
0\to Q_1\to Q\to Q|_{D_0}\to 0
\]
Then we have $c_1(Q_1)= C_1$ and $\text{supp}(Q_1)=C_1$, we have the following surjection
\[
T_X(-C_1-C_3)|_{C_1}\to Q_1\to 0.
\]
 Since $T_X(-C_1-C_3)|_{C_1}\simeq T_X|_{C_1}$, we have $\rho_{-K_X}(Q_1)\geq -\frac{1}{2t}$. By Proposition \ref{prop:desesdeg2} case (iii), $\rho_{-K}(Q|_{D_0})\geq -\frac{1}{2t}$, hence $\rho_{-K_X}(Q)\geq -\frac{1}{2t}$. For 
 $D=C_1+C_3$, 
A similar argument as in Proposition \ref{prop:desesdeg2} gives the following diagram
\[
\begin{tikzcd}
0\arrow[r] & T_X|_{C_1}\otimes \cO_{C_1}(-1)\arrow[r]\arrow[d] & T_{X}|_D\arrow[r]\arrow[d] &T_X|_{C_3}\arrow[r]\arrow[d] & 0\\
0\arrow[r] & Q_2 \arrow[r] &Q  \arrow[r] &Q|_{C_3}\arrow[r] & 0\\
\end{tikzcd}
\]
The coefficient of $C_1$ in $c_1(Q)$ is two, so $Q_2$ is generically of rank two on $C_1$. Since 
\[T_X|_{C_1}\otimes \cO_{C_1}(-1)\simeq \cO_{C_1}(1)\oplus \cO_{C_1}(-2),
\]
and $T_X|_{C_1}\otimes \cO_{C_1}(-1)\to Q$ is surjective, hence $Q_2\simeq \cO_{C_1}(1)\oplus \cO_{C_1}(-2)$.
Then $\ch_2(Q_2)=0$.  The proof of Lemma \ref{lem:destabnonred} implies that $\ch_2(Q_3)\geq -\frac{1}{2}$. Hence $\ch_2(Q)\geq -\frac{1}{2}$ and $\rho_{-tK}(Q)\geq -\frac{1}{2t}$.

For Case (iii) and (iv), $D$ is reduced. We consider Case (iii).
Consider the short exact sequences 
\[
\begin{split}
&0\to \cO_X(C_2+C_3)\to \cO_X(C_1+C_2+C_3)\to \cO_{C_1}\to 0\\
&0\to \cO_X(C_3)\to \cO_X(C_2+C_3)\to \cO_{C_2}\to 0\\
&0\to \cO_X\to \cO_X(C_3)\to \cO_{C_3}(-1)\to 0.
\end{split}
\]
By a similar argument as in Proposition \ref{prop:desesdeg2}, we have $H^0(\cO_X(C_2+C_3))\simeq \mathbb{C}^2$, which further implies that $H^0(\cO_X(C_1+C_2+C_3))\simeq \mathbb{C}^3$, hence $|C_1+C_2+C_3|$ is $2$ dimensional. Furthermore, the generic smooth element in $|C_1+C_2+C_3|$ is isomorphic to $\mathbb{P}^1$ with self intersection $1$. Also, the singular elements are reducible. If there are two irreducible components, we can write $D=D'\cup D''$, where $D'$ is an exceptional curve, and $-K\cdot D''=2$. Then using the analysis of Proposition \ref{prop:desesdeg2}, we have $D''\simeq \mathbb{P}^1$ in a linear system $|E_1+E_2|$ for some exceptional curves $E_1$ and $E_2$. Since $2\leq n\leq 6$, $E_1\cdot E_2=0$ or $1$. If $E_1\cdot E_2=0$, $D''$ is disconnected, contradicting it is irreducible. Hence $E_1\cdot E_2=1$, which implies that $(D'')^2=0$. Finally, if $D$ has three irreducible components, we can write $D=D_1\cup D_2\cup D_3$, where $D_i's$ are exceptional curves. Since $D^2=1$, we must have $D_i\cdot D_j=0$ for some $i, j$, and other pairwise intersection numbers are $1$.

The smooth and two irreducible component cases follow from  similar arguments as in the proof of Proposition \ref{prop:desesdeg2}. We only need to consider the three components case $D=D_1\cup D_2\cup D_3$; without loss of generality, we can assume $D=C_1\cup C_2\cup C_3$.
Apply $T_X\otimes -$ to 
\[
0\to \cO_{C_1}(-1)\oplus \cO_{C_2}(-1)\to \cO_D\to \cO_{C_3}\to 0.
\]
We obtain the following diagram:
\[
\begin{tikzcd}
0\arrow[r] & T_X\otimes \cO_{C_1}(-1)\oplus T_X\otimes \cO_{C_2}(-1)\arrow[r]\arrow[d] & T_{X}|_D\arrow[r]\arrow[d] &T_X|_{C_3}\arrow[r]\arrow[d] & 0\\
0\arrow[r] & Q_1\oplus Q_2 \arrow[r]\arrow[d] &Q  \arrow[r]\arrow[d] &Q|_{C_3}\arrow[r]\arrow[d] & 0\\
& 0 & 0& 0&
\end{tikzcd}
\]
where $Q_1$ is supported on $C_1$, $Q_2$ is supported on $C_2$. 
Then $Q_i\simeq \cO_{C_i}(k_i)$, where $k_i=-2$ or $k_i\geq 1$. While $Q|_{C_3}\in \langle\mathcal{P}, \cO_{C_3}(k_3)\rangle_{\text{ex}}$ where $\mathcal{P}$ is $0$-dimensional if nonzero, $k_3=-1$ or $k_3\geq 2$. There are only two cases in which we could obtain $\rho_{-tK}(Q)<-\frac{1}{2t}$: 
\begin{itemize}
    \item[(1)] $k_3=-1$, $k_1=k_2=-2$, and $\ell(\mathcal{P})=0$. This case cannot occur, since locally it is the same situation as Case (iii) in Proposition \ref{prop:desesdeg2}.
    \item[(2)]  $k_3=-1$, $k_1=k_2=-2$, and $\ell(\mathcal{P})=1$. In this case, $\mathcal{P}$ is supported at either $C_1\cap C_3$ or $C_2\cap C_3$. 
    Taking quotient, there exists a nonzero map $Q\to E_{-2, -1}$, contradicting Lemma \ref{lem:norcrossquot}.
\end{itemize}

Finally we consider Case (iv). Note that this can only happen when $n=6$, $C_1=e_1$, $C_2=l-e_1-e_2$ and $C_3=2l-e_1-e_3-e_4-e_5-e_6$. But when $n=6$, $-K\cdot c_1(Q)<3$, so this case cannot happen.  
\end{proof}

\begin{rem}
\label{rem:deg3caseiii}
    Note that, in Case (iii) of Proposition \ref{prop:desesdeg3}, the equality
$\rho_{-tK}(Q)=-\frac{1}{2t}$ can occur only when $\ell(\mathcal{P})=2$ and $\mathcal{P}\simeq \mathcal{O}_p\oplus\mathcal{O}_q$, where $p$ and $q$ lie on $C_1\cap C_3$ and $C_2\cap C_3$, respectively. If this is the case, we have 
\[
0\to \cO_{C_1}(-2)\oplus \cO_{C_2}(-2)\to Q\to \cO_p\oplus\cO_q\oplus\cO_{C_3}(-1)\to 0.
\]

Then $Q\simeq \cO_{C_1}(-1)\oplus \cO_{C_2}(-1)\oplus \cO_{C_3}(-1)$. Indeed, we can consider locally at one node. Consider the following diagram
\[
0\to \cO_{C_1}(-2)\to Q'\to \cO_p\oplus\cO_{C_3}(-1)\to 0.
\]
Projection to the pure quotient, we obtain the following diagram
\[
\begin{tikzcd}
& 0\ar{d} & 0\ar{d} & &\\
&\cO_{C_1}(-2)\ar{r}{\simeq}\ar{d}&\cO_{C_1}(-2)\ar{d}& &\\
0\ar{r}&\cO_{C_1}(-1)\ar{r}\ar{d}&Q'\ar{r}\ar{d}&\cO_{C_3}(-1)\ar{r}\ar{d}& 0\\
0\ar{r}&\cO_p\ar{r}\ar{d}&Q''\ar{r}\ar{d}&\cO_{C_3}(-1)\ar{r}&0\\
&0&0&&
\end{tikzcd}
\]

Applying the connecting homomorphism, we obtain the commutative diagram
\[
\begin{tikzcd}
\Hom(\cO_{C_3}(-1), \cO_{C_3}(-1))\ar{r}{\delta_1}\ar{d}&\Ext^1(\cO_{C_3}(-1), \cO_{C_1}(-1))\ar{d}{a}\\
\Hom(\cO_{C_3}(-1), \cO_{C_3}(-1))\ar{r}{\delta_2} & \Ext^1(\cO_{C_3}(-1), \cO_p).
\end{tikzcd}
\]
Since $\delta_2=0$, if $\delta_1\neq 0$, we have $a=0$.
Apply $\Hom(\cO_{C_3}(-1),\_)$ to the left most column, we obtain
\[
\begin{split}
&0\to \Hom(\cO_{C_3}(-1), \cO_p)\to \Ext^1(\cO_{C_3}(-1), \cO_{C_1}(-2))\to \\
&\Ext^1(\cO_{C_3}(-1), \cO_{C_1}(-1))\xrightarrow{a}\Ext^1(\cO_{C_3}(-1), \cO_p)\to... 
\end{split}
\]
Since all the vector spaces in the sequence are isomorphic to $\mathbb{C}$, we have $a$ is injective. This leads to a contradiction. Hence $\delta_1=0$, and $Q'\simeq \cO_{C_1}(-1)\oplus \cO_{C_3}(-1)$.

\end{rem}
\begin{prop}
\label{prop:destdeg4}
Let $X=X_n$ where $2\leq n\leq 6$ and let $0\to A\to T_X\to Q\to 0$ be a short exact sequence in $\Coh(X)\cap \mathcal{B}_{-K}$, where $Q$ is a pure $1$-dimensional sheaf.  Assume that $-K\cdot c_1(Q)=4$, then $\rho_{-tK}(Q)\geq -\frac{1}{2t}$.   
\end{prop}

\begin{proof}
Since $A\in \mathcal{B}_{-K}$, we have $-Kc_1(A)>0$. Then $K^2-4>0$, which implies that $n\leq 4$.

Using the same notation as in Proposition \ref{prop:desesdeg3}, and denoting $c_1(Q)=C_1+C_2+C_3+C_4$ where $C_i$'s are the divisor classes of the exceptional curves. 
If $c_1(Q)-D$ is a nonzero effective divisor, $Q$ is generically of rank $2$ in at least one irreducible component. Denoting that component by $C'$, then the surjection $T_X|_{C'}\to Q|_{C'}$ must be an isomorphism. Hence $\ch_2(Q|_{C'})=\ch_2(T_X|_{C'})=2$. 
If $Q\simeq T_X|_{C'}$, we are done, otherwise consider the short exact sequence 
\[
0\to M\to Q\to Q|_{C'}\to 0.
\]
Then $-K\cdot c_1(M)\leq 2$.
Since $Q$ is a quotient of $T_X$, we obtain $T_X(-C')\twoheadrightarrow M$, which implies a surjection $T_X\twoheadrightarrow M(C')$. Then Proposition \ref{prop:desesdeg2} implies that $\rho_{-tK}(M(C'))\geq -\frac{1}{2t}$.  Hence $\ch_2(M(C'))\geq -1$. Furthermore, $c_1(M)\cdot C'\leq 2$. Hence $\ch_2(M)\geq \ch_2(M(C'))-2$. Putting every bound together, we have $\ch_2(Q)\geq -1$, hence $\rho_{-tK}(Q)\geq -\frac{1}{2t}$.

Hence we only need to consider when $D\sim c_1(Q)$. We have the following possible connected cases: 

(i) $C_1= C_2= C_3= C_4$;

(ii) $C_1= C_2$, $C_3= C_4$ and $C_1\cdot C_3=1$; 

(iii) $C_2= C_3= C_4$, $C_2\neq C_1$ and $C_1\cdot C_2=1$; 

(iv) $C_3= C_4$, $C_i\neq C_j$ for $1\leq i\neq j\leq 3$, and $C_1\cdot C_3=1$, $C_2\cdot C_3=1$, $C_1\cdot C_2$=0;

(v) $C_3= C_4$, $C_i\neq C_j$ for $1\leq i\neq j\leq 3$, and $C_1\cdot C_2=1$, $C_1\cdot C_3=1$, $C_2\cdot C_3$=0.

(vi) $C_i\neq C_j$ for all $i\neq j$, $C_i\cdot C_{i+1}=1$ for $1\leq i\leq 3$, all other intersections are zero. 

(vii) $C_i\neq C_j$ for all $i\neq j$, $C_1\cdot C_i=1$ for $2\leq i\leq 4$. $C_i\cdot C_j=0$ for $2\leq i\neq j\leq 4$.

Case (i) follows from Lemma \ref{lem:destabnonred}. We only need to consider the rest of the cases. Let $D_0$ be the reduced support of $D$.

\textbf{Case (ii):}
Computing the intersection numbers, we may assume $D_0\sim C_1+C_3$. Similar to the argument in Proposition \ref{prop:desesdeg3}, we have the following diagram:
\begin{equation}
\label{diag:deg4supp}
\begin{tikzcd}
&&&0\arrow[d]&\\
&0\arrow[d]&&\cO_{D_1}(-C_1-C_3)\arrow[d]&\\
0\arrow[r]&\cO_X(-2C_1-2C_3)\arrow[r]\arrow[d]&\cO_X\arrow[r]\arrow[d, "\simeq"]&\cO_D \arrow[r]\arrow[d]&0\\
0\arrow[r]&\cO_X(-C_1-C_3)\arrow[r]\arrow[d]&\cO_X\arrow[r]&\cO_{D_0}\arrow[r]\arrow[d]&0\\
&\cO_{D_1}(-C_1-C_3)\arrow[d]&&0&\\
&0&&&
\end{tikzcd}
\end{equation}

Consider the following diagram 
\begin{equation}
\label{eq:deg4}
\begin{tikzcd}
0\arrow[r] & T_X\otimes \cO_{D_1}(-C_1-C_3)\arrow[r]\arrow[d] & T_{X}|_D\arrow[r]\arrow[d] &T_X|_{D_0}\arrow[r]\arrow[d] & 0\\
0\arrow[r] & Q_1 \arrow[r]\arrow[d] &Q  \arrow[r]\arrow[d] &Q|_{D_0}\arrow[r]\arrow[d] & 0\\
& 0 & 0&0 &
\end{tikzcd}
\end{equation}
Let $Q'$ be the torsion free quotient of $Q|_{D_0}$. Since $D_0\in |C_1+C_3|$,
Proposition \ref{prop:desesdeg2} implies that $\rho_{-tK}(Q')\geq -\frac{1}{2t}$, hence $\rho_{-tK}(Q|_{D_0})\geq -\frac{1}{2t}$. 
Furthermore, the same argument as in  Proposition \ref{prop:desesdeg2} case (iii) implies that $\rho_{-tK}(Q_1)\geq -\frac{1}{2t}$. Combining the slopes of $Q_1$ and $Q|_{D_0}$, we obtain $\rho_{-tK}(Q)\geq -\frac{1}{2t}$.

\textbf{Case (iii):} In this case $D_0= C_1+C_2$, i.e. $D_0$ has two irreducible components, each with multiciplity $1$. Consider the same diagram as \eqref{diag:deg4supp}, where we replace the left column by 
\[
0\to \cO_X(-C_1-3C_2)\to \cO_X(-C_1-C_2)\to \cO_{D_2}(-C_1-C_2)\to 0,
\]
and $D_1$ by $D_2$. Then $D_2\in |2C_2|$. Also consider the same diagram as \eqref{eq:deg4}, where we replace the left column by $T_X\otimes \cO_{D_2}(-C_1-C_2)\to Q_2\to 0$. 
In this case $D_0\in|C_1+ C_2|$. Replacing $Q|_{D_0}$ by its pure quotient if necessary,  Proposition \ref{prop:desesdeg2} implies that $\rho_{-tK}(Q|_{D_0})\geq -\frac{1}{2t}$. 
Note that $\cO_{D_2}$ fits into 
\begin{equation}
\label{eq:deg4(iii)}
0\to \cO_{C_2}(-C_2)\to \cO_{D_2}\to \cO_{C_2}\to 0.
\end{equation}

Tensoring $T_X\otimes \cO_{D_2}(-C_1-C_2)$ with the short exact sequence \eqref{eq:deg4(iii)}, we have $\rho_{-tK}(Q_2)\geq -\frac{1}{2t}$.

\textbf{Cases (iv):} If $D$ is nonreduced, we may assume $D_0\sim C_1+C_2+C_3$;
we replace the first row of diagram \eqref{diag:deg4supp} by 
\[
0\to \cO_X(-C_1-C_2-2C_3)\to \cO_X(-C_1-C_2-C_3)\to \cO_{D_3}(-C_1-C_2-C_3)\to 0,
\]
and $D_1$ by $D_3$. Also replace the first row of diagram \eqref{eq:deg4} by $T_X\otimes \cO_{D_3}(-C_1-C_2-C_3)\to Q_3$. Then $D_3\simeq C_3$ and $D_0\in |C_1+C_2+C_3|$. By Proposition \ref{prop:desesdeg3} case (iii), $\rho_{-tK}(Q|_{D_0})\geq -\frac{1}{2t}$. By the intersection numbers, $\cO_{C_3}(-C_1-C_2-C_3)\simeq \cO_{C_3}(-1)$. Then $\ch_2(Q_3)\geq -\frac{3}{2}$. If $\rho_{-tK}(Q|_{D_0})>-\frac{1}{2t}$, i.e. $\ch_2(Q|_{D_0})>-\frac{3}{2}$, we have $\chi(Q|_{D_0})>0$. This implies that $\chi(Q|_{D_0})\geq 1$, hence $\ch_2(Q|_{D_0})\geq -\frac{1}{2}$. Hence $\ch_2(Q)\geq -2$ and $\rho_{-tK}(Q)\geq -\frac{1}{2t}$. By Remark \ref{rem:deg3caseiii}, the equality $\rho_{-tK}(Q|_{D_0})=-\frac{1}{2t}$ can occur only if $Q|_{D_0}\simeq \cO_{C_1}(-1)\oplus\cO_{C_2}(-1)\oplus\cO_{C_3}(-1)$. The inequality $\rho_{-tK}(Q)<-\frac{1}{2t}$ can hold only when $Q_3\simeq \cO_{C_3}(-2)$. Since $\Ext^1(\cO_{C_3}(-1), \cO_{C_3}(-2))=0$, the proof of Lemma \ref{lem:norcrossquot} implies that such a map $T_X\to Q$ cannot exist.

Suppose $D_0=D\sim C_1+C_2+2C_3$. Consider the short exact sequences 
\[
\begin{split}
&0\to \cO_X(C_1+C_2+C_3)\to \cO_X(C_1+C_2+2C_3)\to \cO_{C_3}\to 0\\
&0\to \cO_X(C_1+C_2)\to \cO_X(C_1+C_2+C_3)\to \cO_{C_3}(1)\to 0\\
&0\to \cO_X(C_1)\to \cO_X(C_1+C_2)\to \cO_{C_2}(-1)\to 0\\
&0\to \cO_X\to \cO_X(C_1)\to \cO_{C_1}(-1)\to 0.
\end{split}
\]
Hence $|D|$ is three-dimensional. Using a similar argument as in Proposition \ref{prop:desesdeg3}, the generic element in the linear system is isomorphic to $\mathbb{P}^1$ with self-intersection $2$, the singular elements are reducible. The only reducible case which cannot be covered by the arguments of Proposition \ref{prop:desesdeg2} and \ref{prop:desesdeg3} is when $D=D_1\cup D_2\cup D_3\cup D_4$, where $D_i$'s are exceptional curves. Since $D^2=2$, $D_i^2=-1$ and $D_i\cdot D_j\leq 1$ if $i\neq j$, there are two possibilities, (a) $D_i\cdot D_{i+1}=1$, all other intersections are zero, (b) $D=D_1\cdot D_i=1$ for $2\leq i\leq 4$, all other intersections are zero.  
For (a), we have the following short exact sequence
\[
0\to \cO_{D_1}(-1)\to \cO_D\to \cO_{D'}\to 0,
\]
where $D'=D_2\cup D_3\cup D_4$ is as in Proposition \ref{prop:desesdeg3} Case (iii). Then the same argument as in the subcase $D_0\sim C_1+C_2+C_3$ implies that $\rho_{-tK}(Q)\geq -\frac{1}{2t}$. 
For (b), we have the short exact sequence 
\[
0\to \cO_{D_2}(-1)\to \cO_D\to \cO_D'\to 0
\]
where $D'=D_1\cup D_3\cup D_4$. This is again the same as in Proposition \ref{prop:desesdeg3} Case (iii). Hence the same argument implies that $\rho_{-tK}(Q)\geq -\frac{1}{2t}$.

\textbf{Case (v):} In this case, we may assume $D_0\sim C_1+C_2+C_3$. The argument is the same as the case (iv), with the only difference being $\cO_{C_3}(-C_1-C_2-C_3)\simeq \cO_{C_3}$, which will also give the desired bound for Bridgeland slope. 

\textbf{Case (vi) and Case (vii):} The argument is completely analogous to the argument of Case (vi). 
\end{proof}

Now we are ready to prove Theorem \ref{thm:mainthm}.

\begin{cor}
\label{cor:X2stable}
Let $X=X_2$.  Then $T_X$ is $\sigma_{-tK}$-stable for all $t>\sqrt{\frac{2}{7}}$, and $T_X$ is unstable when $t<\sqrt{\frac{2}{7}}$.
\end{cor}

\begin{proof}

Let 
\begin{equation}
\label{eq:X2}
0\to A\to T_X\to Q\to 0
\end{equation}
be a destabilizing short exact sequence in $\Coh(X)\cap \mathcal{B}_{-K}$ at the outermost wall $t$.
We first consider possibilities where $t>\sqrt\frac 27$.  Let us write $\beta$ to denote the morphism of sheaves $A \to T_X$ in the long exact sequence as usual.  Based on the  output of our MATLAB script \url{wall_est.m}, there are three cases we need to consider: 
\[
(\ch_0(\image \beta), -Kc_1(\image \beta), -Kc_1(Q^{-1})) = (1, 3, 0), (2, 5, 0), (2, 4, 0).
\]
By Lemma \ref{lem:AG54-112-1}, in each  of these three cases, we have  $Q^{-1}=0$ meaning $A$ is a subsheaf of $T_X$.  Hence the three cases are: 
\begin{itemize}
    \item[(i)] $\ch_0(\beta)=2, Q^{-1}=0$, and $-K\cdot c_1(Q^0)=3$;
    \item[(ii)] $\ch_0(\beta)=2, Q^{-1}=0$, and $-K\cdot c_1(Q^0)=2$;
    \item[(iii)] $\ch_0(\image \beta)=1, Q^{-1}=0$, and $-K\cdot c_1(Q^0)=4$.
\end{itemize}

Since the wall determined by the short exact sequence \eqref{eq:X2} is the outermost wall, we may assume that $Q$ does not have a zero-dimensional torsion subsheaf. Otherwise, we could replace $Q$ by its quotient by the zero-dimensional torsion subsheaf and obtain a higher wall. 
 Then cases (i) and (ii) follow from Propositions \ref{prop:desesdeg3} and \ref{prop:desesdeg2} respectively. 

Now we consider case (iii).  Since $-K\cdot c_1(Q)=4$, we have $-K\cdot c_1(Q^*)\geq -4$.  
Consider the short exact sequence obtained by dualizing:
\[
0\to Q^*\to \Omega_X\to M\to 0
\]
If $-K\cdot c_1(Q^*)> -4$, we have $\mu_{-K}(Q^*)>\mu_{-K}(\Omega_X)$, which contradicts that $\Omega_X$ is $\mu_{-K}$-stable. Hence $-K\cdot c_1(Q^*)=-4$ and $Q$ cannot have 1-dimensional torsion, and so it is torsion-free. Hence $Q^*$ is locally free. Note that $T_X$ cannot be strictly Bridgeland semistable for all $t>0$ by Proposition \ref{pro:Establehasoutermostwall}. Then
Lemma \ref{lem:AG54-105-1}(i) implies that $A$ is locally free. Hence $M$ is a subsheaf of $A^*$ which must also be torsion-free. 
Then Proposition \ref{prop:K25678} implies that $Q^*\simeq \cO_X(-2l+2e_i)$ and 
$c_1(A)=l+e_i-e_j$.
Computing the Bridgeland slopes, we have 
\[
\rho_{-tK}(A)=-\frac{1}{6t}-\frac{7}{6}t>\rho_{-tK}(T_X)=-\frac{3}{14t}-t.
\]
precisely when $t<\sqrt{\frac{2}{7}}$.

Proposition  \ref{pro:Establehasoutermostwall} together with  the MATLAB output now together imply that $T_X$ is $\sigma_{-tK}$-stable when $t> \sqrt\frac 27$, and case (iii) shows the outermost wall is at $t=\sqrt\frac 27$.
Note that there are other destabilizing short exact sequences on the wall $t=\sqrt\frac 27$. By the same argument as in Proposition \ref{prop:X1rank2wall}, the following short exact sequence in $\Coh(X)$ is also a short exact sequence in $\B_{-K}$:
\[
0\to M\to T_{X}\to \cO_{e}(-1)\to 0,
\]
where $e$ is an exceptional curve, i.e. $e=e_1, e_2$ or $l-e_1-e_2$.
Comparing the Bridgeland slope, we have $\rho_{-tK}(\cO_{e}(-1))<\rho_{-tK}(T_X)$ if and only if $t<\sqrt\frac 27$.
\end{proof}


\begin{cor}
\label{cor:X3stable}
Let $X=X_3$.  Then $T_X$ is $\sigma_{-tK}$-stable for all $t>0$.
\end{cor}
\begin{proof}
Equation \eqref{eq:slopeTX}  implies that $\rho_{-tK}(T_X)=-\frac{1}{2t}-t$.
Let 
\[
0\to A\to T_X\to Q\to 0
\]
be a destabilizing short exact sequence in $\mathcal{B}_{-K}$ at the outermost wall.  Let $\beta$ denote the map $A \to T_X$ in the long exact sequence in $\Coh (X)$.


If $Q^{-1}=0$ and $\ch_0(\text{im}\beta)=2$, the MATLAB output (see Table \ref{tab:wall-est-positive-n3}) gives that $-Kc_1(Q^0)$ could be $2, 3, 4$. Since the wall determined by this short exact sequence is the outermost wall, we may assume that $Q$ is pure of dimension one. Then these cases cannot happen by Propositions \ref{prop:desesdeg2}, \ref{prop:desesdeg3}, \ref{prop:destdeg4}, respectively.  Hence the remaining cases in Table \ref{tab:wall-est-positive-n3} are: 
\begin{itemize}
    \item If $\ch_0(\image \beta)=2$:
     \begin{itemize}
         \item[(i)] $-K\cdot c_1(Q^{-1})=-1$, $-K\cdot c_1(Q^0)=3$;
         \item[(ii)] $-K\cdot c_1(Q^{-1})=-2$, $-K\cdot c_1(Q^0)=3$;
         \item[(iii)] $-K\cdot c_1(Q^{-1})=-1$, $-K\cdot c_1(Q^0)=2$;
         \item[(iv)] $-K\cdot c_1(Q^{-1})=-2$, $-K\cdot c_1(Q^0)=2$;
         \item[(v)] $-K\cdot c_1(Q^{-1})=-3$, $-K\cdot c_1(Q^0)=2$. 
     \end{itemize}
     \item If $\ch_0(\image \beta)=1$:
    \begin{itemize}
        \item[(vi)] $Q^{-1}=0$, $-K\cdot c_1(Q^0)=4$.
    \end{itemize}
\end{itemize}

We first consider cases (i)-(v). Consider the short exact sequence in $\Coh(X)$
\[
0\to M\to T_X\to Q^0\to 0
\]
where $M:= \image \beta$.  Since $M$ is a quotient of $A$ in $\Coh(X)$, this is also a short exact sequence in $\mathcal{B}_{-K}$. If $Q^0$ is pure $1$-dimensional, then Proposition \ref{prop:desesdeg2}, \ref{prop:desesdeg3} and \ref{prop:destdeg4} imply that $\rho_{-tK}(Q^0)\geq -\frac{1}{2t}$. If $Q^0$ is not pure, we can apply the same propositions to its pure quotient, and still obtain $\rho_{-tK}(Q^0)\geq -\frac{1}{2t}$.
It is therefore enough to show $\rho_{-tK}(Q^{-1}[1])\geq -\frac{1}{2t}$ for the above cases.

Note that 
\[
\ch_2(Q^{-1})\leq \frac{(-K\cdot c_1(Q^{-1}))^2}{2K^2}.
\]
Indeed, denoting the $\mu$-HN-factors of $Q^{-1}$ by $H_i:=\HN_i(Q^{-1})$. Since $Q^{-1}\in \mathcal{F}_{-K}$, we have $\mu(H_i)\leq 0$ for all $i$. For each $H_i$, Bogomolov-Gieseker inequality and Hodge index theorem implies that 
\[
\ch_2(H_i)\leq \frac{\ch_1(H_i)^2}{2\ch_0(H_i)}\leq \frac{(-K\cdot\ch_1(H_i))^2}{2\ch_0(H_i)K^2}.
\]
Then we have 
\[
\begin{split}
\ch_2(Q^{-1})&=\sum_i\ch_2(H_i)\leq \sum_i\frac{(-K\cdot\ch_1(H_i))^2}{2\ch_0(H_i)K^2}\leq \frac{1}{2K^2}\sum_i(-K\cdot\ch_1(H_i))^2\\
&\leq \frac{1}{2K^2}(\sum_i-K\cdot\ch_1(H_i))^2=\frac{1}{2K^2}(-K\cdot c_1(Q^{-1}))^2.
\end{split}
\]
Hence 
\[
\begin{split}
    \rho_{-tK}(Q^{-1})&\geq \frac{\frac{1}{2K^2}(-K\cdot c_1(Q^{-1}))^2-\frac{1}{2}\ch_0(Q^{-1})K^2t^2}{t(-K\cdot c_1(Q^{-1}))}\\
    &\geq\frac{1}{2tK^2}(-K\cdot c_1(Q^{-1}))+\frac{\frac{1}{2}K^2t}{K\cdot c_1(Q^{-1})}.
    \end{split}
\]
When $n=3$, $K^2=6$. In each of  cases (i) through (v) we have  $-K\cdot c_1(Q^{-1})=-1, -2, -3$ and hence $\rho_{-tK}(Q^{-1})>0$.

Now consider case (vi). Using Lemma \ref{lem:AG54-105-1}(i), we know $A$ is a line bundle. Taking derived dual, we have the exact sequence
\[
0\to Q^\ast \to \Omega_X\to A^\ast\to...
\]
Since $\deg (Q)=4$, we have $\deg (Q^\ast)\geq -4$. By \cite[Tag 0AY4]{stacks-project}, $Q^\ast$ is a line bundle, and so  Proposition \ref{prop:K25678} applies and we must have $\deg(Q^\ast)=-4$ (hence $Q$ is torsion-free) and  $Q^\ast\simeq \cO_X(-2l+2e_i)$. Then $c_1(Q)=2l-2e_i$ and  $c_1(A)=c_1(T_X) - c_1(Q)=l+e_i-\sum_{j\neq i}e_j$.
Computing the Bridgeland slopes, we have 
\[
\rho_{-tK}(A)=-\frac{1}{2t}-\frac{3t}{2}<\rho_{-tK}(T_X)=-\frac{1}{2t}-t.
\]
\end{proof}

 \begin{cor}
 \label{cor:X4stable}
Let $X:=X_4$. Then $T_X$ is $\sigma_{-tK}$-stable for all $t>0$.
\end{cor}
\begin{proof}
Let $0\to A\to T_X\to Q\to 0$ be a destabilizing short exact sequence in $\mathcal{B}_{-K}$ at the outermost wall in the sense of \ref{para:outermostwall-existence}.  Let $\beta$ denote the map $A \to T_X$ in the long exact sequence in $\Coh (X)$.

From the MATLAB output (see Table \ref{tab:wall-est-positive-n4}), we need to consider the cases 
\begin{itemize}
\item If $\ch_0(\text{im}\beta)=2$:
\begin{itemize}
\item[(i)] $Q^{-1}=0$, $-K\cdot c_1(Q^0)=2$;

\item[(ii)] $Q^{-1}=0$, $-K\cdot c_1(Q^0)=3$; 

\item[(iii)] $-K\cdot c_1(Q^{-1})=-1$, $-K\cdot c_1(Q^0)=3$;

\item[(iv)] $-K\cdot c_1(Q^{-1})=-1$, $-K\cdot c_1(Q^0)=2$;

\item[(v)] $-K\cdot c_1(Q^{-1})=-2$, $-K\cdot c_1(Q^0)=2$;
\end{itemize}

\item If $\ch_0(\text{im}\beta)=1$:
\begin{itemize}
    \item [(vi)] $Q^{-1}=0$, $-K\cdot c_1(Q^0)=3$;

\item[(vii)] $Q^{-1}=0$, $-K\cdot c_1(Q^0)=4$;

\item[(viii)] $-K\cdot c_1(Q^{-1})=-1$, $-K\cdot c_1(Q^0)=3$.
\end{itemize}
\end{itemize}
For Cases (i) and (ii), such a destabilizing short exact sequence cannot exist by Propositions \ref{prop:desesdeg2}, \ref{prop:desesdeg3}, respectively. For case (vi), Lemma \ref{lem:AG54-105-1} implies that $A$ is locally free and $Q$ is a rank $1$ torsion free sheaf. Then $Q^*$ is a saturated invertible subsheaf of $\Omega_X$. The degree bound in Proposition \ref{prop:K25678} implies that case (vi) cannot occur. 
 
 We next consider cases (iii), (iv) and (v). 
The argument in the second half of the proof of Corollary \ref{cor:X3stable} also shows that
 \[
 \rho_{-tK}(Q^{-1})>0. 
 \]
Then since $\rho_{-tK}(Q^0)\geq -\frac{1}{2t}$ we have  $\rho_{-tK}(Q)>-\frac{1}{2t}>\rho_{-tK}(T_X)$.

Next we consider case (viii). Since $Q^{-1}\in\mathcal{F}_{-K}$ and $-K\cdot c_1(Q^{-1})=-1$, it has exactly one
$\mu_{-K}$-HN factor $R$ of negative degree, with $-K\cdot c_1(R)=-1$.
There may also be a factor $R_0$ of degree zero, but the Hodge index theorem and the BG inequality give $\ch_2(R_0)\le0$ if such a factor exists. For $R$, the Hodge index theorem gives $c_1(R)^2\le \tfrac15$;
by Lemma \ref{lem:nKDDsmod2}, however,  $c_1(R)^2$ must be odd and so $c_1(R)^2\le-1$.
By the BG inequality, we now have
\[
 \ch_2(R)\le-\frac{1}{2\ch_0(R)}<0
\]
implying  $\ch_2(R)\le-\tfrac12$  and hence $\ch_2(Q^{-1})\le-\frac12$.
On the other hand, from the Hodge index theorem we have $c_1(\image \beta)^2\le \tfrac45$; by Lemma \ref{lem:nKDDsmod2} again we know $c_1(\image \beta)^2$ must be even, and so it is at most 0.  Then the BG inequality gives $\ch_2(\image \beta) \leq 0$.
By Lemma \ref{lem:AG55-59-1}, $Q^{-1}$ must be zero, contradicting $-K\cdot c_1(Q^{-1})\neq 0$.

Finally we consider case (vii). Since $Q^{-1}=0$, it follows that  $A = \image \beta$ is a rank-one subsheaf of $T_X$ with $-K\cdot c_1(A)=1$.  Then by Lemma \ref{lem:AG54-105-1}, $A$ must be a line bundle and $Q$ must be a rank-one torsion-free sheaf.  Applying $\HHom (-, \OO_X)$ to the short exact sequence of sheaves $0 \to A \to T_X \to Q \to 0$ now gives the long exact sequence
\[
0 \to Q^\ast \to \Omega_X \overset{\alpha}{\to} A^\ast \to \EExt^1 (Q, \OO_X) \to 0
\]
from which we see $Q^\ast$ is a saturated invertible subsheaf of $\Omega_X$. Since $Q$ is torsion-free, $-Kc_1(Q^\ast)=-(-Kc_1(Q))=4$ and so, by  Proposition \ref{prop:K25678}(iii),  $Q^\ast$ is either of the form $\OO_X (-4l + 2 \sum_j e_j)$ or $\OO_X (-2l+2e_i)$ for some $i$, and so 
\[
c_1(A) = \begin{cases} -l + \sum_j e_j &\text{\quad if } Q^\ast = \OO_X (-4l + 2\sum_j e_j) \\
l + e_i - \sum_{j\neq i} e_j  &\text{\quad if } Q^\ast = \OO_X (-2l + 2e_i)
\end{cases}.
\]
In both cases we have $c_1(A)^2 = -3$ and, since $A$ is a line bundle, $\ch_2(A) = -\tfrac32$.  Now the wall equation $\rho_{-tK} (A) = \rho_{-tK} (T_X)$ can be rewritten as 
\[
-\tfrac32-\tfrac52 t^2 = -\tfrac{9}{10} - t^2
\]
which gives $t^2 = -\tfrac35$, a contradiction.  So this case cannot occur, either.
\end{proof}

\begin{cor}
\label{cor:X5stable}
Let $X:=X_5$. Then $T_X$ is $\sigma_{-tK}$-stable for all $t>0$.
\end{cor}
\begin{proof}
Let $0\to A\to T_X\to Q\to 0$ be a destabilizing  sequence in $\mathcal{B}_{-K}$ at the outermost wall in the sense of \ref{para:outermostwall-existence}. We can assume $Q$ is a pure $1$-dimensional sheaf. Let $\beta$ denote the map $A \to T_X$ in the long exact sequence in $\Coh (X)$.
From the MATLAB output (see Table \ref{tab:wall-est-positive-n5}), we need to consider the cases
\begin{itemize}
\item If $\ch_0(\text{im}\beta)=2$:
\begin{itemize}

\item[(i)] $Q^{-1}=0$, $-K\cdot c_1(Q^0)=2$; 

\item[(ii)] $Q^{-1}=0$, $-K\cdot c_1(Q^0)=3$; 

\item[(iii)] $-K\cdot c_1(Q^{-1})=-1$, $-K\cdot c_1(Q^0)=2$,
\end{itemize}
\item If $\ch_0(\text{im}\beta)=1$:
\begin{itemize}
\item[(iv)] $Q^{-1}=0$, $-K\cdot c_1(Q^0)=3$.
\end{itemize}
\end{itemize}

We first consider the cases when $Q^{-1}=0$. Since $\rho_{-tK}(T_X)=-\frac{3}{2t}-t$. If $A$ is of rank 2, such a destabilizing sequence cannot exist by Lemma \ref{lem:AG54-73-1}, Propositions \ref{prop:desesdeg2}, \ref{prop:desesdeg3} respectively.
When $A$ is of rank $1$, one has $-Kc_1(A)=1$. By Lemma \ref{lem:AG54-105-1}, we know that $Q$ is torsion-free. Applying $\mathcal{H}om(-, \cO_X)$, we get  
\[
0\to Q^*\to \Omega_X\to A^*\to...
\]
Since $Q$ is torsion-free, $Q^*$ is locally free. Furthermore $A^*$ being torsion-free implies that $Q^*$ is a saturated invertible subsheaf of $\Omega_X$ with $-K\cdot c_1(Q^*)\geq -3$. But Proposition \ref{prop:K25678} implies that $-K\cdot c_1(Q^*)\leq -4$,
hence such a destabilizing sequence cannot exist. 

We are only left with case (iii). The same estimate as in Corollary \ref{cor:X3stable} implies that 
$\rho_{-tK}(Q^{-1})>0$. Hence $\rho_{-tK}(Q)>-\frac{1}{2t}>\rho_{-tK}(T_X)$.
\end{proof}

\begin{cor}
\label{cor:X6stable}
Let $X=X_6$. Then $T_X$ is $\sigma_{-tK}$-stable for all $t>0$.
\end{cor}
\begin{proof}
Let  
\[
0\to A\to T_X\to Q\to 0
\]
be a destabilizing short exact sequence in $\mathcal{B}_{-K}$ at the outermost wall in the sense of \ref{para:outermostwall-existence}. Recall that Proposition \ref{pro:AG54-104-1} implies 
 $A$ is a rank-$2$ subsheaf of $T_X$ and $Q$ is a torsion sheaf with $-K\cdot c_1(Q)=2$. We can assume $Q$ is a pure $1$-dimensional sheaf. Since $\rho_{-tK}(T_X)=-\frac{5}{2t}-t$,
Proposition \ref{prop:desesdeg2} implies that such destabilizing short exact sequence cannot exist. 
\end{proof}

\newpage

\appendix

\section{Upper bounds of outermost wall from  Matlab scripts}
\label{app:wall-est-positive}

The Matlab script \url{wall_est.m} runs for del Pezzo surfaces $X=X_n$ where $3 \leq K^2 \leq 8$.  The tables below record all  triples 
\[
(r, d, q)  = (\ch_0 (\image \beta), -Kc_1 (\image \beta), -Kc_1(Q^{-1}) ) 
\]
under which an outermost wall $t$ in the sense of  \ref{para:outermostwall-existence} can potentially arise, and give an upper bound $b$ for $t^2$.  Triples that give rise to the same $b$ by the Matlab script are grouped together; the entries are also arranged in decreasing order of $b$ in each table.  Note that by Lemma \ref{lem:AG54-112-1}, if $-Kc_1(Q^{-1})=0$ then  $Q^{-1}=0$ in the destabilizing sequence at the outermost wall.

\noindent
\textbf{Caveat.}  If a triple $(r, d, q)$ is listed together with a bound $b$, it does not necessarily mean that there is an outermost wall $b$; it only means that the outermost wall $t$ could potentially occur with a destabilizing sequence with numerical data $(r, d, q)$, and that $t^2 \leq b$.

For $n=1$, the triple with an asterisk $(1, 4, 0)$ in the row with $b=\tfrac14$ requires separate treatment: When $(r, d, q)$ is equal to this triple, the coefficient of $t^2$ in the wall equation \eqref{eq:walleq_C123}, and the Matlab script displays a message saying that this case cannot be ruled out when $\ch_0(A)=1$.  This special case is dealt with in Lemma \ref{lem:Xeq1-specialcase}.

\begingroup
\small
\setlength{\tabcolsep}{5pt}
\renewcommand{\arraystretch}{1.2}

\begin{table}[htbp]
\centering
\caption{Bounds for $t^2$ when $n=1$, i.e.\  $K^2=8$.}
\label{tab:wall-est-positive-n1}
\begin{tabular}{rl}
\hline
$(r, d, q)$ & $b$ \\
\hline
 $(1,3, 0), (1, 4, -1), (1,4,-2), (1,4,-3), (2,6,0), (2,7,0)$ & $\frac12=0.5$ \\
 $(2,6,-5),(2,7,-6),(2,8,-7)$ & $\frac{5}{14}\approx 0.357$ \\
 $(2,6,-4),(2,7,-5),(2,8,-6)$ & $\frac13\approx 0.333$ \\
 $(2,6,-3),(2,7,-4),(2,8,-5)$ & $\frac{3}{10}=0.3$ \\
 $(1,4,0)^\ast, (2,6,-2), (2,7,-3), (2,8,-4)$ & $\frac14=0.25$\\
 $(2,5,-4)$ & $\frac{3}{14} \approx 0.214$ \\
 $(1,3,-2), (2,5,0), (2,5,-3), (2,6,-1), (2,7,-2), (2,8,-3)$ & $\frac16\approx 0.167$ \\
 $(2,4,-3)$ & $\frac17\approx 0.143$ \\
 $(2,4,0)$ & $\frac18=0.125$ \\
 $(2,5,-2)$ & $\frac{1}{10}=0.1$ \\
 $(2,4,-2)$ & $\frac{1}{12} \approx 0.083$\\
\hline
\end{tabular}
\end{table}

\begin{table}[htbp]
\centering
\caption{Bounds for $t^2$ when $n=2$, i.e.\  $K^2=7$.}
\label{tab:wall-est-positive-n2}
\begin{tabular}{rc}
\hline
$(r, d, q)$ & $b$ \\
\hline
 $(1,3, 0)$ & $\frac{16}{7} \approx 2.286$ \\ 
 $(2,5, 0)$& $\frac{11}{14} \approx 0.786$ \\
 $(2,4, 0)$ & $\frac{19}{42}\approx 0.452$ \\
 $(1,2,0), (1, 3, -1), (1, 3, -2), (2, 5, -1), (2, 5, -2), (2, 5, -3), (2, 5, -4), (2,6, 0)$ & $\frac27\approx 0.286$ \\
 $(2,4, -3)$ & $\frac{17}{84}\approx 0.202$ \\
 $(2, 4, -2)$ & $\frac{13}{70} \approx 0.186$ \\
 $(2,3, 0), (2,4, -1)$ & $\frac{9}{56} \approx 0.161$ \\
 $(2, 6, -5), (2,7, -6)$ & $\frac{5}{42} \approx 0.119$ \\
 $(2,2, 0), (2, 6, -4), (2, 7, -5)$ & $\frac{3}{35} \approx 0.086$ \\
 $(2, 3, -2), (2, 6, -3), (2, 7, -4)$ & $\frac{1}{28}\approx 0.036$\\
\hline
\end{tabular}
\end{table}

\begin{table}[htbp]
\centering
\caption{Bounds for $t^2$ when $n=3$, i.e.\  $K^2=6$.}
\label{tab:wall-est-positive-n3}
\begin{tabular}{rl}
\hline
$(r, d, q)$ & $b$ \\
\hline
 $(2,4, 0)$ & $\frac54=1.25$ \\ 
 $(1, 2, 0)$& $1$ \\ 
 $(2, 3, 0), (2, 4, -1)$ & $\frac12=0.5$ \\
 $(2, 4, -2)$ & $\frac38=0.375$ \\
 $(2, 4, -3)$ & $\frac{3}{10}=0.3$ \\
 $(2, 2, 0)$ & $\frac14=0.25$ \\
 $(2, 3, -1)$ & $\frac18=0.125$ \\
 $(2, 3, -2)$ & $\frac{1}{10}=0.1$ \\
\hline
\end{tabular}
\end{table}

\begin{table}[htbp]
\centering
\caption{Bounds for $t^2$ when $n=4$, i.e.\  $K^2=5$.}
\label{tab:wall-est-positive-n4}
\begin{tabular}{rl}
\hline
$(r, d, q)$ & $b$ \\
\hline
 $(1,2, 0)$ & $\frac{18}{5}=3.6$ \\ 
 $(2,3, 0)$& $\frac{27}{20} = 1.35$ \\
 $(2,2,0)$ & $\frac35=0.6$ \\
 $(2,3, -1)$ & $\frac{13}{30}\approx 0.433$ \\
 $(1, 1, 0), (1, 2, -1)$ & $\frac{4}{15} \approx 0.267$ \\
 $(2, 3, -2)$ & $\frac{9}{40}=0.225$ \\
 $(2, 2, -1)$ & $\frac{1}{10}=0.1$ \\
\hline
\end{tabular}
\end{table}

\begin{table}[htbp]
\centering
\caption{Bounds for $t^2$ when $n=5$, i.e.\  $K^2=4$.}
\label{tab:wall-est-positive-n5}
\begin{tabular}{rl}
\hline
$(r, d, q)$ & $b$ \\
\hline
 $(2,2, 0)$ & $\frac32=1.5$ \\ 
 $(1,1, 0)$& $1$ \\
 $(2,2,-1)$ & $\frac13\approx 0.333$ \\
 $(2,1, 0)$ & $\frac16 \approx 0.167$ \\
\hline
\end{tabular}
\end{table}

\begin{table}[htbp]
\centering
\caption{Bounds for $t^2$ when $n=6$, i.e.\  $K^2=3$.}
\label{tab:wall-est-positive-n6}
\begin{tabular}{rl}
\hline
$(r, d, q)$ & $b$ \\
\hline
 $(1,1, 0)$ & $4$ \\ 
 $(2,1, 0)$& $\frac34 = 0.75$ \\ 
\hline
\end{tabular}
\end{table}

\endgroup

\bibliography{refs}{}
\bibliographystyle{plain}

\end{document}